\documentclass[12pt]{article}
\usepackage{amsfonts, amssymb, amsmath, amsthm}

\makeatletter
\renewcommand{\subsection}{\@startsection{subsection}{2}{0pt}{-3.5ex plus -1ex minus -.2ex}{1ex}{\bf}}

\newtheorem{theorem}{Theorem}
\@addtoreset{theorem}{section}

\newtheorem{lemma}{Lemma}
\@addtoreset{lemma}{section}

\newtheorem{definition}{Definition}
\@addtoreset{definition}{section}

\newtheorem{remark}{Remark}
\@addtoreset{remark}{section}

\newtheorem{example}{Example}
\@addtoreset{example}{section}

\newtheorem{condition}{Condition}
\@addtoreset{condition}{section}

\@addtoreset{corollary}{section}

\@addtoreset{equation}{section}

\renewcommand{\Im}{{\rm Im\,}}
\renewcommand{\Re}{{\rm Re\,}}
\newcommand{\loc}{{\rm loc}}
\newcommand{\ind}{{\rm ind\,}}
\newcommand{\supp}{{\rm supp\,}}
\newcommand{\Dom}{{\rm Dom\,}}
\renewcommand{\ker}{{\rm ker\,}}
\renewcommand{\dim}{{\rm dim\,}}
\newcommand{\codim}{{\rm codim\,}}

\makeatother
\begin{document}

\begin{center}
\Large\bf Nonlocal problems for elliptic equations in dihedral
angles and the Green formula\footnotemark
\end{center}
\footnotetext{This work was partially supported by RFBR, grant
No~01-01-01030 and by grant No~E00-1-195 of Ministry of
Education.}
\begin{center}
Pavel Gurevich

\end{center}

\begin{abstract}
Nonlocal problems for higher-order elliptic operators in dihedral
and plane angles are considered. The Green formula is obtained,
which leads to  adjoint problems that take the form of nonlocal
transmission problems in dihedral and plane angles. This allows us
to establish necessary conditions of Fredholm solvability and
sufficient conditions of one-valued solvability for nonlocal
problems in dihedral angles.
\end{abstract}

\tableofcontents

\section*{Introduction}
In the theory of nonlocal elliptic boundary value prolems in
bounded domains, the most difficult case deals with the situation
when support of nonlocal terms intersects with boundary of a
domain (see~\cite{BitzDAN85}--\cite{SkDu91}). This leads to
apearance of degree singularities for solutions near some set.
Therefore it is natural to consider nonlocal elliptic problems in
weighted spaces (see~\cite{KondrTMMO67}--\cite{NP}). In order to
establish a priori estimates of solutions and construct a right
regularizer for nonlocal problems in bounded domains, one must
study nonlocal problems in dihedral and plane angles
(see~\cite{SkDu90, SkDu91}).

In paper~\cite{SkDu90}, A.L. Skubachevskii found sufficient
conditions of Fredholm solvability\footnote{ A closed operator
${\cal A}$ acting from a Hilbert space $H_1$ into a Hilbert space
$H_2$ is said to be {\it Fredholm} if its range ${\cal R}({\cal
A})$ is closed, dimension of its kernel $\dim\ker({\cal A})$ and
 codimension of its range $\codim{\cal R}({\cal A})$ are finite. The number
$\ind{\cal A}=\dim\ker({\cal A})-\codim{\cal R}({\cal A})$ is
called {\it index} of the Fredholm operator ${\cal A}.$ } for
auxiliary nonlocal problems with parameter $\theta$ in plane
angles and sufficient conditions of one--valued solvability for
model nonlocal problems in dihedral angles. His consideration was
based on a priori estimates of solutions and on using a right
regularizer which needed some additional conditions on a
corresponding ``local" model problem.

In the present work, we use another approach. Instead of
constructing a right regularizer, we obtain the Green formula and
study adjoint nonlocal problems. This leads to nonlocal
transmission problems in dihedral and plane angles. Similar
problems were studied in~\cite{Sh65, RSh} for the case of smooth
boundary of a domain, in~\cite{IM} for the one-dimensional case,
etc.

Our approach allows to establish 1) a necessary and sufficient
condition of Fredholm solvability for auxiliary nonlocal problems
with parameter $\theta$ in plane angles
(Theorem~\ref{thLThetaFred}); 2) necessary conditions of Fredholm
solvability and sufficient conditions of one--valued solvability
for model nonlocal problems in dihedral angles
(Theorems~\ref{thLSolv}, \ref{thLNecessCond}).

The paper is organized as follows.
In~\S\S\ref{sectBoundStatement}--\ref{sectBoundApr}, we consider
nonlocal boundary value problems in plane and dihedral angles. A
priori estimates in weighted spaces are established. For reader's
convinience, we formulate a number of results from the
paper~\cite{SkDu90}. In~\S\ref{sectGr}, we obtain the Green
formulas for nonlocal elliptic problems. The Green formulas
generate nonlocal transmission problems, which are formally
adjoint to nonlocal boundary value problems. Nonlocal transmission
problems are studied
in~\S\S\ref{sectTransStatement}--\ref{sectTransApr}. We prove the
results that are analogous to those
from~\S\S\ref{sectBoundStatement}--\ref{sectBoundApr}.
\S\ref{sectLMAdj} deals with operators that are adjoint to
operators of nonlocal boundary value problems. Connection between
adjoint operators and formally adjoint nonlocal transmission
problems is considered. The main results are collected
in~\S\ref{sectLSolv} where we study solvability of nonlocal
boundary value problems in plane and dihedral angles.
\S\ref{appendPoisson} illustrates the results obtained in this
work: we investigate the one--valued solvability of nonlocal
problems for the Poisson equation in dihedral angles. The paper
has two appendices. Appendix~\ref{appendL*Rn} deals with the
operator that is adjoint to the operator of elliptic problem in
${\mathbb R}^n$ with additional conditions on the hyperplane
$\{x_n=0\}.$ We prove a theorem concerning smoothness of solutions
for the corresponding problem. This result is used
in~\S\ref{sectLMAdj}. In Appendix~\ref{appendWeightSpace}, we
prove some auxiliary properties of weighted spaces that are needed
in the main part of the paper.

\section{Nonlocal elliptic boundary value problems. Reduction to problems with homogeneous nonlocal 
conditions}\label{sectBoundStatement}

\subsection{Nonlocal problems in dihedral angles.}
Introduce the sets
$$
   M=\{x=(y,\ z):\ y=0,\ z\in{\mathbb R}^{n-2}\},
$$
$$
  \Omega_j=\{x=(y,\ z):\ r>0,\ b_{j1}<\varphi<b_{j,R_j+1},\ z\in{\mathbb R}^{n-2}\},
$$
$$
  \Omega_{jt}=\{x=(y,\ z):\ r>0,\ b_{jt}<\varphi<b_{j,t+1},\ z\in{\mathbb R}^{n-2}\}\ (t=1,\ \dots,\ R_j),
$$
$$
  \Gamma_{jq}=\{x=(y,\ z):\ r>0,\ \varphi=b_{jq},\ z\in{\mathbb R}^{n-2}\}\ (q=1,\ \dots,\ R_j+1).
$$
Here  $x=(y,\ z)\in{\mathbb R}^n,\ y\in{\mathbb R}^2,\ z\in{\mathbb R}^{n-2};$ $r,\ \varphi$ are the polar coordinates of a point $y;$ 
$R_j\ge 1$ is an integers; $0<b_{j1}<\dots<b_{j,R_j+1}<2\pi;$ $j=1,\ \dots,\ N.$

Denote by ${\cal P}_j(D_y,\ D_z),$ $B_{j\sigma\mu}(D_y,\ D_z)$, and $B_{j\sigma\mu kqs}(D_y,\ D_z)$ 
homogeneous differential operators with constant complex coefficients of orders $2m,$  $m_{j\sigma\mu}\le 2m-1$ and
$m_{j\sigma\mu}\le 2m-1$ correspondingly
$(j,\ k=1,\ \dots,\ N;$ $\sigma=1,\ R_j+1;$ $\mu=1,\ \dots,\ m;$ $q=2,\ \dots,\ R_j$ $s=1,\ \dots,\ S_{j\sigma kq}).$ 

We shall assume that the following conditions hold (see~\cite[Chapter 2, \S\S 1.2, 1.4]{LM}).
\begin{condition}\label{condEllipP}
 For all $j=1,\ \dots,\ N$, the operators ${\cal P}_j(D_y,\ D_z)$ are properly elliptic.
\end{condition}
\begin{condition}\label{condComplB}
 For all  $j=1,\ \dots,\ N;\ \sigma=1,\ R_j+1$,
 the system $\{B_{j\sigma\mu}(D_y,\ D_z)\}_{\mu=1}^m$ is normal and covers the operator ${\cal P}_j(D_y,\ D_z)$ 
 on $\Gamma_{j\sigma}.$
\end{condition}

Consider the $N$ equations for functions $U_1,\ \dots,\ U_N$
\begin{equation}\label{eqP}
  {\cal P}_j(D_y,\ D_z)U_j=f_j(x) \quad (x\in\Omega_j)
\end{equation}
with the nonlocal conditions
\begin{equation}\label{eqB}
  \begin{array}{c}
  {\cal B}_{j\sigma\mu}(D_y,\ D_z)U=B_{j\sigma\mu}(D_y,\ D_z)U_j|_{\Gamma_{j\sigma}}+\\
  +\sum\limits_{k,q,s}
                   (B_{j\sigma\mu kqs}(D_y,\ D_z)U_k)({\cal G}_{j\sigma kqs}y,\ z)|_{\Gamma_{j\sigma}}
    =g_{j\sigma\mu}(x) \quad (x\in\Gamma_{j\sigma})
  \end{array}
\end{equation}
$$
 (j=1,\ \dots,\ N;\ \sigma=1,\ R_j+1;\ \mu=1,\ \dots,\ m).
$$
Here and below the summation in the formula for ${\cal B}_{j\sigma\mu}(D_y,\ D_z)$ is taken over
$k=1,\ \dots,\ N$; $q=2,\ \dots,\ R_k$; $s=1,\ \dots,\ S_{j\sigma kq}$;
$U=(U_1,\ \dots,\ U_N)$; $(B_{j\sigma\mu kqs}(D_y,\ D_z)U_k)({\cal G}_{j\sigma kqs}y,\ z)$ means that the expression
$(B_{j\sigma\mu kqs}(D_{y'},\ D_{z'})U_k)(x')$ is calculated for $x'=({\cal G}_{j\sigma kqs}y,\ z);$
${\cal G}_{j\sigma kqs}$ is the operator of rotation by the angle~$\varphi_{j\sigma kq}$ and expansion by~$\chi_{j\sigma kqs}$ 
times in the plane~$\{y\}$ such that $b_{k1}<b_{j\sigma}+\varphi_{j\sigma kq}=b_{kq}<b_{k,R_k+1},$ $0<\chi_{j\sigma kqs}.$

We introduce the space $H_a^l(\Omega)$ as a completion of the set $C_0^\infty(\bar\Omega\backslash M)$ in the norm
$$
 \|w\|_{H_a^l(\Omega)}=\left(
    \sum_{|\alpha|\le l}\int\limits_\Omega r^{2(a-l+|\alpha|)} |D_x^\alpha w(x)|^2 dx
                                       \right)^{1/2},
$$
where $\Omega=\{x=(y,\ z):\ r>0,\ 0<b_1<\varphi<b_2<2\pi,\ z\in{\mathbb R}^{n-2}\},$ 
$C_0^\infty(\bar\Omega\backslash M)$ is the set of infinitely differentiable functions
in $\bar\Omega$ with compact supports belonging to $\bar\Omega\backslash M;$
$a\in\mathbb R,$ $l\ge 0$ is an integer. Denote by $H_a^{l-1/2}(\Gamma)$ (for $l\ge1$) the space of traces on an
$(n-1)\mbox{-dimensional}$ half-plane $\Gamma\subset\bar\Omega$ with the norm
$$
 \|\psi\|_{H_a^{l-1/2}(\Gamma)}=\inf\|w\|_{H_a^l(\Omega)} \quad (w\in H_a^l(\Omega):\  w|_\Gamma = \psi).
$$

Introduce the spaces of vector--functions 
$$
 H_a^{l+2m,\,N}(\Omega)=\prod_{j=1}^N H_a^{l+2m}(\Omega_j),\ 
 H_a^{l,\,N}(\Omega,\ \Gamma)=\prod_{j=1}^N H_a^l(\Omega_j,\ \Gamma_j),
$$
$$
 H_a^l(\Omega_j,\ \Gamma_j)=
 H_a^l(\Omega_j)\times\prod_{\sigma=1,\,R_j+1}\prod_{\mu=1}^m H_a^{l+2m-m_{j\sigma\mu}-1/2}
   (\Gamma_{j\sigma}).
$$
We study solutions $U=(U_1,\ \dots,\ U_N)\in H_a^{l+2m,\,N}(\Omega)$ for problem~(\ref{eqP}), (\ref{eqB}) 
supposing that $f=\{f_j,\ g_{j\sigma\mu}\}\in  H_a^{l,\,N}(\Omega,\ \Gamma).$ Introduce the bounded operator
corresponding to problem~(\ref{eqP}), (\ref{eqB})
$$
 {\cal L}=\{{\cal P}_j(D_y,\ D_z),\ {\cal B}_{j\sigma\mu}(D_y,\ D_z)\}:
 H_a^{l+2m,\,N}(\Omega)\to  H_a^{l,\,N}(\Omega,\ \Gamma).
$$

\begin{lemma}\label{lHomog} 
 For any  $g_{j\sigma\mu}\in H_a^{l+2m-m_{j\sigma\mu}-1/2}(\Gamma_{j\sigma})$ 
 $(j=1,\ \dots,\ N;\ \sigma=1,\ R_j+1;\ \mu=1,\ \dots,\ m)$, there exists a vector--function $U\in  H_a^{l+2m,\,N}(\Omega)$
 such that
 $$
  {\cal B}_{j\sigma\mu}(D_y,\ D_z)U=g_{j\sigma\mu}(x)\ \quad (x\in\Gamma_{j\sigma}),
 $$
 $$
  \|U\|_{H_a^{l+2m,\,N}(\Omega)}\le 
  c\sum_{j,\,\sigma,\,\mu}\|g_{j\sigma\mu}\|_{H_a^{l+2m-m_{j\sigma\mu}-1/2}(\Gamma_{j\sigma})},
 $$
  where $c>0$ is independent of $g_{j\sigma\mu}.$
\end{lemma}
Lemma~\ref{lHomog} is proved in~\cite[\S1]{SkDu90}.

Let $W^l(Q)$ be a Sobolev space, where $Q\subset{\mathbb R}^n$ is an open domain with Lipschitz 
boundary. By $W^{l-1/2}(\Gamma)$ (for $l\ge1$) we denote the space of traces on an $(n-1)\mbox{-dimensional}$ smooth manifold
$\Gamma\subset\bar Q.$ Further we shall need interpolation inequalities for Sobolev and weighted spaces.

\begin{lemma}\label{lInterp1} 
 Let $Q$ be bounded; then for any $w\in W^l(Q)$ and $\lambda\in{\mathbb C}$, we have
 \begin{equation}\label{eqInterp1}
   |\lambda|^{l-s}\,\|w\|_{W^s(Q)}\le c_{ls}(\|w\|_{W^l(Q)}+|\lambda|^l\, \|w\|_{L_2(Q)}).
 \end{equation}
Here $0<s<l;\ c_{ls}>0$ is independent of $w,\ \lambda.$
\end{lemma}

\begin{lemma}\label{lInterp2} 
 Let $Q$ be bounded; then for any $w\in W^1(Q)$ and $\lambda\in{\mathbb C}$, we have
 \begin{equation}\label{eqInterp2}
   |\lambda|^{1/2}\,\|w|_{\partial Q}\|_{L_2(\partial Q)}\le c(\|w\|_{W^1(Q)}+|\lambda|\, \|w\|_{L_2(Q)}).
 \end{equation}
 Here $c>0$ is independent of $w,\ \lambda.$
\end{lemma}
Lemmas~\ref{lInterp1}, \ref{lInterp2} are proved in~\cite[Chapter 1, \S1]{AV}. Using lemma~\ref{lInterp1} and properties of weighted 
spaces, one can establish the following result (see~\cite[\S1]{SkMs86}).

\begin{lemma}\label{lInterpH}
 For any $w\in H_a^l(\Omega)$ and $\lambda\in\mathbb C$, we have
\begin{equation}\label{eqInterpH}
 |\lambda|^s\,\|w\|_{H_{a-s}^{l-s}(\Omega)}\le 
      c_{ls}(\|w\|_{H_a^l(\Omega)}+|\lambda|^l\, \|w\|_{H_{a-l}^0(\Omega)}).
\end{equation}
Here $0<s<l;\ c_{ls}>0$ is independent of $w,\ \lambda.$
\end{lemma}

\subsection{Nonlocal problems with parameter $\theta$ in plane angles.}
Now we consider the case of the space ${\mathbb R}^2$. Put
$K=\{y\in{\mathbb R}^2:\ r>0,\ 0<b_1<\varphi<b_2<2\pi\}.$ As above, we introduce the spaces $H_a^l(K)$ and $H_a^{l-1/2}(\gamma)$, 
where $\gamma\subset \bar K$ is a ray.

Let us also introduce the space $E_a^l(K)$ as a completion of the set
$C_0^\infty(\bar K\backslash \{0\})$ in the norm
$$
 \|w\|_{E_a^l(K)}=\left(
    \sum_{|\alpha|\le l}\int\limits_K r^{2a}(r^{2(|\alpha|-l)}+1) |D_y^\alpha w(y)|^2 dy
                                       \right)^{1/2}.
$$
By $E_a^{l-1/2}(\gamma)$ (for $l\ge1$) we denote the space of traces on a ray $\gamma\subset\bar K$ with the norm
$$
 \|\psi\|_{E_a^{l-1/2}(\gamma)}=\inf\|w\|_{E_a^l(K)} \quad (w\in E_a^l(K):\  w|_\gamma = \psi).
$$
One can find constructive definitions of the spaces $H_a^{l-1/2}(\Gamma)$ and $E_a^{l-1/2}(\gamma)$ in~\cite[\S1]{MP}.

Introduce the spaces of vector--functions
$$
 E_a^{l+2m,\,N}(K)=\prod_{j=1}^N E_a^{l+2m}(K_j),\ 
 E_a^{l,\,N}(K,\ \gamma)=\prod_{j=1}^N E_a^l(K_j,\ \gamma_j),
$$
$$
 E_a^l(K_j,\ \gamma_j)=
 E_a^l(K_j)\times\prod_{\sigma=1,\,R_j+1}\prod_{\mu=1}^m 
    E_a^{l+2m-m_{j\sigma\mu}-1/2}(\gamma_{j\sigma}),
$$
where $K_j=\{y:\ r>0,\ b_{j1}<\varphi<b_{j,R_j+1}\}$, $\gamma_{j\sigma}=\{y:\ r>0,\ \varphi=b_{j\sigma}\}$.

Consider the auxiliary problem for $u=(u_1,\ \dots,\ u_N)\in  E_a^{l+2m,\,N}(K)$
\begin{equation}\label{eqPTheta}
 {\cal P}_j(D_y,\ \theta)u_j=f_j(y) \quad (y\in K_j),
\end{equation}
\begin{equation}\label{eqBTheta}
  \begin{array}{c}
  {\cal B}_{j\sigma\mu}(D_y,\ \theta)u=B_{j\sigma\mu}(D_y,\ \theta)u_j|_{\gamma_{j\sigma}}+
  \sum\limits_{k,q,s}
                  (B_{j\sigma\mu kqs}(D_y,\ \theta)u_k)({\cal G}_{j\sigma kqs}y)|_{\gamma_{j\sigma}}=\\
   =g_{j\sigma\mu}(y) \quad (y\in\gamma_{j\sigma})
  \end{array}
\end{equation}
$$ 
 (j=1,\ \dots,\ N;\ \sigma=1,\ R_j+1;\ \mu=1,\ \dots,\ m),
$$
where $\theta$ is an arbitrary point on a unit sphere $S^{n-3}=\{z\in{\mathbb R}^{n-2}:\ |z|=1\},$
$f=\{f_j,\ g_{j\sigma\mu}\}\in  E_a^{l,\,N}(K,\ \gamma).$

Introduce the bounded operator corresponding to problem~(\ref{eqPTheta}), (\ref{eqBTheta}),
$$
 {\cal L}(\theta)=\{{\cal P}_j(D_y,\ \theta),\ {\cal B}_{j\sigma\mu}(D_y,\ \theta)\}:
 E_a^{l+2m,\,N}(K)\to  E_a^{l,\,N}(K,\ \gamma).
$$ 

\begin{lemma}\label{lHomogTheta}
 For any $g_{j\sigma\mu}\in E_a^{l+2m-m_{j\sigma\mu}-1/2}(\gamma_{j\sigma})$ 
 $(j=1,\ \dots,\ N;\ \sigma=1,\ R_j+1;\ \mu=1,\ \dots,\ m)$ and $\theta\in S^{n-3}$, there exists a vector--function
 $u\in  E_a^{l+2m,\,N}(K)$ such that
 $$
  {\cal B}_{j\sigma\mu}(D_y,\ \theta)u=g_{j\sigma\mu}(y)\ \quad (y\in\gamma_{j\sigma}),
 $$
 $$
  \|u\|_{E_a^{l+2m,\,N}(K)}\le 
  c\sum_{j,\,\sigma,\,\mu}\|g_{j\sigma\mu}\|_{E_a^{l+2m-m_{j\sigma\mu}-1/2}(\gamma_{j\sigma})},
 $$
 where $c>0$ is independent of $g_{j\sigma\mu},\ \theta.$
\end{lemma}
Lemma~\ref{lHomogTheta} is proved in~\cite[\S1]{SkDu90}.
\section{Solvability of nonlocal boundary value problems in plane angles}\label{sectBoundSolv}
We shall need the results of this section (obtained by A.L.~Skubachevskii in~\cite[\S2]{SkDu90}) in~\S\ref{sectBoundApr} 
for study a priori estimates of solutions to nonlocal boundary value problems in dihedral angles. 

\subsection{Reduction of nonlocal problems in plane angles to nonlocal problems on arcs.}
Consider the following nonlocal problem for $U=(U_1,\ \dots,\ U_N)\in H_a^{l+2m,\,N}(K)$
\begin{equation}\label{eqP0}
  {\cal P}_j(D_y,\ 0)U_j=f_j(x) \quad (y\in K_j),
\end{equation}
\begin{equation}\label{eqB0}
  \begin{array}{c}
  {\cal B}_{j\sigma\mu}(D_y,\ 0)U=B_{j\sigma\mu}(D_y,\ 0)U_j|_{\gamma_{j\sigma}}+\\
  +\sum\limits_{k,q,s}
                  (B_{j\sigma\mu kqs}(D_y,\ 0)U_k)({\cal G}_{j\sigma kqs}y)|_{\gamma_{j\sigma}}
   =g_{j\sigma\mu}(y) \quad (y\in\gamma_{j\sigma})
  \end{array}
\end{equation}
$$
 (j=1,\ \dots,\ N;\ \sigma=1,\ R_j+1;\ \mu=1,\ \dots,\ m),
$$
where $f=\{f_j,\ g_{j\sigma\mu}\}\in  H_a^{l,\,N}(K,\ \gamma).$

We write the operators ${\cal P}_j(D_y,\ 0),$ $B_{j\sigma\mu}(D_y,\ 0),$ $B_{j\sigma\mu kqs}(D_y,\ 0)$ in the polar coordinates: 
${\cal P}_j(D_y,\ 0)=r^{-2m}\tilde{\cal P}_j(\varphi,\ D_\varphi,\ rD_r),$
$B_{j\sigma\mu}(D_y,\ 0)=r^{-m_{j\sigma\mu}}\tilde B_{j\sigma\mu}(\varphi,\ D_\varphi,\ rD_r),$
$B_{j\sigma\mu kqs}(D_y,\ 0)=r^{-m_{j\sigma\mu}}\tilde B_{j\sigma\mu kqs}(\varphi,\ D_\varphi,\ rD_r),$ where
$D_\varphi=-i\frac{\displaystyle\partial}{\displaystyle\partial\varphi},\ 
D_r=-i\frac{\displaystyle\partial}{\displaystyle\partial r}.$

Put $\tau=\ln r$ and do the Fourier transform with respect to~$\tau$; then from~(\ref{eqP0}),
(\ref{eqB0}), we get
\begin{equation}\label{eqPLambda}
  \tilde{\cal P}_j(\varphi,\ D_\varphi,\ \lambda)\tilde U_j(\varphi,\ \lambda)=\tilde F_j(\varphi,\ \lambda) \quad 
   (b_{j1}<\varphi<b_{j,R_j+1}),
\end{equation}
\begin{equation}\label{eqBLambda}
 \begin{array}{c}
   \tilde{\cal B}_{j\sigma\mu}(\varphi,\ D_\varphi,\ \lambda)\tilde U(\varphi,\ \lambda)=
   \tilde B_{j\sigma\mu}(\varphi,\ D_\varphi,\ \lambda)\tilde U_j(\varphi,\ \lambda)|_{\varphi=b_{j\sigma}}+\\
           + \sum\limits_{k,q,s}
              e^{(i\lambda-m_{j\sigma\mu})\ln\chi_{j\sigma kqs}}{\tilde B}_{j\sigma\mu kqs}(\varphi,\ D_\varphi,\ \lambda)
              \tilde U_k(\varphi+\varphi_{j\sigma kq},\ \lambda)|_{\varphi=b_{j\sigma}}=\\
    =\tilde G_{j\sigma\mu}(\lambda),
  \end{array} 
\end{equation}
where $F_j(\varphi,\ \tau)=e^{2m\tau}f_j(\varphi,\ \tau),\ 
G_{j\sigma\mu}(\tau)=e^{m_{j\sigma\mu}\tau}g_{j\sigma\mu}(\tau);$
$
 \tilde U_j(\varphi,\ \lambda)=(2\pi)^{-1/2}\int\limits_{-\infty}^{\infty}U_j(\varphi,\ \tau)e^{-i\lambda\tau}d\tau.
$

This problem is a system of $N$ ordinary differential equations~(\ref{eqPLambda}) for 
functions $\tilde U_j\in W^{l+2m}(b_{j1},\ b_{j,R_j+1})$ with nonlocal conditions~(\ref{eqBLambda}) 
connecting values of $\tilde U_j$ and their derivatives at the point $\varphi=b_{j\sigma}$ with values of
$\tilde U_k$ and their derivatives at the points of the intervals $(b_{k1},\ b_{k,R_k+1}).$

\subsection{Solvability of nonlocal problems with parameter $\lambda$ on arcs.}
Let us consider the operator--valued function 
$$
 \begin{array}{c}
 \tilde{\cal L}(\lambda)=\{\tilde{\cal P}_j(\varphi,\ D_\varphi,\ \lambda),\ 
  \tilde{\cal B}_{j\sigma\mu}(\varphi,\ D_\varphi,\ \lambda)\}:\\
  W^{l+2m,\, N}(b_1,\ b_2)\to W^{l,\, N}[b_1,\ b_2]
 \end{array}
$$
corresponding to problem~(\ref{eqPLambda}), (\ref{eqBLambda}). Here
$$
 \begin{array}{c}
 W^{l+2m,\, N}(b_1,\ b_2)=\prod_{j=1}^N W^{l+2m}(b_{j1},\ b_{j,R_j+1}),\\
 W^{l,\, N}[b_1,\ b_2]=\prod_{j=1}^N W^{l}[b_{j1},\ b_{j,R_j+1}], 
 \end{array}
$$
$$
 W^{l}[b_{j1},\ b_{j,R_j+1}]=W^l(b_{j1},\ b_{j,R_j+1}) \times{\mathbb C}^m \times{\mathbb C}^m.
$$
Introduce the equivalent norms depending on the parameter $\lambda$ ($|\lambda|\ge1$) in the Hilbert spaces
$W^{l}(b_{j1},\ b_{j,R_j+1})$ and $W^{l}[b_{j1},\ b_{j,R_j+1}]$:
$$
 |||\tilde U_j|||_{W^{l}(b_{j1},\ b_{j,R_j+1})}=
 \bigl(\|\tilde U_j\|_{ W^{l}(b_{j1},\ b_{j,R_j+1})}^2+|\lambda|^{2l}\,\|\tilde U_j\|
   _{L_2(b_{j1},\ b_{j,R_j+1})}^2\bigr)^{1/2},
$$
$$
 \begin{array}{c}
 |||\{\tilde F_j,\ \tilde G_{j\sigma\mu}\}|||_{W^{l}[b_{j1},\ b_{j,R_j+1}]}=
 \Bigl(|||\tilde F_j|||_{ W^l(b_{j1},\ b_{j,R_j+1})}^2+\\
 +\sum\limits_{\sigma,\,\mu}(1+|\lambda|^{2(l+2m-m_{j\sigma\mu}-1/2)}) |\tilde G_{j\sigma\mu}|^2\Bigr)^{1/2},
 \end{array}
$$
where $\tilde U_j\in W^l(b_{j1},\ b_{j,R_j+1}),\ \{\tilde F_j,\ \tilde G_{j\sigma\mu}\}\in W^l[b_{j1},\ b_{j,R_j+1}].$
And therefore we have
$$
  |||\tilde U|||_{W^{l+2m,\,N}(b_1,\ b_2)}=
  \Bigl(\sum_j |||\tilde U_j|||^2_{W^{l+2m}(b_{j1},\ b_{j,R_j+1})}\Bigr)^{1/2},
$$
$$
  |||\tilde \Phi|||_{W^{l,\,N}[b_1,\ b_2]}=
  \Bigl(\sum_j |||\tilde \Phi_j|||^2_{W^l[b_{j1},\ b_{j,R_j+1}]}\Bigr)^{1/2},
$$
where $\tilde U=(\tilde U_1,\ \dots,\ \tilde U_N)\in W^{l+2m,\,N}(b_1,\ b_2),$
$\tilde\Phi=(\tilde\Phi_1,\ \dots,\ \tilde\Phi_N)\in W^{l,\,N}[b_1,\ b_2]$.

The next two statements are proved in~\cite[\S2]{SkDu90}.
\begin{lemma}\label{lSolvLambda1} 
 For all $\lambda\in{\mathbb C}$,
 the operator $\tilde{\cal L}(\lambda):W^{l+2m,\, N}(b_1,\ b_2)\to W^{l,\, N}[b_1,\ b_2]$ is Fredholm,
 $\ind\tilde{\cal L}(\lambda)=0;$ for any $h\in\mathbb R$, there exists a $q_0>0$ such that for 
 $\lambda\in J_{h,\,q_0}=\{\lambda\in{\mathbb C}:\ \Im\,\lambda=h,\ |\Re\,\lambda|\ge q_0\}$, 
 the operator $\tilde{\cal L}(\lambda)$ has the bounded inverse
 $\tilde{\cal L}^{-1}(\lambda):W^{l,\, N}[b_1,\ b_2]\to W^{l+2m,\, N}(b_1,\ b_2)$ and
 $$
  |||\tilde{\cal L}^{-1}(\lambda)\tilde \Phi|||_{W^{l+2m,\, N}(b_1,\ b_2)}\le c|||\tilde \Phi|||_{W^{l,\, N}[b_1,\ b_2]}
 $$
 for all $\tilde \Phi\in W^{l,\, N}[b_1,\ b_2],$ where $c>0$ is independent of $\lambda$ and $\tilde\Phi$. 
 The operator--valued function $\tilde{\cal L}^{-1}(\lambda):W^{l,\, N}[b_1,\ b_2]\to W^{l+2m,\, N}(b_1,\ b_2)$ 
 is finitely meromorphic.
\end{lemma}

\begin{lemma}\label{lSolvLambda2} 
 For any $0<\varepsilon<1/d$, there exists a $q>1$ such that the set
 $\{\lambda\in{\mathbb C}:\ |\Im\,\lambda|\le \varepsilon\ln\,|\Re\,\lambda|,\ |\Re\,\lambda|\ge q\}$ contains no poles of the
 operator--valued function $\tilde{\cal L}^{-1}(\lambda),$ where $d=\max\,|\ln\,\chi_{j\sigma kqs}|;$ for every pole
 $\lambda_0$ of the operator--valued function $\tilde{\cal L}^{-1}(\lambda)$, there exists a
 $\delta>0$ such that the set $\{\lambda\in{\mathbb C}:\ 0<|\Im\lambda - \Im\lambda_0|<\delta\}$ contains no poles of
 the operator--valued function $\tilde{\cal L}^{-1}(\lambda).$
\end{lemma}

\subsection{One--valued solvability of nonlocal problems in plane angles.}
The following theorem is obtained from Lemma~\ref{lSolvLambda1}.
\begin{theorem}\label{thSolvFlat} 
 Suppose the line $\Im\,\lambda=a+1-l-2m$ contains no poles of the operator--valued function $\tilde{\cal L}^{-1}(\lambda)$; 
 then nonlocal boundary value problem~(\ref{eqP0}), (\ref{eqB0}) has a unique solution $U\in H_a^{l+2m,\,N}(K)$ for every
 right-hand side $f\in H_a^{l,\,N}(K,\ \gamma)$ and
 $$
  \|U\|_{H_a^{l+2m,\,N}(K)}\le c\|f\|_{H_a^{l,\,N}(K,\ \gamma)},
 $$
 where $c>0$ does not depend on $f.$
\end{theorem}

One can find the proof of Theorem~\ref{thSolvFlat} in~\cite[\S2]{SkDu90}.
\section{A priori estimates of solutions for nonlocal boundary value problems}\label{sectBoundApr}

\subsection{A priori estimates in dihedral angles.}
Denote $d_1=\min\{1,\ \chi_{j\sigma kqs}\}/2,\ d_2=2 \max\{1,\ \chi_{j\sigma kqs}\},$
$\Omega_{j}^p=\Omega_j\cap\{r_1d_1^{6-p}<r<r_2d_2^{6-p},\ |z|<2^{-p-1}\},$ 
where $j=1,\ \dots,\ N;\ p=0,\ \dots,\ 6;\ 0<r_1<r_2.$

\begin{lemma}\label{lAprW} 
 Suppose $U_j\in W^{2m}(\Omega_j^0),$
 \begin{equation}\label{eqLUl}
  \begin{array}{c}
   {\cal P}_j(D_y,\ D_z)U_j\in W^l(\Omega_j^0),\ 
   {\cal B}_{j\sigma\mu}(D_y,\ D_z)U\in 
      W^{l+2m-m_{j\sigma\mu}-1/2}(\Gamma_{j\sigma}\cap\bar\Omega_j^0)
  \end{array}
 \end{equation}
 $$
 ( j=1,\ \dots,\ N;\ \sigma=1,\ R_j+1;\ \mu=1,\ \dots,\ m);
 $$
 then  $U\in \prod\limits_j W^{l+2m}(\Omega_j^3)$ and for $|\lambda|\ge 1$,
 \begin{equation}\label{eqAprW}
  \begin{array}{c}
   \sum\limits_j \|U_j\|_{W^{l+2m}(\Omega_j^6)}\le
        c \sum\limits_j \bigl\{\|{\cal P}_j(D_y,\ D_z)U_j\|_{W^l(\Omega_j^3)}+\\
   +\sum\limits_{\sigma,\,\mu}\|{\cal B}_{j\sigma\mu}(D_y,\ D_z)U\|
           _{W^{l+2m-m_{j\sigma\mu}-1/2}(\Gamma_{j\sigma}\cap\bar \Omega_j^3)}
   +|\lambda|^{-1}\|U_j\|_{W^{l+2m}(\Omega_j^3)}+\\
             +|\lambda|^{l+2m-1}\|U_j\|_{L_2(\Omega_j^3)}\bigr\},
  \end{array}
 \end{equation}
 where $c>0$ is independent of $\lambda$ and $U.$
\end{lemma}
\begin{proof}
Denote
\begin{equation}\label{eqEpsilon}
  \varepsilon=\min\{b_{j,q+1}-b_{jq}\}/4\ (j=1,\ \dots,\ N;\ q=1,\ \dots,\ R_j)
\end{equation}
and introduce the functions $\zeta_{jq}\in C^\infty({\mathbb R})$ such that
\begin{equation}\label{eqZeta}
 \zeta_{jq}(\varphi)=1\ \mbox{for } |b_{jq}-\varphi|<\varepsilon/2,\ 
     \zeta_{jq}(\varphi)=0\ \mbox{for } |b_{jq}-\varphi|>\varepsilon
\end{equation}
$$
 (j=1,\ \dots,\ N;\ q=1,\ \dots,\ R_j+1). 
$$

Put $\zeta_j(\varphi)=\zeta_{j1}(\varphi)+\zeta_{j,R_j+1}(\varphi).$
Since the functions $\zeta_j$ are the multiplicators in $W^l(\Omega_j^p),$ we have
$(1-\zeta_j)U_j\in W^{2m}(\Omega_j^0).$
Apply theorem~5.1 \cite[Chapter 2, \S 5.1]{LM} to the function $(1-\zeta_j)U_j$ and to the operator ${\cal P}_j(D_y,\ D_z)$;
then from~(\ref{eqLUl}) and Leibniz' formula, we get
\begin{equation}\label{eqAprW2}
  (1-\zeta_j)U_j\in W^{l+2m}(\Omega_j^1).
\end{equation}
Denote $V_{j\sigma\mu}=
  \sum\limits_{k,\,q,\,s}(B_{j\sigma\mu kqs}(D_y,\ D_z)((1-\zeta_k)U_k))({\cal G}_{j\sigma kqs}y,\ z).$
Clearly, we have
\begin{equation}\label{eqAprW3}
 V_{j\sigma\mu}|_{\Gamma_{j\sigma}\cap\bar\Omega_j^2}=
  \sum\limits_{k,\,q,\,s}(B_{j\sigma\mu kqs}(D_y,\ D_z)U_k)({\cal G}_{j\sigma kqs}y,\ z)|
         _{\Gamma_{j\sigma}\cap\bar\Omega_j^2}.
\end{equation}
From equality~(\ref{eqAprW3}) and relations~(\ref{eqLUl}), (\ref{eqAprW2}), it follows that
\begin{equation}\label{eqAprW4}
 \begin{array}{c}
  B_{j\sigma\mu}(D_y,\ D_z)U_j|_{\Gamma_{j\sigma}\cap\bar\Omega_j^2}=
  {\cal B}_{j\sigma\mu}(D_y,\ D_z)U - V_{j\sigma\mu}|_{\Gamma_{j\sigma}\cap\bar\Omega_j^2}\in\\
   \in W^{l+2m-m_{j\sigma\mu}-1/2}(\Gamma_{j\sigma}\cap\bar\Omega_j^2).
 \end{array}
\end{equation}
Again applying theorem~5.1 \cite[Chapter 2, \S 5.1]{LM} to the function $U_j$ and to the operator
$\{{\cal P}_j(D_y,\ D_z),\ B_{j\sigma\mu}(D_y,\ D_z)|_{\Gamma_{j\sigma}\cap\bar\Omega_j^2}\}$ 
$(\sigma=1,\ R_j+1;\ \mu=1,\ \dots,\ m)$ from~(\ref{eqLUl}),
(\ref{eqAprW4}), we obtain $U_j\in W^{l+2m}(\Omega_j^3).$

Now estimate~(\ref{eqAprW}) follows from lemma~3.1 \cite[\S3]{SkDu90}.
\end{proof}

Let $W_\loc^l(\bar\Omega_j\backslash M)$ be a set of functions belonging to the space
$W^l$ on any compactum in $\bar\Omega_j$ that does not intersect with $M.$

\begin{theorem}\label{thAprH}
 Let $U\in\prod\limits_j W_\loc^{2m}(\bar\Omega_j\backslash M)$ be a solution for nonlocal boundary value 
 problem~(\ref{eqP}), (\ref{eqB}) such that $U\in H_{a-l-2m}^{0,\,N}(\Omega)$ and 
 $f\in H_a^{l,\,N}(\Omega,\ \Gamma);$ then $U\in H_a^{l+2m,\,N}(\Omega)$ and
 \begin{equation}\label{eqAprH}
  \|U\|_{H_a^{l+2m,\,N}(\Omega)}\le c \bigl(\|f\|_{H_a^{l,\,N}(\Omega,\ \Gamma)}+
      \|U\|_{H_{a-l-2m}^{0,\,N}(\Omega)}\bigr),
 \end{equation}
 where $c>0$ is independent of $U.$
\end{theorem}
\begin{proof}
From Lemma~\ref{lAprW}, it follows that $U\in\prod\limits_j W_\loc^{l+2m}(\bar\Omega_j\backslash M).$
Now lemma~3.2 \cite[\S3]{SkDu90} implies that $U\in H_a^{l+2m,\,N}(\Omega)$ and a priori estimate~(\ref{eqAprH}) 
is valid.
\end{proof}

\subsection{A priori estimates in plane angles.}
Put $K_j^{ps}=K_j\cap\{r_1d_1^{6-p}\cdot 2^s<r<r_2d_2^{6-p}\cdot 2^s\},$ where $0<r_1<r_2;$ 
$s\ge 1;\   j=1,\ \dots,\ N;\ p=0,\ \dots,\ 6.$
\begin{lemma}\label{lAprWTheta}
 Suppose $s\ge 1,\ \theta\in S^{n-3}.$ Assume that $u_j\in W^{2m}(K_j^{0s}),$
 $$
  \begin{array}{c}
   {\cal P}_j(D_y,\ \theta)u_j\in W^l(K_j^{0s}),\ 
   {\cal B}_{j\sigma\mu}(D_y,\ \theta)u = 0 \quad (y\in \gamma_{j\sigma}\cap \bar K_j^{0s})
  \end{array}
 $$
 $$
 ( j=1,\ \dots,\ N;\ \sigma=1,\ R_j+1;\ \mu=1,\ \dots,\ m);
 $$
 then $u\in \prod\limits_j W^{l+2m}(K_j^{3s})$ and for $|\lambda|\ge 1$,
 \begin{equation}\label{eqAprWTheta}
  \begin{array}{c}
   \sum\limits_j  2^{sa}\|u_j\|_{W^{l+2m}(K_j^{6s})}\le
           c \sum\limits_j \bigl\{2^{sa}\|{\cal P}_j(D_y,\ \theta)u_j\|_{W^l(K_j^{3s})}+\\
   +|\lambda|^{-1}2^{sa}\|u_j\|_{W^{l+2m}(K_j^{3s})}+
             |\lambda|^{l+2m-1}2^{s(a-l-2m)}\|u_j\|_{L_2(K_j^{3s})} \bigr\},
  \end{array}
 \end{equation}
 where $c>0$ is independent of $u,$ $\theta$, $\lambda$, and $s.$
\end{lemma}
\begin{proof}
Repeating the proof of Lemma~\ref{lAprW} and substituting  $K_j^{ps}$ for $\Omega_j^p$ and $\theta$ for $D_z$, we obtain
$u\in \prod\limits_j W^{l+2m}(K_j^{3s}).$
Now a priori estimate~(\ref{eqAprWTheta}) follows from lemma~3.3 \cite[\S3]{SkDu90}.
\end{proof}

\begin{theorem}\label{thAprE1}
 Let $u\in\prod\limits_j W_\loc^{2m}(\bar K_j\backslash\{0\})$ be a solution for problem~(\ref{eqPTheta}),  
 (\ref{eqBTheta}) such that $u\in E_{a-l-2m}^{0,\,N}(K)$ and
  $f\in E_a^{l,\,N}(K,\ \gamma)$; then $u\in E_a^{l+2m,\,N}(K)$ and
  \begin{equation}\label{eqAprE1}
   \|u\|_{E_a^{l+2m,\,N}(K)}\le c \bigl(\|f\|_{E_a^{l,\,N}(K,\ \gamma)}+
       \|u\|_{E_{a-l-2m}^{0,\,N}(K)}\bigr),
  \end{equation}
  where $c>0$ is independent of $u,\ \theta\in S^{n-3}.$
\end{theorem}

\begin{proof}
1) By Lemma~\ref{lHomogTheta}, it suffices to consider the case $g_{j\sigma\mu}=0.$ Since
$f_j\in E_a^l(K_j)\subset W_\loc^l(\bar K_j\backslash\{0\}),$ as above, one can show that
$u\in\prod\limits_j W_\loc^{l+2m}(\bar K_j\backslash\{0\}).$ 
Put $r_1=d_1$, $r_2=d_2$ and denote $K_j^{ps}=K_j\cap\{d_1^{7-p}\cdot 2^s<r<d_2^{7-p}\cdot 2^s\},$ where
$s\ge 1;\   j=1,\ \dots,\ N;\ p=0,\ \dots,\ 6.$ Let us also denote $K_j^{60}=K_j\cap\{r<d_2\}$.
Introduce the functions $\psi\in C^\infty({\mathbb R}),$ $\psi(r)=1$ for
$r<d_2,$ $\psi(r)=0$ for $r>2d_2;$ $\hat\psi\in C^\infty({\mathbb R}),$ $\hat\psi(r)=1$ for
$r<2d_2^2,$ $\hat\psi(r)=0$ for $r>3d_2^2.$ 

Applying Theorem~\ref{thAprH} to the operator $\{{\cal P}_j(D_y,\ 0),\ {\cal B}_{j\sigma\mu}(D_y,\ 0)\}$ (for $n=2$), we get
\begin{equation}\label{eqAprE11}
 \begin{array}{c}
  \sum\limits_j \|u_j\|_{E_a^{l+2m}(K_j^{60})}\le  k_1\sum\limits_j \|\psi u_j\|_{H_a^{l+2m}(K_j)}\le \\
  \le k_2\sum\limits_j\bigl\{\|{\cal P}_j(D_y,\ 0)(\psi u_j)\|_{H_a^l(K_j)}+\\
  +\sum\limits_{\sigma,\,\mu}\|{\cal B}_{j\sigma\mu}(D_y,\ 0)(\psi u_j)\|_
             {H_a^{l+2m-m_{j\sigma\mu}-1/2}(\gamma_{j\sigma})}
   +\|\psi u_j\|_{H_{a-l-2m}^0(K_j)}\bigr\}.
 \end{array} 
\end{equation}
Let us estimate $\|{\cal P}_j(D_y,\ 0)(\psi u_j)\|_{H_a^l(K_j)}$. Using Leibniz' formula, the condition $\theta\in S^{n-3}$, and
limitations for supports of the functions $\psi,\ \hat\psi,$ we obtain
\begin{equation}\label{eqAprE12}
 \begin{array}{c}
   \|{\cal P}_j(D_y,\ 0)(\psi u_j)\|_{H_a^l(K_j)}\le \\
  \le k_3  (\|{\cal P}_j(D_y,\ \theta)(\psi u_j)\|_{H_a^l(K_j)}+
     \|\psi u_j\|_{H_{a-1}^{l+2m-1}(K_j)})\le\\
   \le k_4(\|{\cal P}_j(D_y,\ \theta)u_j\|_{E_a^l(K_j)}+ \|\hat\psi u_j\|_{H_{a-1}^{l+2m-1}(K_j)}).
 \end{array}
\end{equation}
Let us estimate $\|{\cal B}_{j\sigma\mu}(D_y,\ 0)(\psi u_j)\|_{H_a^{l+2m-m_{j\sigma\mu}-1/2}(\gamma_{j\sigma})}$. 
Using Leibniz' formula, the condition $\theta\in S^{n-3}$,
limitations for supports of the functions $\psi,\ \hat\psi$, and the condition $g_{j\sigma\mu}=0,$ we get
\begin{equation}\label{eqAprE13}
 \begin{array}{c}
   \|{\cal B}_{j\sigma\mu}(D_y,\ 0)(\psi u_j)\|_{H_a^{l+2m-m_{j\sigma\mu}-1/2}(\gamma_{j\sigma})}\le \\
   \le k_5( \|{\cal B}_{j\sigma\mu}(D_y,\ \theta)(\psi u_j)\|_{H_a^{l+2m-m_{j\sigma\mu}-1/2}(\gamma_{j\sigma})}
    +\|\hat \psi u_j\|_{H_{a-1}^{l+2m-1}(K_j)})\le\\
   \le k_6\bigl(\|\psi {\cal B}_{j\sigma\mu}(D_y,\ \theta)u_j\|_{H_a^{l+2m-m_{j\sigma\mu}-1/2}(\gamma_{j\sigma})}
    +\sum\limits_{k,\,q,\,s}\|(\psi(\chi_{j\sigma kqs}y)-\\
  -\psi(y))  
      (B_{j\sigma\mu kqs}(D_y,\ \theta)u_k)({\cal G}_{j\sigma kqs}y)|_{\gamma_{j\sigma}}\|_
     {H_a^{l+2m-m_{j\sigma\mu}-1/2}(\gamma_{j\sigma})}+\\
  + \|\hat \psi u_j\|_{H_{a-1}^{l+2m-1}(K_j)}\bigr)
   \le k_7\bigl(\sum\limits_k\|u_k\|_{W^{l+2m}(K_j\cap S_0)}+ \|\hat \psi u_j\|_{H_{a-1}^{l+2m-1}(K_j)}\bigr),
 \end{array}
\end{equation}
where $S_0=\{y\in{\mathbb R}^2:\ 1<r<2d_2/d_1\}.$

Inequalities~(\ref{eqAprE11})--(\ref{eqAprE13}), Lemma~\ref{lAprW}, and interpolation inequality~(\ref{eqInterpH}) yield
\begin{equation}\label{eqAprE14}
 \begin{array}{c}
   \sum\limits_j \|u_j\|_{E_a^{l+2m}(K_j^{60})}\le  k_8\sum\limits_j\bigl\{\|f_j\|_{E_a^l(K_j)}+\\
     + |\lambda|^{-1}\|u_j\|_{E_a^{l+2m}(K_j)}+|\lambda|^{l+2m-1}\|u_j\|_{E_{a-l-2m}^0(K_j)}\bigr\}
 \end{array}
\end{equation}

2) By virtue of Lemma~\ref{lAprWTheta}, for $s\ge 1$, we have
\begin{equation}\label{eqAprE15}
  \begin{array}{c}
   \sum\limits_j  \|u_j\|_{E_a^{l+2m}(K_j^{6s})}\le
           k_9 \sum\limits_j \bigl\{\|f_j\|_{E_a^l(K_j^{3s})}+\\
   +|\lambda|^{-1}\|u_j\|_{E_a^{l+2m}(K_j^{3s})}+
             |\lambda|^{l+2m-1}\|u_j\|_{E_{a-l-2m}^0(K_j^{3s})} \bigr\}.
  \end{array}
\end{equation}
Summing up (\ref{eqAprE14}), (\ref{eqAprE15}) for all $s\ge 1$ and taking a sufficiently large $|\lambda|$,
we obtain~(\ref{eqAprE1}).
\end{proof}

From Theorem~\ref{thSolvFlat} and Lemma~\ref{lAprWTheta}, one can also get the following result
(see theorem~3.1 \cite[\S3]{SkDu90}).
\begin{theorem}\label{thAprE2}
 Suppose the line $\Im\,\lambda=a+1-l-2m$ contains no poles of the operator--valued function $\tilde{\cal L}^{-1}(\lambda)$;
 then for all solutions $u\in E_a^{l+2m,\,N}(K)$ to nonlocal boundary value problem~(\ref{eqPTheta}), (\ref{eqBTheta}) 
 and all $\theta\in S^{n-3}$, we have
 \begin{equation}\label{eqAprE2}
   \|u\|_{E_a^{l+2m,\,N}(K)}\le c \bigl(\|f\|_{E_a^{l,\,N}(K,\ \gamma)}+
      \sum\limits_j \|u_j\|_{L_2(K_j\cap S)}\bigr),
 \end{equation}
 where $S=\{y\in{\mathbb R}^2:\ 0<R_1<r<R_2\}$, $c>0$ is independent of $\theta$ and $u.$

 If for any $\theta\in S^{n-3}$, estimate~(\ref{eqAprE2}) holds for all solutions to nonlocal boundary value 
 problem~(\ref{eqPTheta}), (\ref{eqBTheta}), then the line $\Im\,\lambda=a+1-l-2m$ contains no poles
 of the operator--valued function $\tilde{\cal L}^{-1}(\lambda).$
\end{theorem}

Theorem~\ref{thAprE2} implies that kernel of ${\cal L}(\theta)$ is of finite dimension and range of ${\cal L}(\theta)$
is closed. In order to prove that cokernel of ${\cal L}(\theta)$ is also of finite dimension, we shall obtain the Green formula
for nonlocal problems and study problems that are adjoint to nonlocal boundary value problems with respect to the
Green formula.
\section{The Green formula for nonlocal elliptic problems}\label{sectGr}
In this section, we obtain the Green formula, which connects nonlocal boundary value problems and nonlocal transmission
problems in dihedral angles, plane angles, and on arcs. Nonlocal transmission problems will be studied
in~\S\S\ref{sectTransStatement}--\ref{sectTransApr}.

\subsection{The Green formula in dihedral angles.}
Consider nonlocal boundary value problem~(\ref{eqP}), (\ref{eqB}).

Let $n_{kq}$ be the unit normal vector to $\Gamma_{kq}$ directed inside $\Omega_{kq}$
$(q=1,\ \dots,\ R_k),$ $n_{k,R_k+1}$ be the unit normal vector to $\Gamma_{k,R_k+1},$ 
directed inside $\Omega_{kR_k}.$

Denote by $C^\infty(\bar\Omega_{jt}\backslash M)$ ($C^\infty(\bar\Omega_j\backslash M),$ 
$C^\infty(\Gamma_{jq}\backslash M)$) the set of infinitely differentiable in
$\bar\Omega_{jt}\backslash M$ (in $\bar\Omega_j\backslash M,$
in $\Gamma_{jq}\backslash M$)
functions. We also denote by $C_0^\infty(\bar\Omega_{jt}\backslash M)$ 
($C_0^\infty(\bar\Omega_j\backslash M),$ $C_0^\infty(\Gamma_{jq}\backslash M)$) the set of infinitely differentiable in
$\bar\Omega_{jt}$ (in $\bar\Omega_j,$ in $\Gamma_{jq}$) functions with compact support from
$\bar\Omega_{jt}\backslash M$ (from $\bar\Omega_j\backslash M,$ from 
$\Gamma_{jq}\backslash M$) ($j=1,\ \dots\ N;$ $t=1,\ \dots,\ R_j;$ $q=1,\ \dots,\ R_j+1$).

For $U_{jt}\in C_0^\infty(\bar\Omega_{jt}\backslash M),\ V_{jt}\in C^\infty(\bar\Omega_{jt}\backslash M)$
(or $U_{jt}\in C^\infty(\bar\Omega_{jt}\backslash M),\ V_{jt}\in C_0^\infty(\bar\Omega_{jt}\backslash M)$),
put
$$
 (U_{jt},\ V_{jt})_{\Omega_{jt}}=\int\limits_{\Omega_{jt}}U_{jt}\cdot \bar V_{jt}\, dx\quad (j=1,\ \dots,\ N;\ t=1,\ \dots,\ R_j).
$$
For $U_{\Gamma_{jq}}\in C_0^\infty(\Gamma_{jq}),\ V_{\Gamma_{jq}}\in C^\infty(\Gamma_{jq})$
(or $U_{\Gamma_{jq}}\in C^\infty(\Gamma_{jq}),\ V_{\Gamma_{jq}}\in C_0^\infty(\Gamma_{jq})$),
put
$$
 (U_{\Gamma_{jq}},\ V_{\Gamma_{jq}})_{\Gamma_{jq}}=
          \int\limits_{\Gamma_{jq}} U_{\Gamma_{jq}}\cdot \bar V_{\Gamma_{jq}}\, d\Gamma\quad 
 (j=1,\ \dots,\ N;\ q=1,\ \dots,\ R_j+1).
$$

If we have functions $V_{jt}(x)$ defined in $\Omega_{jt}$, then denote by $V_j(x)$ the function given by
$V_j(x)\equiv V_{jt}(x)$ for $x\in\Omega_{jt}$.

For short, let us omit the arguments $(D_y,\ D_z)$ of differential operators. Denote by ${\cal Q}_j$ the operator that is
formally adjoint to ${\cal P}_j$. 

\begin{theorem}\label{thGrP}
 For the operators ${\cal P}_j,$ $B_{j\sigma\mu}$, and $B_{j\sigma\mu kqs}$ defined in \S1, there exist (not unique)

  1) a system $\{B'_{j\sigma\mu}\}_{\mu=1}^m$ of normal on $\Gamma_{j\sigma}$ operators of orders $2m-1-m'_{j\sigma\mu}$
 with constant coefficients such that the system
 $\{B_{j\sigma\mu},\ B'_{j\sigma\mu}\}_{\mu=1}^m$ is a Dirichlet one on $\Gamma_{j\sigma}$\footnote{
 See~\cite[Chapter 2, \S 2.2]{LM} for the definition of a Dirichlet system.} of order $2m$
 $(\sigma=1,\ R_j+1);$ 

  2) a Dirichlet system $\{B_{jq\mu},\ B'_{jq\mu}\}_{\mu=1}^m$ on $\Gamma_{jq}$ of order $2m$ such that
  the operators $B_{jq\mu}$ and $B'_{jq\mu}$ are of orders $2m-\mu$ and $m-\mu$ correspondingly $(q=2,\ \dots,\ R_j).$

 If the choice has been done, then there exist operators $C_{j\sigma\mu},\ F_{j\sigma\mu},\ T_{jq\nu}$, and
 $T_{jq\nu k\sigma s}$  $(j,\ k=1,\ \dots,\ N;$ $\sigma=1,\ R_j+1$ for the operators $C_{j\sigma\mu}$ and
 $F_{j\sigma\mu},$ $\sigma=1,\ R_k+1$ for the operators $T_{jq\nu k\sigma s};$  $\mu=1,\ \dots,\ m;$ 
 $q=2,\ \dots,\ R_j;$  $\nu=1,\ \dots,\ 2m;$ $s=1,\ \dots,\ S'_{jqk\sigma}=S_{k\sigma jq})$ with constant coefficients such
 that

  I) the operators $C_{j\sigma\mu},\ F_{j\sigma\mu},\ T_{jq\nu}$, and $T_{jq\nu k\sigma s}$ are of orders
 $m'_{j\sigma\mu},\ 2m-1-m_{j\sigma\mu},\ \nu-1$, and $\nu-1$ correspondingly;
 
  II) the system $\{C_{j\sigma\mu},\ F_{j\sigma\mu}\}_{\mu=1}^m$ is a Dirichlet one on $\Gamma_{j\sigma}$ of order
  $2m$ $(\sigma=1,\ R_j+1),$ 

  the system $\{C_{j\sigma\mu}\}_{\mu=1}^m$ covers the operator ${\cal Q}_j$ on $\Gamma_{j\sigma}$  $(\sigma=1,\ R_j+1),$ 

  the system $\{T_{jq\nu}\}_{\nu=1}^{2m}$ is a Dirichlet one on $\Gamma_{jq}$ of order $2m$ $(q=2,\ \dots,\ R_j);$ 

  III) for all $U_{j}\in C_0^\infty(\bar\Omega_{j}\backslash M),$ 
 $V_{jt}\in C^\infty(\bar\Omega_{jt}\backslash M)$ 
 (or $U_j\in C^\infty(\bar\Omega_j\backslash M),$ 
 $V_{jt}\in C_0^\infty(\bar\Omega_{jt}\backslash M)$), the following Green formula is valid:
 \begin{equation}\label{eqGrP}
 \begin{array}{c}
  \sum\limits_{j}\bigl\{\sum\limits_{t}({\cal P}_jU_j,\ V_{jt})_{\Omega_{jt}}+\sum\limits_{\sigma,\mu}
    ({\cal B}_{j\sigma\mu}U,\ F_{j\sigma\mu}V_j|_{\Gamma_{j\sigma}})_{\Gamma_{j\sigma}}+\\
 +\sum\limits_{q,\mu} (B_{jq\mu}U_j|_{\Gamma_{jq}},\ {\cal T}_{jq\mu}V)_
      {\Gamma_{jq}}\bigr\}=
\sum\limits_{j}\bigl\{\sum\limits_{t}(U_j,\ {\cal Q}_jV_{jt})_{\Omega_{jt}}+\\
 +\sum\limits_{\sigma,\mu}
    (B'_{j\sigma\mu}U_j|_{\Gamma_{j\sigma}},\ C_{j\sigma\mu}V_j|_{\Gamma_{j\sigma}})_{\Gamma_{j\sigma}}
 +\sum\limits_{q,\mu} (B'_{jq\mu}U_j|_{\Gamma_{jq}},\ {\cal T}_{jq,m+\mu}V)_
        {\Gamma_{jq}}\bigr\}.
 \end{array} 
 \end{equation}
 In the Green formulas (here and below), the summation is taken over $j=1,\ \dots,\ N$; $t=1,\ \dots,\ R_j$; $\sigma=1,\ R_j+1$;
 $q=2,\ \dots,\ R_j$; $\mu=1,\ \dots,\ m$;
 ${\cal B}_{j\sigma\mu}$ is given by~(\ref{eqB});
 $$
 \begin{array}{c}
 {\cal T}_{jq\nu}V=T_{jq\nu}V_{j,q-1}|_{\Gamma_{jq}}-T_{jq\nu}V_{jq}|_{\Gamma_{jq}}+
  \sum\limits_{k,\sigma,s}
    (T_{jq\nu k\sigma s}V_k)({\cal G}'_{jqk\sigma s}y,\ z)|_{\Gamma_{jq}}\\
  (\nu=1,\ \dots,\ 2m),
 \end{array}
 $$
 in the formula for ${\cal T}_{jq\nu}$ (here and below), the summation is taken 
 over $k=1,\ \dots,\ N$; $\sigma=1,\,R_k+1$; $s=1,\ \dots,\ S'_{jqk\sigma}=S_{k\sigma jq}$;
 ${\cal G}'_{jqk\sigma s}$ is the operator of rotation by the angle $\varphi'_{jqk\sigma}=-\varphi_{k\sigma jq}$ and expansion by 
 $\chi'_{jqk\sigma s}=1/\chi_{k\sigma j qs}$ times in the plane $\{y\}.$ 
\end{theorem}

\begin{proof}
For $j=1,\ \dots,\ N$, put
$
 B'_{j\sigma\mu}=
 \left(-i\frac{\displaystyle\partial}
           {\displaystyle\partial n_{j\sigma}}\right)^{2m-1-m'_{j\sigma\mu}},\ 
 B_{jq\mu}=\left(-i\frac{\displaystyle\partial}{\displaystyle\partial n_{jq}}\right)^{2m-\mu},\ 
 B'_{jq\mu}=\left(-i\frac{\displaystyle\partial}{\displaystyle\partial n_{jq}}\right)^{m-\mu}
$
$
 (\sigma=1,\ R_j+1;\ q=2,\ \dots, R_j;\ \mu=1,\ \dots,\ m),
$
where $m'_{j\sigma\mu}$ are chosen so that the numbers $m_{j\sigma\mu}$ and $2m-1-m'_{j\sigma\mu}$ run over the
set $0,\ 1,\ \dots,\ 2m-1,$ while $\mu$ changes from $1$ to $2m.$

By theorem~2.1 \cite[Chapter 2, \S 2.2]{LM}, there exist uniquely defined differential operators
$F_{j\sigma\mu},$ $F'_{j\sigma\mu},$  $F_{jq\mu}$, and $F'_{jq\mu}$ 
$(j=1,\ \dots,\ N;\ \sigma=1,\ R_j+1;\ q=2,\ \dots, R_j;\ \mu=1,\ \dots,\ m)$ of orders 
$2m-1-m_{j\sigma\mu},$ $m'_{j\sigma\mu},$ $\mu-1,$ and $m+\mu-1$ correspondingly
with constant coefficients such that

 the system $\{F_{j\sigma\mu},\ F'_{j\sigma\mu}\}_{\mu=1}^m$ is a Dirichlet one on $\Gamma_{j\sigma}$ of order $2m$ 
 $(\sigma=1,\ R_j+1),$ 

 the system  $\{F'_{j\sigma\mu}\}_{\mu=1}^m$ covers the operator $Q_j$ on $\Gamma_{j\sigma}$  $(\sigma=1,\ R_j+1),$ 

 the system $\{F_{jq\mu},\ F'_{jq\mu}\}_{\mu=1}^m$ is a Dirichlet one on $\Gamma_{jq}$ of order $2m$ $(q=2,\ \dots,\ R_j),$ 

 for any $U_j\in C_0^\infty(\bar\Omega_j\backslash M),$ 
 $V_{jt}\in C^\infty(\bar\Omega_{jt}\backslash M)$ (or $U_j\in C^\infty(\bar\Omega_j\backslash M),$ 
 $V_{jt}\in C_0^\infty(\bar\Omega_{jt}\backslash M)$), the following Green formulas are valid:
\begin{equation}\label{eqGrP1}
  \begin{array}{c}
 ({\cal P}_jU_{j},\ V_{j1})_{\Omega_{j1}}+\sum\limits_{\mu=1}^m
    (B_{j1\mu}U_{j}|_{\Gamma_{j1}},\ F_{j1\mu}V_{j1}|_{\Gamma_{j1}})_{\Gamma_{j1}}+\\
     +\sum\limits_{\mu=1}^m (B_{j2\mu}U_{j}|_{\Gamma_{j2}},\ F_{j2\mu}V_{j1}|_{\Gamma_{j2}})_
            {\Gamma_{j2}}
 = (U_{j},\ {\cal Q}_jV_{j1})_{\Omega_{j1}}+\\
  +\sum\limits_{\mu=1}^m
    (B'_{j1\mu}U_{j}|_{\Gamma_{j1}},\ F'_{j1\mu}V_{j1}|_{\Gamma_{j1}})_{\Gamma_{j1}}+
     \sum\limits_{\mu=1}^m (B'_{j2\mu}U_{j}|_{\Gamma_{j2}},\ F'_{j2\mu}V_{j1}|_{\Gamma_{j2}})_
           {\Gamma_{j2}},\\
\\
 ({\cal P}_jU_{j},\ V_{j2})_{\Omega_{j2}}-\sum\limits_{\mu=1}^m
    (B_{j2\mu}U_{j}|_{\Gamma_{j2}},\ F_{j2\mu}V_{j2}|_{\Gamma_{j2}})_{\Gamma_{j2}}+\\
     +\sum\limits_{\mu=1}^m (B_{j3\mu}U_{j}|_{\Gamma_{j3}},\ F_{j3\mu}V_{j2}|_{\Gamma_{j3}})_
              {\Gamma_{j3}}
 = (U_{j},\ {\cal Q}_jV_{j2})_{\Omega_{j2}}-\\
   -\sum\limits_{\mu=1}^m
    (B'_{j2\mu}U_{j}|_{\Gamma_{j2}},\ F'_{j2\mu}V_{j2}|_{\Gamma_{j2}})_{\Gamma_{j2}}+
     \sum\limits_{\mu=1}^m (B'_{j3\mu}U_{j}|_{\Gamma_{j3}},\ F'_{j3\mu}V_{j2}|_{\Gamma_{j3}})_
             {\Gamma_{j3}},\\
\\
    \cdots,\\
\\
 ({\cal P}_jU_{j},\ V_{jR_j})_{\Omega_{jR_j}}
 -\sum\limits_{\mu=1}^m
    (B_{jR_j\mu}U_{j}|_{\Gamma_{jR_j}},\ F_{jR_j\mu}V_{jR_j}|_{\Gamma_{jR_j}})_{\Gamma_{jR_j}}+\\
     +\sum\limits_{\mu=1}^m
     (B_{j,R_j+1,\mu}U_{j}|_{\Gamma_{j,R_j+1}},\ F_{j,R_j+1,\mu}V_{jR_j}|_{\Gamma_{j,R_j+1}})_
             {\Gamma_{j,R_j+1}}=\\
 = (U_{j},\ {\cal Q}_jV_{jR_j})_{\Omega_{jR_j}}
 -\sum\limits_{\mu=1}^m
    (B'_{jR_j\mu}U_{j}|_{\Gamma_{jR_j}},\ F'_{jR_j\mu}V_{jR_j}|_{\Gamma_{jR_j}})_
            {\Gamma_{jR_j}}+\\
     +\sum\limits_{\mu=1}^m
    (B'_{j,R_j+1,\mu}U_{j}|_{\Gamma_{j,R_j+1}},\ F'_{j,R_j+1,\mu}V_{jR_j}|_{\Gamma_{j,R_j+1}})_
           {\Gamma_{j,R_j+1}}.
 \end{array} 
\end{equation}
Adding equalities~(\ref{eqGrP1}) together, we get
\begin{equation}\label{eqGrP2}
 \begin{array}{c}
   \sum\limits_t({\cal P}_jU_j,\ V_{jt})_{\Omega_{jt}}+\sum\limits_{\sigma=1,\,R_j+1}\sum\limits_{\mu=1}^m
  (B_{j\sigma\mu}U_j|_{\Gamma_{j\sigma}},\ F_{j\sigma\mu}V_j|_{\Gamma_{j\sigma}})_
          {\Gamma_{j\sigma}}+\\
    +\sum\limits_{q=2}^{R_j}\sum\limits_{\mu=1}^m 
  (B_{jq\mu}U_j|_ {\Gamma_{jq}},\ F_{jq\mu}V_{j,q-1}|_ {\Gamma_{jq}}-F_{jq\mu}V_{jq}|_ {\Gamma_{jq}})_
          {\Gamma_{jq}}=\\
 = \sum\limits_t(U_j,\ {\cal Q}_jV_{jt})_{\Omega_{jt}}+\sum\limits_{\sigma=1,\,R_j+1}\sum\limits_{\mu=1}^m
    (B'_{j\sigma\mu}U_j|_{\Gamma_{j\sigma}},\ F'_{j\sigma\mu}V_j|_{\Gamma_{j\sigma}})_
          {\Gamma_{j\sigma}}+\\
    +\sum\limits_{q=2}^{R_j}\sum\limits_{\mu=1}^m 
  (B'_{jq\mu}U_j|_{\Gamma_{jq}},\ F'_{jq\mu}V_{j,q-1}|_{\Gamma_{jq}}-F'_{jq\mu}V_{jq}|_{\Gamma_{jq}})_
          {\Gamma_{jq}}.
 \end{array}
\end{equation}

Add $\sum\limits_{k=1}^N\sum\limits_{q=2}^{R_k}\sum\limits_{s=1}^{S_{j\sigma kq}}
(B_{j\sigma\mu kqs}U_k)({\cal G}_{j\sigma kqs}\cdot)$ to $B_{j\sigma\mu}U_j$ and subtract it
in~(\ref{eqGrP2}); then using change of variables  $x'=({\cal G}_{j\sigma kqs}y,\ z)$ in the integrals over $\Gamma_{j\sigma},$ 
we obtain
\begin{equation}\label{eqGrP3}
 \begin{array}{c}
  (B_{j\sigma\mu}U_j|_{\Gamma_{j\sigma}},\ F_{j\sigma\mu}V_j|_{\Gamma_{j\sigma}})_{\Gamma_{j\sigma}}=
   ({\cal B}_{j\sigma\mu}U,\ F_{j\sigma\mu}V_j|_{\Gamma_{j\sigma}})_{\Gamma_{j\sigma}}-\\
   -(\sum\limits_{k,q}\sum\limits_{s=1}^{S_{j\sigma kq}}
 (B_{j\sigma\mu kqs}U_k)({\cal G}_{j\sigma kqs}\cdot)|_{\Gamma_{j\sigma}},\ F_{j\sigma\mu}V_j|_{\Gamma_{j\sigma}})_
         {\Gamma_{j\sigma}}=\\
 =({\cal B}_{j\sigma\mu}U,\ F_{j\sigma\mu}V_j|_{\Gamma_{j\sigma}})_{\Gamma_{j\sigma}}+\\
  +\sum\limits_{k,q}\sum\limits_{s=1}^{S'_{kqj\sigma}}
     (-\frac{\displaystyle 1}{\displaystyle\chi_{j\sigma kqs}}
 B_{j\sigma\mu kqs}U_k|_{\Gamma_{kq}},\ (F_{j\sigma\mu}V_j)({\cal G}'_{kqj\sigma s}\cdot)|_{\Gamma_{kq}})_
        {\Gamma_{kq}}.
 \end{array} 
\end{equation}
Here $S'_{kqj\sigma}=S_{j\sigma kq};$ ${\cal G}'_{kqj\sigma s}$ is the operator of rotation by the angle
$\varphi'_{kqj\sigma}=-\varphi_{j\sigma k q}$ and expansion by 
$\chi'_{kqj\sigma s}=1/\chi_{j\sigma k qs}$ times in the plane $\{y\}.$ 

Clearly, we have
\begin{equation}\label{eqGrP4}
 -\frac{\displaystyle 1}{\displaystyle\chi_{j\sigma kqs}}B_{j\sigma\mu kqs}=
 \sum_{\alpha=1}^m \Lambda_{j\sigma\mu kqs\alpha}B_{kq\alpha}-
 \sum_{\alpha=1}^m \Lambda'_{j\sigma\mu kqs\alpha}B'_{kq\alpha}\footnotemark{}.
\end{equation} 
\footnotetext{We choose the sign ``minus'' in right hand side of relation~(\ref{eqGrP4}) just for convenience.}
Here  
$$
 \begin{array}{c}
  \Lambda_{j\sigma\mu kqs\alpha}=\sum\limits_{|\beta|+l=0}^{m_{j\sigma\mu}-(2m-\alpha)} a_{j\sigma\mu kqs\alpha}^{\beta l}
  D_z^\beta\left(\frac{\displaystyle\partial}{\displaystyle\partial y_{kq}}\right)^l,\\
  \Lambda'_{j\sigma\mu kqs\alpha}=\sum\limits_{|\beta|+l=0}^{m_{j\sigma\mu}-(m-\alpha)} {a'}_{j\sigma\mu kqs\alpha}^{\beta l}
  D_z^\beta\left(\frac{\displaystyle\partial}{\displaystyle\partial y_{kq}}\right)^l,
 \end{array}
$$
$a_{j\sigma\mu kqs\alpha}^{\beta l}$, ${a'}_{j\sigma\mu kqs\alpha}^{\beta l}\in{\mathbb C}$, $y_{kq}$ is the coordinate on the
half-axis $\Gamma_{kq}\cap\{z=0\}$. If $m_{j\sigma\mu}-(2m-\alpha)<0$ 
$(m_{j\sigma\mu}-(m-\alpha)<0),$ then we put
$\Lambda_{j\sigma\mu kqs\alpha}=0$ $(\Lambda'_{j\sigma\mu kqs\alpha}=0).$

Denote by $(\Lambda_{j\sigma\mu kqs\alpha})^*,$ $(\Lambda'_{j\sigma\mu kqs\alpha})^*$ the operators that are formally
adjoint to $\Lambda_{j\sigma\mu kqs\alpha},$ $\Lambda'_{j\sigma\mu kqs\alpha}$ correspondingly. 
Then~(\ref{eqGrP3}) and (\ref{eqGrP4}) imply
\begin{equation}\label{eqGrP5}
 \begin{array}{c}
  (B_{j\sigma\mu}U_j|_{\Gamma_{j\sigma}},\ F_{j\sigma\mu}V_j|_{\Gamma_{j\sigma}})_
      {\Gamma_{j\sigma}}=
  ({\cal B}_{j\sigma\mu}U,\ F_{j\sigma\mu}V_j|_{\Gamma_{j\sigma}})_{\Gamma_{j\sigma}}+\\
  +\sum\limits_{k,q}\sum\limits_{s=1}^{S'_{kqj\sigma}}\sum\limits_{\alpha=1}^m
     (B_{kq\alpha}U_k|_{\Gamma_{kq}},\ 
     (\Lambda_{j\sigma\mu kqs\alpha})^*[(F_{j\sigma\mu}V_j)({\cal G}'_{kqj\sigma s}\cdot)|_{\Gamma_{kq}}])_
          {\Gamma_{kq}}-\\
  -\sum\limits_{k,q}\sum\limits_{s=1}^{S'_{kqj\sigma}}\sum\limits_{\alpha=1}^m
     (B'_{kq\alpha}U_k|_{\Gamma_{kq}},\ 
     (\Lambda'_{j\sigma\mu kqs\alpha})^*[(F'_{j\sigma\mu}V_j)({\cal G}'_{kqj\sigma s}\cdot)|_{\Gamma_{kq}}])_
         {\Gamma_{kq}}.
 \end{array}
\end{equation}
Substituting~(\ref{eqGrP5}) into~(\ref{eqGrP2}), summing over~$j$, and grouping the summands containing $B_{jq\mu}U_j,$ 
we get
\begin{equation}\label{eqGrP6}
 \begin{array}{c}
  \sum\limits_j\bigl\{\sum\limits_t({\cal P}_jU_j,\ V_{jt})_{\Omega_{jt}}+
    \sum\limits_{\sigma=1,\,R_j+1}\sum\limits_{\mu=1}^m
    ({\cal B}_{j\sigma\mu}U,\ F_{j\sigma\mu}V_j|_{\Gamma_{j\sigma}})_{\Gamma_{j\sigma}}+\\
    +\sum\limits_{q=2}^{R_j}\sum\limits_{\mu=1}^m (B_{jq\mu}U_j|_{\Gamma_{jq}},\ 
    F_{jq\mu}V_{j,q-1}|_{\Gamma_{jq}}-F_{jq\mu}V_{jq}|_{\Gamma_{jq}}+\\
     +\sum\limits_k\sum\limits_{\sigma=1,\,R_k+1}\sum\limits_{s=1}^{S'_{jqk\sigma}}
     (\{\sum\limits_{\alpha=1}^m
     (\hat\Lambda_{k\sigma\alpha jqs\mu})^*F_{k\sigma\alpha}\}V_k)({\cal G}'_{jqk\sigma s}\cdot)|_{\Gamma_{jq}})_
          {\Gamma_{jq}}\bigr\}=\\
\\
 = \sum\limits_j\bigl\{\sum\limits_t (U_j,\ {\cal Q}_jV_{jt})_{\Omega_{jt}}+
    \sum\limits_{\sigma=1,\,R_j+1}\sum\limits_{\mu=1}^m
 (B'_{j\sigma\mu}U_j|_{\Gamma_{j\sigma}},\ F'_{j\sigma\mu}V_j|_{\Gamma_{j\sigma}})_{\Gamma_{j\sigma}}+\\
    +\sum\limits_{q=2}^{R_j}\sum\limits_{\mu=1}^m (B'_{jq\mu}U_j|_ {\Gamma_{jq}},\ 
    F'_{jq\mu}V_{j,q-1}|_ {\Gamma_{jq}}-F'_{jq\mu}V_{jq}|_ {\Gamma_{jq}}+\\
     +\sum\limits_k\sum\limits_{\sigma=1,\,R_k+1}\sum\limits_{s=1}^{S'_{jqk\sigma}}
     (\{\sum\limits_{\alpha=1}^m
     (\hat\Lambda'_{k\sigma\alpha jqs\mu})^*F'_{k\sigma\alpha}\}V_k)({\cal G}'_{jqk\sigma s}\cdot)|_ {\Gamma_{jq}})_
          {\Gamma_{jq}}\bigr\},
 \end{array}
\end{equation}
where the operators $\hat\Lambda_{k\sigma\alpha jqs\mu}$ and $\hat\Lambda'_{k\sigma\alpha jqs\mu}$ are obtained from
the operators $\Lambda_{k\sigma\alpha jqs\mu}$ and $\Lambda'_{k\sigma\alpha jqs\mu}$ by substituting 
$a_{k\sigma\alpha jqs\mu}^{\beta l}(\chi'_{jqk\sigma s})^l$ and ${a'}_{k\sigma\alpha jqs\mu}^{\beta l}(\chi'_{jqk\sigma s})^l$
for $a_{k\sigma\alpha jqs\mu}^{\beta l}$ and ${a'}_{k\sigma\alpha jqs\mu}^{\beta l}$ correspondingly.

Denoting
$$
 \begin{array}{c}
 C_{j\sigma\mu}=F'_{j\sigma\mu}\ (j=1,\ \dots,\ N;\ \sigma=1,\ R_j+1;\ \mu=1,\ \dots,\ m),\\
 T_{jq\nu}=F_{jq\nu}\ \mbox{for } \nu=1,\ \dots,\ m;\ T_{jq\nu}=F'_{jq,\nu-m}\ \mbox{for }\nu=m+1,\ \dots,\ 2m;\\
  T_{jq\nu k\sigma s}=\sum\limits_{\alpha=1}^m
     (\hat\Lambda_{k\sigma\alpha jqs\nu})^*F_{k\sigma\alpha}\ \mbox{for } \nu=1,\ \dots,\ m,\\
 T_{jq\nu k\sigma s}=\sum\limits_{\alpha=1}^m
     (\hat\Lambda'_{k\sigma\alpha jqs,\nu-m})^*F'_{k\sigma\alpha}\ \mbox{for } \nu=m+1,\ \dots,\ 2m\\
 (j,\ k=1,\ \dots,\ N;\ q=2,\ \dots,\ R_j;\ \sigma=1,\ R_k+1;\ s=1,\ \dots,\ S'_{jqk\sigma}),\\
 \end{array}
$$
we complete the proof.
\end{proof}

\begin{remark}\label{rGr}
 Formula~(\ref{eqGrP}) can be extended by continuity for the case $U_j\in H_a^{2m}(\Omega_j),$
 $V_{jt}\in H_{-a+2m}^{2m}(\Omega_{jt})$.
 Indeed, $C_0^\infty(\bar\Omega_j\backslash M)$ is dense in $H_a^{2m}(\Omega_j),$ $C_0^\infty(\bar\Omega_{jt}\backslash M)$ 
 is dense in $H_{-a+2m}^{2m}(\Omega_{jt})$; therefore there exist sequences
 $\{U_j^p\}_{p=1}^\infty\subset C_0^\infty(\bar\Omega_j\backslash \{0\})$ and
 $\{V_{jt}^q\}_{q=1}^\infty\subset C_0^\infty(\bar\Omega_{jt}\backslash \{0\})$ that converge to
 $U_j$ and $V_{jt}$ in $H_a^{2m}(\Omega_j)$ and $H_{-a+2m}^{2m}(\Omega_{jt})$ correspondingly. Green 
 formula~(\ref{eqGrP}) is valid for the functions $U_j^p$ and $V_{jt}^q$; passing to the limit as $p,\ q\to\infty,$ we obtain the
 Green formula for  $U_j$ and $V_{jt}$ (we can pass to the limit by virtue of  the Schwarz inequality and 
 Theorem~\ref{thTraceH}).
\end{remark}

The following two examples illustrate the Green formula.
\begin{example}\label{exGr1}
{\rm 
For simplicity we assume that $n=2,\ N=1.$ Put 
$
  K=\{y:\ r>0,\ b_1<\varphi<b_3\},
$
$
 K_{t}=\{y:\ r>0,\ b_t<\varphi<b_{t+1}\}\ (t=1,\ 2),
$
$
  \gamma_{q}=\{y:\ r>0,\ \varphi=b_{q}\}\ (q=1,\ 2,\ 3),
$
where $y=(y_1,\ y_2)\in{\mathbb R}^2;$ $0<b_1<b_2<b_3<2\pi.$ 

Let $n_1$ be the unit normal vector to $\gamma_1$ directed inside
$K_1$ and $n_2$, $n_3$ be the unit normal vectors to $\gamma_2$, $\gamma_3$ 
correspondingly directed inside $K_2.$

Consider the nonlocal problem
\begin{equation}\label{eqPEx1}
 -\triangle U=f(y) \quad (y\in K),
\end{equation}
\begin{equation}\label{eqBEx1}
 \begin{array}{c}
 U|_{\gamma_{1}}+\alpha U(\chi_{12}r,\ \varphi+\varphi_{12})|_{\gamma_{1}}=
  g_{1}(y)\ \quad (y\in\gamma_{1}),\\
 U|_{\gamma_{3}}=g_{3}(y) \quad (y\in\gamma_{3}).
 \end{array}
\end{equation}
Here $U(r,\ \varphi)$ is the function $U(y)$ written in the polar coordinates; 
$b_1+\varphi_{12}=b_2,$ $\chi_{12}>0;$ $\alpha\in{\mathbb R}.$

Take $U\in C_0^\infty(\bar K\backslash \{0\}),$ $V_t\in C^\infty(\bar K_t\backslash \{0\})$. 
Multiply $-\triangle U$ by $\bar V_t$ and
integrate over $K_t,$ $t=1,\ 2$; then using the formula of integration by parts, we get 
$$
 \begin{array}{c}
 \int\limits_{K_1}(-\triangle U)\cdot \bar V_1\,dy+\int\limits_{\gamma_1}U|_{\gamma_1}\cdot
  \frac{\displaystyle\partial \bar V_1}{\displaystyle\partial n_1}\Bigr|_{\gamma_1}\,d\gamma-
  \int\limits_{\gamma_2}U|_{\gamma_2}\cdot
  \frac{\displaystyle\partial \bar V_1}{\displaystyle\partial n_2}\Bigr|_{\gamma_2}\,d\gamma=\\
   =\int\limits_{K_1} U\cdot(-\triangle \bar V_1)\,dy+\int\limits_{\gamma_1}
 \frac{\displaystyle\partial U}{\displaystyle\partial n_1}\Bigr|_{\gamma_1}\cdot \bar V_1|_{\gamma_1}
 \,d\gamma-\int\limits_{\gamma_2}
  \frac{\displaystyle\partial U}{\displaystyle\partial n_2}\Bigr|_{\gamma_2}
  \cdot \bar V_1|_{\gamma_2}\,d\gamma,
 \end{array}
$$
$$
 \begin{array}{c}
 \int\limits_{K_2}(-\triangle U)\cdot \bar V_2\,dy+\int\limits_{\gamma_2}U|_{\gamma_2}\cdot
  \frac{\displaystyle\partial \bar V_2}{\displaystyle\partial n_2}\Bigr|_{\gamma_2}\,d\gamma+
  \int\limits_{\gamma_3}U|_{\gamma_3}\cdot
  \frac{\displaystyle\partial \bar V_2}{\displaystyle\partial n_3}\Bigr|_{\gamma_3}\,d\gamma=\\
   =\int\limits_{K_2} U\cdot(-\triangle \bar V_2)\,dy+\int\limits_{\gamma_2}
 \frac{\displaystyle\partial  U}{\displaystyle\partial n_2}\Bigr|_{\gamma_2}\cdot \bar V_2|_{\gamma_2}
 \,d\gamma+\int\limits_{\gamma_3}
  \frac{\displaystyle\partial U}{\displaystyle\partial n_3}\Bigr|_{\gamma_3}
  \cdot \bar V_2|_{\gamma_3}\,d\gamma.
 \end{array}
$$
Adding the last two equalities together, we obtain
\begin{equation}\label{eqGrEx11}
 \begin{array}{c}
   \sum\limits_t\int\limits_{K_t}(-\triangle U)\cdot \bar V_t\,dy+\int\limits_{\gamma_1}U|_{\gamma_1}\cdot
  \frac{\displaystyle\partial \bar V_1}{\displaystyle\partial n_1}\Bigr|_{\gamma_1}\,d\gamma+\\
  +\int\limits_{\gamma_3}U|_{\gamma_3}\cdot
  \frac{\displaystyle\partial \bar V_2}{\displaystyle\partial n_3}\Bigr|_{\gamma_3}\,d\gamma
 +\int\limits_{\gamma_2}U|_{\gamma_2}\cdot
  \Bigl(\frac{\displaystyle\partial \bar V_2}{\displaystyle\partial n_2}\Bigr|_{\gamma_2}-
      \frac{\displaystyle\partial \bar V_1}{\displaystyle\partial n_2}\Bigr|_{\gamma_2}\Bigr)\,d\gamma=\\
 =\sum\limits_t\int\limits_{K_t}U\cdot (-\triangle \bar V_t)\,dy+\int\limits_{\gamma_1}
 \frac{\displaystyle\partial U}{\displaystyle\partial n_1}\Bigr|_{\gamma_1}\cdot \bar V_1|_{\gamma_1}\,d\gamma
  +\int\limits_{\gamma_3}
  \frac{\displaystyle\partial U}{\displaystyle\partial n_3}\Bigr|_{\gamma_3}
  \cdot \bar V_2|_{\gamma_3}\,d\gamma+\\
 +\int\limits_{\gamma_2}
 \frac{\displaystyle\partial  U}{\displaystyle\partial n_2}\Bigr|_{\gamma_2}\cdot (\bar V_2|_{\gamma_2}-
 \bar V_1|_{\gamma_2})\,d\gamma.
 \end{array}
\end{equation}
But we have 
$$
 \begin{array}{c}
 \int\limits_{\gamma_1}U|_{\gamma_1}\cdot
  \frac{\displaystyle\partial \bar V_1}{\displaystyle\partial n_1}\Bigr|_{\gamma_1}\,d\gamma
 =\int\limits_{\gamma_1}(U|_{\gamma_1}+\alpha U(\chi_{12}r,\ \varphi+\varphi_{12})|_{\gamma_1})\cdot
  \frac{\displaystyle\partial \bar V_1}{\displaystyle\partial n_1}\Bigr|_{\gamma_1}\,d\gamma-\\
  -\int\limits_{\gamma_1}\alpha U(\chi_{12}r,\ \varphi+\varphi_{12})|_{\gamma_1}\cdot
  \frac{\displaystyle\partial \bar V_1}{\displaystyle\partial n_1}\Bigr|_{\gamma_1}\,d\gamma=\\
 =\int\limits_{\gamma_1}(U|_{\gamma_1}+\alpha U(\chi_{12}r,\ \varphi+\varphi_{12})|_{\gamma_1})\cdot
  \frac{\displaystyle\partial \bar V_1}{\displaystyle\partial n_1}\Bigr|_{\gamma_1}\,d\gamma-\\
  -\int\limits_{\gamma_2} U|_{\gamma_2}\cdot
  \alpha\chi'_{21}\frac{\displaystyle\partial \bar V_1}{\displaystyle\partial n_1}
   (\chi'_{21}r,\ \varphi+\varphi'_{21})\Bigr|_{\gamma_2}\,d\gamma,
 \end{array}
$$
where $\chi'_{21}=1/\chi_{12},\ \varphi'_{21}=-\varphi_{12}.$ This and~(\ref{eqGrEx11}) finally yield
$$ 
  \begin{array}{c}
   \sum\limits_t\int\limits_{K_t}(-\triangle U)\cdot \bar V_t\,dy+\int\limits_{\gamma_1}(U|_{\gamma_1}+
   \alpha U(\chi_{12}r,\ \varphi+\varphi_{12})|_{\gamma_1})\cdot
  \frac{\displaystyle\partial \bar V_1}{\displaystyle\partial n_1}\Bigr|_{\gamma_1}\,d\gamma+\\
  +\int\limits_{\gamma_3}U|_{\gamma_3}\cdot
  \frac{\displaystyle\partial \bar V_2}{\displaystyle\partial n_3}\Bigr|_{\gamma_3}\,d\gamma
 + \int\limits_{\gamma_2}
 \frac{\displaystyle\partial  U}{\displaystyle\partial n_2}\Bigr|_{\gamma_2}\cdot (\bar V_1|_{\gamma_2}-
 \bar V_2|_{\gamma_2})\,d\gamma=\\
\\
 =\sum\limits_t\int\limits_{K_t}U\cdot (-\triangle \bar V_t)\,dy+\int\limits_{\gamma_1}
 \frac{\displaystyle\partial U}{\displaystyle\partial n_1}\Bigr|_{\gamma_1}\cdot \bar V_1|_{\gamma_1}\,d\gamma+
  \int\limits_{\gamma_3}
  \frac{\displaystyle\partial  U}{\displaystyle\partial n_3}\Bigr|_{\gamma_3}
  \cdot\bar V_2|_{\gamma_3}\,d\gamma+\\
 \int\limits_{\gamma_2}U|_{\gamma_2}\cdot
  \Bigl(\frac{\displaystyle\partial \bar V_1}{\displaystyle\partial n_2}\Bigr|_{\gamma_2}-
      \frac{\displaystyle\partial \bar V_2}{\displaystyle\partial n_2}\Bigr|_{\gamma_2}+
 \alpha\chi'_{21}\frac{\displaystyle\partial \bar V_1}{\displaystyle\partial n_1}
   (\chi'_{21}r,\ \varphi+\varphi'_{21})\Bigr|_{\gamma_2}\Bigr)\,d\gamma.
 \end{array}
$$
}
\end{example}

\begin{example}\label{exGr2}
{\rm
 Using denotations of Example~\ref{exGr1}, consider the nonlocal problem
\begin{equation}\label{eqPEx2}
 -\triangle U=f(y) \quad (y\in K),
\end{equation}
\begin{equation}\label{eqBEx2}
 \begin{array}{c}
 \frac{\displaystyle\partial U}{\displaystyle\partial n_1}|_{\gamma_{1}}+
\alpha \frac{\displaystyle\partial U}{\displaystyle\partial r}(\chi_{12}r,\ \varphi+\varphi_{12})|_{\gamma_{1}}=
  g_{1}(y)\ \quad (y\in\gamma_{1}),\\
  \frac{\displaystyle\partial U}{\displaystyle\partial n_3}|_{\gamma_{3}}=g_{3}(y) \quad (y\in\gamma_{3}).
 \end{array}
\end{equation}
From formula~(\ref{eqGrEx11}) and equality
$$
 \begin{array}{c}
 \int\limits_{\gamma_1}
 \frac{\displaystyle\partial  U}{\displaystyle\partial n_1}\Bigr|_{\gamma_1}\cdot \bar V_1|_{\gamma_1}
 \,d\gamma
 =\int\limits_{\gamma_1}\Bigl(\frac{\displaystyle\partial U}{\displaystyle\partial n_1}\Bigr|_{\gamma_1}+
 \alpha\frac{\displaystyle\partial U}{\displaystyle\partial r}(\chi_{12}r,\ \varphi+\varphi_{12})\Bigr|_{\gamma_1}
  \Bigr)\cdot \bar V_1|_{\gamma_1}\,d\gamma+\\
 +\int\limits_{\gamma_2}U|_{\gamma_2}\cdot \alpha(\chi'_{21})^2
  \frac{\displaystyle\partial \bar V_1}{\displaystyle\partial r}(\chi'_{21}r,\ \varphi+\varphi'_{21})\Bigr|_{\gamma_2}
  \,d\gamma
 \end{array}
$$
(where $\chi'_{21}=1/\chi_{12},\ \varphi'_{21}=-\varphi_{12}$), we get the following Green formula:
$$
  \begin{array}{c}
   \sum\limits_t\int\limits_{K_t}(-\triangle U)\cdot \bar V_t\,dy+\int\limits_{\gamma_1}
 \Bigl(\frac{\displaystyle\partial U}{\displaystyle\partial n_1}\Bigr|_{\gamma_1}+
 \alpha\frac{\displaystyle\partial U}{\displaystyle\partial r}(\chi_{12}r,\ \varphi+\varphi_{12})\Bigr|_{\gamma_1}
  \Bigr)\cdot (-\bar V_1)|_{\gamma_1}\,d\gamma+\\
 +\int\limits_{\gamma_3}
  \frac{\displaystyle\partial  U}{\displaystyle\partial n_3}\Bigr|_{\gamma_3}\cdot(-\bar V_2)|_{\gamma_3}\,d\gamma+
 \int\limits_{\gamma_2}
 \frac{\displaystyle\partial  U}{\displaystyle\partial n_2}\Bigr|_{\gamma_2}\cdot (\bar V_1|_{\gamma_2}-
 \bar V_2|_{\gamma_2})\,d\gamma=\\
\\ 
=\sum\limits_t\int\limits_{K_t}U\cdot (-\triangle \bar V_t)\,dy+\int\limits_{\gamma_1} (-U)|_{\gamma_1}\cdot
  \frac{\displaystyle\partial \bar V_1}{\displaystyle\partial n_1}\Bigr|_{\gamma_1}\,d\gamma+
 \int\limits_{\gamma_3}(-U)|_{\gamma_3}\cdot
  \frac{\displaystyle\partial \bar V_2}{\displaystyle\partial n_3}\Bigr|_{\gamma_3}\,d\gamma+\\
+\int\limits_{\gamma_2}U|_{\gamma_2}\cdot
  \Bigl(\frac{\displaystyle\partial \bar V_1}{\displaystyle\partial n_2}\Bigr|_{\gamma_2}-
      \frac{\displaystyle\partial \bar V_2}{\displaystyle\partial n_2}\Bigr|_{\gamma_2}+
 \alpha(\chi'_{21})^2\frac{\displaystyle\partial \bar V_1}{\displaystyle\partial r}
   (\chi'_{21}r,\ \varphi+\varphi'_{21})\Bigr|_{\gamma_2}\Bigr)\,d\gamma.
 \end{array}
$$
}
\end{example}

\subsection{The Green formula with parameter $\eta$ in plane angles.}
For $n=2$, $j=1,\ \dots,\ N$, put 
$$
  K_j=\{y:\ r>0,\ b_{j1}<\varphi<b_{j,R_j+1}\},
$$
$$
  K_{jt}=\{y:\ r>0,\ b_{jt}<\varphi<b_{j,t+1}\}\ (t=1,\ \dots,\ R_j),
$$
$$
  \gamma_{jq}=\{y:\ r>0,\ \varphi=b_{jq}\}\ (q=1,\ \dots,\ R_j+1).
$$

Replace $D_z$ by $\eta$ in differential operators and consider the auxiliary nonlocal boundary value problem with 
parameter $\eta\in{\mathbb R}^{n-2}$ for $u=(u_1,\ \dots,\ u_N)$
\begin{equation}\label{eqPEta}
 {\cal P}_j(D_y,\ \eta)u_j=f_j(y) \quad (y\in K_j),
\end{equation}
\begin{equation}\label{eqBEta}
  \begin{array}{c}
  {\cal B}_{j\sigma\mu}(D_y,\ \eta)u=B_{j\sigma\mu}(D_y,\ \eta)u_j|_{\gamma_{j\sigma}}+\\
  +\sum\limits_{k,q,s}
                  (B_{j\sigma\mu kqs}(D_y,\ \eta)u_k)({\cal G}_{j\sigma kqs}y)|_{\gamma_{j\sigma}}
   =g_{j\sigma\mu}(y) \quad (y\in\gamma_{j\sigma})
  \end{array}
\end{equation}
$$ 
 (j=1,\ \dots,\ N;\ \sigma=1,\ R_j+1;\ \mu=1,\ \dots,\ m).
$$

For $u_{jt}\in C_0^\infty(\bar K_{jt}\backslash\{0\}),\ v_{jt}\in C^\infty(\bar K_{jt}\backslash \{0\})$
(or $u_{jt}\in C^\infty(\bar K_{jt}\backslash \{0\}),\ v_{jt}\in C_0^\infty(\bar K_{jt}\backslash \{0\})$), put
$$
 (u_{jt},\ v_{jt})_{K_{jt}}=\int\limits_{K_{jt}}u_{jt}\cdot \bar v_{jt}\, dy\quad (j=1,\ \dots,\ N;\ t=1,\ \dots,\ R_j).
$$
For $u_{\gamma_{jq}}\in C_0^\infty(\gamma_{jq}),\ 
v_{\gamma_{jq}}\in C^\infty(\gamma_{jq})$
(or $u_{\gamma_{jq}}\in C^\infty(\gamma_{jq}),\ 
v_{\gamma_{jq}}\in C_0^\infty(\gamma_{jq})$), put
$$
 (u_{\gamma_{jq}},\ v_{\gamma_{jq}})_{\gamma_{jq}}=
          \int\limits_{\gamma_{jq}}u_{\gamma_{jq}}\cdot \bar v_{\gamma_{jq}}\, d\gamma\quad (j=1,\ \dots,\ N;\ q=1,\ \dots,\ R_j+1).
$$
If we have functions $v_{jt}(y)$ defined in $K_{jt}$, then denote by $v_j(y)$ the function given by
$v_j(y)\equiv v_{jt}(y)$ for $y\in K_{jt}$.

\begin{theorem}\label{thGrPEta}
 Let ${\cal P}_j,$  $B_{j\sigma\mu},$ etc., be the operators from Theorem~\ref{thGrP}. Then for all
 $u_j\in C_0^\infty(\bar K_j\backslash \{0\}),$ 
 $v_{jt}\in C^\infty(\bar K_{jt}\backslash \{0\})$ (or $u_j\in C^\infty(\bar K_j\backslash \{0\}),$ 
 $v_{jt}\in C_0^\infty(\bar K_{jt}\backslash \{0\})$), the following Green formula with parameter $\eta$ is valid:
 \begin{equation}\label{eqGrPEta}
 \begin{array}{c}
  \sum\limits_j\bigl\{\sum\limits_t({\cal P}_j(D_y,\ \eta)u_j,\ v_{jt})_{K_{jt}}+\\
   + \sum\limits_{\sigma,\mu}
    ({\cal B}_{j\sigma\mu}(D_y,\ \eta)u,\ F_{j\sigma\mu}(D_y,\ \eta)v_j|_{\gamma_{j\sigma}})_
            {\gamma_{j\sigma}}+\\
 +\sum\limits_{q,\mu}
    (B_{jq\mu}(D_y,\ \eta)u_j|_{\gamma_{jq}},\ {\cal T}_{jq\mu}(D_y,\ \eta)v)_
     {\gamma_{jq}}\bigr\}=\\
\\
 =\sum\limits_j\bigl\{\sum\limits_t(u_j,\ {\cal Q}_j(D_y,\ \eta)v_{jt})_{K_{jt}}+\\
    +\sum\limits_{\sigma,\mu}
    (B'_{j\sigma\mu}(D_y,\ \eta)u_j|_{\gamma_{j\sigma}},\ C_{j\sigma\mu}(D_y,\ \eta)v_j|_{\gamma_{j\sigma}})_
    {\gamma_{j\sigma}}+\\
 +\sum\limits_{q,\mu}
 (B'_{jq\mu}(D_y,\ \eta)u_j|_{\gamma_{jq}},\ {\cal T}_{jq,m+\mu}(D_y,\ \eta)v)_
   {\gamma_{jq}}\bigr\}.
 \end{array} 
 \end{equation}
 Here ${\cal B}_{j\sigma\mu}(D_y,\ \eta)$ is given by~(\ref{eqBEta});
 $$
 \begin{array}{c}
 {\cal T}_{jq\nu}(D_y,\ \eta)v=T_{jq\nu}(D_y,\ \eta)v_{j,q-1}|_{\gamma_{jq}}-
 T_{jq\nu}(D_y,\ \eta)v_{jq}|_{\gamma_{jq}}+\\
  +\sum\limits_{k,\sigma,s}
     (T_{jq\nu k\sigma s}(D_y,\ \eta)v_k) ({\cal G}'_{jqk\sigma s}y)|_{\gamma_{jq}}\\
  (\nu=1,\ \dots,\ 2m);
 \end{array}
 $$
 ${\cal G}'_{jqk\sigma s}$ is the transformation defined in Theorem~\ref{thGrP}.
\end{theorem}
\begin{proof}
Introduce the functions $\psi_1,\ \psi_2\in C_0^\infty({\mathbb R}^{n-2})$ such that
$$
 \psi_1(z)=0\ \mbox{for } |z|>1,\    \int\limits_{{\mathbb R}^{n-2}}\psi_1(z)dz=1,
$$ 
$$
 \psi_2(z)=1 \ \mbox{for } |z|<1,\ \psi_2(z)=0\ \mbox{for } |z|>2.
$$
Substituting $U_j(y,\ z)=e^{i(\eta,\ z)}\psi_1(z)u_j(y),$ $V_{jt}(y,\ z)=e^{i(\eta,\ z)}\psi_2(z)v_{jt}(y)$ into equality~(\ref{eqGrP}),
we get~(\ref{eqGrPEta}).
\end{proof}

\begin{remark}\label{rGrEta}
 Replacing in Remark~\ref{rGr} $H_a^{2m}(\cdot)$ and 
 $H_{-a+2m}^{2m}(\cdot)$ by $E_a^{2m}(\cdot)$ and
 $E_{-a+2m}^{2m}(\cdot)$ correspondingly and Theorem~\ref{thTraceH} by Theorem~\ref{thTraceE}, we see that
 formula~(\ref{eqGrPEta}) can be extended by continuity for the 
 case $u_j\in E_a^{2m}(K_j),$ $v_{jt}\in E_{-a+2m}^{2m}(K_{jt})$.
\end{remark}

\subsection{The Green formula with parameter $\lambda$ on arcs.}
Put 
$
  \Pi_j=\{(\varphi,\ \tau):\ b_{j1}<\varphi<b_{j,R_j+1},\ \tau\in{\mathbb R}\},
$
$  
  \Pi_{jt}=\{(\varphi,\ \tau):\ b_{jt}<\varphi<b_{j,t+1},\ \tau\in{\mathbb R}\}\ (t=1,\ \dots,\ R_j).
$

For $u_{jt}\in C_0^\infty(\bar\Pi_{jt}),\  v_{jt}\in C^\infty(\bar\Pi_{jt})$
(or $u_{jt}\in C^\infty(\bar\Pi_{jt}),\ v_{jt}\in C_0^\infty(\bar\Pi_{jt})$), denote
$$
 \begin{array}{c}
 (u_{jt},\ v_{jt})_{\Pi_{jt}}=\int\limits_{-\infty}^\infty\int\limits_{b_{jt}}^{b_{j,t+1}}
  u_{jt}(\varphi,\ \tau)\cdot \overline{v_{jt}(\varphi,\ \tau)}\, d\varphi d\tau\\
  (j=1,\ \dots,\ N;\ t=1,\ \dots,\ R_j).
 \end{array}
$$
For $\psi\in C_0^\infty({\mathbb R}),\ \xi\in C^\infty({\mathbb R})$ (or $\psi\in C^\infty({\mathbb R}),\ \xi\in C_0^\infty({\mathbb R})$), denote
$
 (\psi,\ \xi)_{\mathbb R}=
          \int\limits_{-\infty}^\infty \psi(\tau)\cdot \overline{\xi(\tau)}\, d\tau.
$
For $\tilde U_{jt},\ \tilde V_{jt}\in C^\infty([b_{jt},\ b_{j,t+1}])$, we also denote
$$
 (\tilde U_{jt},\ \tilde V_{jt})_{(b_{jt},\ b_{j,t+1})}=\int\limits_{b_{jt}}^{b_{j,t+1}}
  \tilde U_{jt}(\varphi)\cdot \overline{\tilde V_{jt}(\varphi)}\, d\varphi\quad (j=1,\ \dots,\ N;\ t=1,\ \dots,\ R_j).
$$
And finally for $d,\ e\in{\mathbb C}$, we put $(d,\ e)_{\mathbb C}=d\cdot \bar e.$

If we have functions $\tilde V_{jt}(\varphi)$ defined in $[b_{jt},\ b_{j,t+1}]$, then denote by $\tilde V_j(\varphi)$ the function given by
$\tilde V_j(\varphi)\equiv \tilde V_{jt}(\varphi)$ for $\varphi\in(b_{jt},\ b_{j,t+1})$.

Put $D_z=0$ and write the differential operators in the polar coordinates:
${\cal P}_j(D_y,\ 0)=r^{-2m}\tilde{\cal P}_j(\varphi,\ D_\varphi,\ rD_r),$ 
$B_{j\sigma\mu}(D_y,\ 0)=r^{-m_{j\sigma\mu}}\tilde B_{j\sigma\mu}(\varphi,\ D_\varphi,\ rD_r)$, ets.

Consider nonlocal boundary value problem~(\ref{eqPLambda}), (\ref{eqBLambda}) with parameter $\lambda.$
\begin{theorem}\label{thGrPLambda}
 Let ${\cal P}_j,$  $B_{j\sigma\mu},$ etc., be the operators from Theorem~\ref{thGrP}. Then for all
 $\tilde U_j\in C^\infty([b_{j1},\ b_{j,R_j+1}])$, $\tilde V_{jt}\in C^\infty([b_{jt},\ b_{j,t+1}])$, the following Green formula with
 parameter $\lambda$ is valid:
 \begin{equation}\label{eqGrPLambda}
 \begin{array}{c}
  \sum\limits_j\bigl\{\sum\limits_t
   (\tilde{\cal P}_j(\varphi,\ D_\varphi,\ \lambda)\tilde U_j,\ \tilde V_{jt})_{(b_{jt},\ b_{j,t+1})}+\\
  +\sum\limits_{\sigma,\mu}
    (\tilde{\cal B}_{j\sigma\mu}(\varphi,\ D_\varphi,\ \lambda)\tilde U,\ 
    \tilde F_{j\sigma\mu}(\varphi,\ D_\varphi,\ \lambda')\tilde V_j|_{\varphi=b_{j\sigma}})_{\mathbb C}+\\
   +\sum\limits_{q,\mu}
  (\tilde B_{jq\mu}(\varphi,\ D_\varphi,\ \lambda)\tilde U_j|_{\varphi=b_{jq}},\ 
          \tilde{\cal T}_{jq\mu}(\varphi,\ D_\varphi,\ \lambda')\tilde V)_{\mathbb C}\bigr\}=\\
 \\
    =\sum\limits_j\bigl\{\sum\limits_t
    (\tilde U_j,\ \tilde{\cal Q}_j(\varphi,\ D_\varphi,\ \lambda')\tilde V_{jt})_{(b_{jt},\ b_{j,t+1})}+\\
   +\sum\limits_{\sigma,\mu}
    (\tilde B'_{j\sigma\mu}(\varphi,\ D_\varphi,\ \lambda)\tilde U_j|_{\varphi=b_{j\sigma}},
    \tilde C_{j\sigma\mu}(\varphi,\ D_\varphi,\ \lambda')\tilde V_j|_{\varphi=b_{j\sigma}})_{\mathbb C}+\\
  +\sum\limits_{q,\mu} 
  (\tilde B'_{jq\mu}(\varphi,\ D_\varphi,\ \lambda)\tilde U_j|_{\varphi=b_{jq}},\ 
   \tilde{\cal T}_{jq,m+\mu}(\varphi,\ D_\varphi,\ \lambda')\tilde V)_{\mathbb C}\bigr\}.
 \end{array} 
 \end{equation}
 Here $\tilde{\cal B}_{j\sigma\mu}(\varphi,\ D_\varphi,\ \lambda)$ is given by~(\ref{eqBLambda});
 $$
 \begin{array}{c}
 \tilde{\cal T}_{jq\nu}(\varphi,\ D_\varphi,\ \lambda')\tilde V=
 \tilde T_{jq\nu}(\varphi,\ D_\varphi,\ \lambda')\tilde V_{j,q-1}(\varphi)|_{\varphi=b_{jq}}-\\
 -\tilde T_{jq\nu}(\varphi,\ D_\varphi,\ \lambda')\tilde V_{jq}(\varphi)|_{\varphi=b_{jq}}+\\
  +\sum\limits_{k,\sigma,s}
     e^{(i\lambda'-(\nu-1))\ln\chi'_{jqk\sigma s}}\tilde T_{jq\nu k\sigma s}(\varphi,\ D_\varphi,\ \lambda')
            \tilde V_k(\varphi+\varphi'_{jqk\sigma})|_{\varphi=b_{jq}};
 \end{array}
 $$
 $\lambda'=\bar\lambda-2i(m-1);$ $\varphi'_{jqk\sigma}$ and $\chi'_{jqk\sigma s}$ are the rotation angles and 
 the expansion coefficients correspondingly defined in Theorem~\ref{thGrP}.
\end{theorem}
\begin{proof}
Put $r=e^\tau,\ v_{jt}=r^{2m-2}w_{jt}$, $w_j(\varphi,\ \tau)\equiv w_{jt}(\varphi,\ \tau)$ for $(\varphi,\ \tau)\in\Pi_{jt}$. Then
from formula~(\ref{eqGrPEta}) for $\eta=0$, we obtain
\begin{equation}\label{eqGrPLambda1}
 \begin{array}{c}
  \sum\limits_j\Bigl\{\sum\limits_t
   \Bigl(\tilde{\cal P}_j(\varphi,\ D_\varphi,\ D_\tau) u_j,\ w_{jt}\Bigr)_{\Pi_{jt}}+\\
  +\sum\limits_{\sigma,\mu}
    \Bigl( \tilde{\cal B}_{j\sigma\mu}(\varphi,\ D_\varphi,\ D_\tau)u,
 \tilde F_{j\sigma\mu}(\varphi,\ D_\varphi,\ D_\tau-2i(m-1))w_j|_{\varphi=b_{j\sigma}}\Bigr)_{\mathbb R}+\\
  +\sum\limits_{q,\mu}
   \Bigl(\tilde B_{jq\mu}(\varphi,\ D_\varphi,\ D_\tau)u_j|_{\varphi=b_{j\sigma}},
 \tilde {\cal T}_{jq\mu}(\varphi,\ D_\varphi,\ D_\tau-2i(m-1))w\Bigr)_{\mathbb R}\Bigr\}=\\
 \\
    =\sum\limits_j\Bigl\{\sum\limits_t
    \Bigl(u_j,\ \tilde{\cal Q}_j(\varphi,\ D_\varphi,\ D_\tau-2i(m-1))w_{jt}\Bigr)_{\Pi_{jt}}+\\
   \sum\limits_{\sigma,\mu}
    \Bigl(\tilde B'_{j\sigma\mu}(\varphi,\ D_\varphi,\ D_\tau)u_j|_{\varphi=b_{j\sigma}},
    \tilde C_{j\sigma\mu}(\varphi,\ D_\varphi,\ D_\tau-2i(m-1))w_j|_{\varphi=b_{j\sigma}}\Bigr)_{\mathbb R}+\\
  + \sum\limits_{q,\mu}
  \Bigl(\tilde B'_{jq\mu}(\varphi,\ D_\varphi,\ D_\tau)w_j|_{\varphi=b_{jq}},
    \tilde {\cal T}_{jq,m+\mu}(\varphi,\ D_\varphi,\ D_\tau-2i(m-1))w\Bigr)_{\mathbb R}\Bigr\},
 \end{array} 
\end{equation}
where
$$
 \begin{array}{c}
 \tilde{\cal B}_{j\sigma\mu}(\varphi,\ D_\varphi,\ D_\tau)u=
 \tilde B_{j\sigma\mu}(\varphi,\ D_\varphi,\ \ D_\tau)u_j|_{\varphi=b_{j\sigma}}+\\
 +\sum\limits_{k,q,s}e^{-m_{j\sigma\mu}\ln\chi_{j\sigma kqs}}\tilde B_{j\sigma\mu kqs}(\varphi,\ D_\varphi,\ D_\tau)
  u_k(\varphi+\varphi_{j\sigma kq},\ \tau+\ln\chi_{j\sigma kqs})|_{\varphi=b_{j\sigma}},
 \end{array}
$$
$$
 \begin{array}{c}
 \tilde {\cal T}_{jq\nu}(\varphi,\ D_\varphi,\ D_\tau-2i(m-1))w=\\
 =\tilde T_{jq\nu}(\varphi,\ D_\varphi,\ D_\tau-2i(m-1))w_{j,q-1}|_{\varphi=b_{jq}}-\\
 -\tilde T_{jq\nu}(\varphi,\ D_\varphi,\ D_\tau-2i(m-1))w_{jq}|_{\varphi=b_{jq}}
   +\sum\limits_{k,\sigma,s}
     e^{(2(m-1)-(\nu-1))\ln\chi'_{jqk\sigma s}}\times\\
 \times \tilde T_{jq\nu k\sigma s}(\varphi,\ D_\varphi,\ D_\tau-2i(m-1))
   w_k(\varphi+\varphi'_{jqk\sigma},\ \tau+\ln\chi'_{jqk\sigma s})|_{\varphi=b_{jq}}.
 \end{array}
$$
Introduce the functions $\psi_1,\ \psi_2\in C_0^\infty({\mathbb R})$ such that
$$
 \psi_1(\tau)=0\ \mbox{for } |\tau|>1,\    \int\limits_{-\infty}^\infty\psi_1(\tau)d\tau=1,
$$ 
$$
 \psi_2(\tau)=1 \ \mbox{for } |\tau|<1,\ \psi_2(\tau)=0\ \mbox{for } |\tau|>2.
$$
Substituting $u_j(\varphi,\ \tau)=e^{i\lambda\tau}\psi_1(\tau)\tilde U_j(\varphi)$,
$w_{jt}(\varphi,\ \tau)=e^{i\bar\lambda\tau}\psi_2(\tau)\tilde V_{jt}(\varphi)$ into equality~(\ref{eqGrPLambda1}), we 
obtain~(\ref{eqGrPLambda}).
\end{proof}

\begin{remark}\label{rGrLambda}
 Formula~(\ref{eqGrPLambda}) can be extended by continuity for the case
 $\tilde U_j\in W^{2m}(b_{j1},\ b_{j,R_j+1}),$ $\tilde V_{jt}\in W^{2m}(b_{jt},\ b_{j,t+1})$
 (see remark~2.2 \cite[Chapter 2, \S 2.3]{LM}).
\end{remark}
\section{Nonlocal elliptic transmission problems. Reduction to problems with homogeneous nonlocal and boundary 
conditions}\label{sectTransStatement}
\subsection{Nonlocal problems in dihedral angles.}
Put $V=(V_1,\ \dots,\ V_N)$, $f=(f_1,\ \dots,\ f_N)$. Here the functions $V_j(x)$ ($f_j(x)$) 
are defined in $\Omega_j$ ($j=1,\ \dots,\ N$).
As before, we shall denote by $V_{jt}$ ($f_{jt}$) the restriction of $V_j$ ($f_j$) to $\Omega_{jt}$. Then we see that 
Green formula~(\ref{eqGrP}) generates the problem, which is formally adjoint to problem~(\ref{eqP}), (\ref{eqB})
\begin{equation}\label{eqQ}
  {\cal Q}_j(D_y,\ D_z)V_{jt}=f_{jt}(x) \quad (x\in\Omega_{jt};\ t=1,\ \dots,\ R_j),
\end{equation}
\begin{equation}\label{eqC}
  \begin{array}{c}
  {\cal C}_{j1\mu}(D_y,\ D_z)V=
      C_{j1\mu}(D_y,\ D_z)V_{j1}(x)|_{\Gamma_{j1}}=g_{j1\mu}(x) \quad (x\in\Gamma_{j1}),\\
  {\cal C}_{j,R_j+1,\mu}(D_y,\ D_z)V=
     C_{j,R_j+1,\mu}(D_y,\ D_z)V_{jR_j}(x)|_{\Gamma_{j,R_j+1}}=\\
  =g_{j,R_j+1,\mu}(x) 
      \quad (x\in\Gamma_{j,R_j+1}),
  \end{array}
\end{equation}
\begin{equation}\label{eqT}
  \begin{array}{c}
  {\cal T}_{jq\nu}(D_y,\ D_z)V=
     T_{jq\nu}(D_y,\ D_z)V_{j,q-1}(x)|_{\Gamma_{jq}}-T_{jq\nu}(D_y,\ D_z)V_{jq}(x)|_{\Gamma_{jq}}+\\
   +\sum\limits_{k,\sigma,s}
                   (T_{jq\nu k\sigma s}(D_y,\ D_z)V_k)({\cal G}'_{jqk\sigma s}y,\ z)|_{\Gamma_{jq}}=h_{jq\nu}(x) \quad (x\in\Gamma_{jq})
  \end{array}
\end{equation}
$$
 (j=1,\ \dots,\ N;\ \mu=1,\ \dots,\ m;\ q=2,\ \dots,\ R_j;\ \nu=1,\ \dots,\ 2m).
$$
Here ${\cal Q}_j$ is formally adjoint to ${\cal P}_j$; the operators $C_{j\sigma\mu}$, $T_{jq\nu}$, $T_{jq\nu k\sigma s}$
are of orders $m'_{j\sigma\mu}$, $\nu-1$, $\nu-1$ correspondingly;
${\cal G}'_{jqk\sigma s}$ is the operator of rotation by the angle~$\varphi'_{jqk\sigma}=-\varphi_{k\sigma jq}$ and expansion 
by~$\chi'_{jqk\sigma s}=1/\chi_{k\sigma j qs}$ times in the plane~$\{y\}$ such that $b_{jq}+\varphi'_{jqk\sigma}=b_{k\sigma},$
$0<\chi'_{jqk\sigma s}$; $j,\ k=1,\ \dots,\ N;$ $q=2,\ \dots,\ R_j;$ $\sigma=1,\ R_k+1;$ 
$s=1,\ \dots,\ S'_{jqk\sigma}=S_{k\sigma jq}.$

Problem~(\ref{eqQ})--(\ref{eqT}) is a system of $R_1+\dots+R_N$ equations for functions $V_{jt}$ with
boundary conditions~(\ref{eqC}) and nonlocal transmission conditions~(\ref{eqT}). We shall say that 
problem~(\ref{eqQ})--(\ref{eqT}) is {\it a nonlocal transmission problem}. 

Let us write the nonlocal transmission problems, which are formally adjoint to nonlocal boundary value problems of 
Examples~\ref{exGr1} and~\ref{exGr2}.
\begin{example}
{\rm 
 From Example~\ref{exGr1}, it follows that the problem
$$ 
 -\triangle V_t=f_t(y)\quad (y\in K_t;\ t=1,\ 2),
$$
$$ 
  \begin{array}{c}
 V_1|_{\gamma_{1}}=g_{1}(y) \quad (y\in\gamma_{1}),\\
 V_2|_{\gamma_{3}}=g_{3}(y) \quad (y\in\gamma_{3}),
 \end{array}
$$
$$ 
 \begin{array}{c}
 V_1|_{\gamma_2}-V_2|_{\gamma_2}=h_{21}(y)\quad (y\in\gamma_{2}),\\
 \frac{\displaystyle\partial V_1}{\displaystyle\partial n_2}\Bigr|_{\gamma_2}-
      \frac{\displaystyle\partial V_2}{\displaystyle\partial n_2}\Bigr|_{\gamma_2}+
 \alpha\chi'_{21}\frac{\displaystyle\partial V_1}{\displaystyle\partial n_1}
   (\chi'_{21}r,\ \varphi+\varphi'_{21})\Bigr|_{\gamma_2}=h_{22}(y)\quad (y\in\gamma_{2})
 \end{array}
$$
is formally adjoint to problem~(\ref{eqPEx1}), (\ref{eqBEx1}).
}
\end{example}

\begin{example}
{\rm 
 From Example~\ref{exGr2}, it follows that the problem
$$ 
 -\triangle V_t=f_t(y)\quad (y\in K_t;\ t=1,\ 2),
$$
$$ 
 \begin{array}{c}
  \frac{\displaystyle\partial V_1}{\displaystyle\partial n_1}\Bigr|_{\gamma_1}=g_{1}(y) \quad (y\in\gamma_{1}),\\
  \frac{\displaystyle\partial V_2}{\displaystyle\partial n_3}\Bigr|_{\gamma_1}=g_{3}(y) \quad (y\in\gamma_{3}),
 \end{array}
$$
$$ 
 \begin{array}{c}
 V_1|_{\gamma_2}-V_2|_{\gamma_2}=h_{21}(y)\quad (y\in\gamma_{2}),\\
 \frac{\displaystyle\partial V_1}{\displaystyle\partial n_2}\Bigr|_{\gamma_2}-
      \frac{\displaystyle\partial V_2}{\displaystyle\partial n_2}\Bigr|_{\gamma_2}+
 \alpha(\chi'_{21})^2\frac{\displaystyle\partial V_1}{\displaystyle\partial r}
   (\chi'_{21}r,\ \varphi+\varphi'_{21})\Bigr|_{\gamma_2}=h_{22}(y)\quad (y\in\gamma_{2})
 \end{array}
$$
is formally adjoint to problem~(\ref{eqPEx2}), (\ref{eqBEx2}).
}
\end{example}

From Theorem~\ref{thGrP}, it follows that the
following conditions hold (see~\cite[Chapter 2, \S\S 1.2, 1.4]{LM}).
\begin{condition}\label{condEllipQ}
 For all $j=1,\ \dots,\ N$, the operators ${\cal Q}_j(D_y,\ D_z)$ are properly elliptic.
\end{condition}
\begin{condition}\label{condComplC}
 For all $j=1,\ \dots,\ N;\ \sigma=1,\ R_j+1$, the system $\{C_{j\sigma\mu}(D_y,\ D_z)\}_{\mu=1}^m$ is normal and 
 covers the operator ${\cal Q}_j(D_y,\ D_z)$ on $\Gamma_{j\sigma}.$
\end{condition} 
\begin{condition}\label{condNormT}
 For all $j=1,\ \dots,\ N;\ q=2,\ \dots,\ R_j$, the system $\{T_{jq\nu}(D_y,\ D_z)\}_{\nu=1}^{2m}$ is normal
 on $\Gamma_{jq}$.
\end{condition}

\begin{remark} 
 One can easily prove that under condition~\ref{condNormT}, 
 the system $\{T_{jq\nu}(D_y,\ D_z),\ T_{jq\nu}(D_y,\ D_z)\}_{\nu=1}^{2m}$ jointly covers the operator ${\cal Q}_j(D_y,\ D_z)$ 
 on $\Gamma_{jq}$ in the sense of~\cite{Sh65}.
\end{remark}

Consider the space ${\cal H}_a^l(\Omega_j)=\bigoplus\limits_{t=1}^{R_j}H_a^l(\Omega_{jt})$ with the norm
$
 \|V_j\|_{{\cal H}_a^l(\Omega_j)}=\left(\sum\limits_{t=1}^{R_j}\|V_{jt}\|_{H_a^l(\Omega_{jt})}^2\right)^{1/2}.
$

Introduce the spaces of vector--functions
$$
 {\cal H}_a^{l+2m,\,N}(\Omega)=\prod_{j=1}^N {\cal H}_a^{l+2m}(\Omega_j),\ 
 {\cal H}_a^{l,\,N}(\Omega,\ \Gamma)=\prod_{j=1}^N {\cal H}_a^l(\Omega_j,\ \Gamma_j),
$$
$$
 \begin{array}{c}
 {\cal H}_a^l(\Omega_j,\ \Gamma_j)={\cal H}_a^l(\Omega_j)\times\\
   \times\prod_{\sigma=1,\,R_j+1}\prod_{\mu=1}^m H_a^{l+2m-m'_{j\sigma\mu}-1/2}(\Gamma_{j\sigma})
   \times\prod_{q=2}^{R_j}\prod_{\nu=1}^{2m} H_a^{l+2m-\nu+1/2}(\Gamma_{jq}).
 \end{array}
$$

We study solutions $V=(V_1,\ \dots,\ V_N)\in{\cal H}_a^{l+2m,\,N}(\Omega)$ for problem~(\ref{eqQ})--(\ref{eqT}) supposing that
$f=\{f_j,\ g_{j\sigma\mu},\ h_{jq\nu}\}\in {\cal H}_a^{l,\,N}(\Omega,\ \Gamma).$ Introduce the bounded operator
${\cal M} : {\cal H}_a^{l+2m,\,N}(\Omega)\to {\cal H}_a^{l,\,N}(\Omega,\ \Gamma)$ 
corresponding to problem~(\ref{eqQ})--(\ref{eqT}) and given by
$$
{\cal M}V=\{W_j,\ {\cal C}_{j\sigma\mu}(D_y,\ D_z)V,\ {\cal T}_{jq\nu}(D_y,\ D_z)V\}
$$ 
Here ${\cal C}_{j\sigma\mu}(D_y,\ D_z)V$ and ${\cal T}_{jq\nu}(D_y,\ D_z)V$ are given by~(\ref{eqC}) and~(\ref{eqT}) 
correspondingly; $W_j(x)\equiv {\cal Q}_j(D_y,\ D_z)V_{jt}(x)$ for $x\in\Omega_{jt}$. (Notice that we cannot write
$W_j\equiv{\cal Q}_j(D_y,\ D_z)V_j$ for $x\in\Omega_j$ because $V_j\in{\cal H}_a^{l+2m}(\Omega_j)$ may have 
discontinuity on $\Gamma_{jq}$, $q=2,\ \dots,\ R_j$.)

\begin{lemma}\label{lHomog'} 
 For any $g_{j\sigma\mu}\in H_a^{l+2m-m'_{j\sigma\mu}-1/2}(\Gamma_{j\sigma}),$ 
 $h_{jq\nu}\in H_a^{l+2m-\nu+1/2}(\Gamma_{jq})$
 $(j=1,\ \dots,\ N;\ \sigma=1,\ R_j+1;\ \mu=1,\ \dots,\ m;\ q=2,\ \dots,\ R_j;\ \nu=1,\ \dots,\ 2m)$,
 there exists a vector--function $V\in  {\cal H}_a^{l+2m,\,N}(\Omega)$ such that
  $$
  {\cal C}_{j\sigma\mu}(D_y,\ D_z)V=g_{j\sigma\mu}(x)\ (x\in\Gamma_{j\sigma}),\quad
  {\cal T}_{jq\nu}(D_y,\ D_z)V=h_{jq\nu}(x)\ (x\in\Gamma_{jq}),
  $$
 $$
  \begin{array}{c}
  \|V\|_{{\cal H}_a^{l+2m,\,N}(\Omega)}\le 
  c\sum\limits_j\Bigl\{\sum\limits_{\sigma,\,\mu}\|g_{j\sigma\mu}\|_{H_a^{l+2m-m'_{j\sigma\mu}-1/2}
      (\Gamma_{j\sigma})}+\\
 +\sum\limits_{q,\,\nu}\|h_{jq\nu}\|_{H_a^{l+2m-\nu+1/2}(\Gamma_{jq})}\Bigr\},
 \end{array}
 $$
 where $c>0$ is independent of $g_{j\sigma\mu},\ h_{jq\nu}.$
\end{lemma}

\begin{proof}
By virtue of condition~\ref{condComplC} and lemma~3.1 \cite{MP},  there exists a 
vector--function $W\in H_a^{l+2m,\,N}(\Omega)$ such that
\begin{equation}\label{eqHomog1'}
  {\cal C}_{j\sigma\mu}(D_y,\ D_z)W=g_{j\sigma\mu}(x)\ \quad (x\in\Gamma_{j\sigma}),
\end{equation} 
\begin{equation}\label{eqHomog2'}
  \|W\|_{H_a^{l+2m,\,N}(\Omega)}\le 
  k_1\sum_{j,\,\sigma,\,\mu}\|g_{j\sigma\mu}\|_{H_a^{l+2m-m'_{j\sigma\mu}-1/2}(\Gamma_{j\sigma})}.
\end{equation}

By virtue of condition~\ref{condNormT} and lemma~3.1 \cite{MP},
for all $j=1,\ \dots,\ N$ and $q=2,\ \dots,\ R_j$ there exists a function
$\hat W_{j,q-1}\in H_a^{l+2m}(\Omega_{j,q-1})$ such that
\begin{equation}\label{eqHomog3'}
  \begin{array}{c}
     T_{jq\nu}(D_y,\ D_z)\hat W_{j,q-1}(x)|_{\Gamma_{jq}}=h_{jq\nu}(x)-\\
     -\sum\limits_{k,\,\sigma,\,s}(T_{jq\nu k\sigma s}(D_y,\ D_z)W_k)({\cal G}'_{jqk\sigma s}y,\ z)|_{\Gamma_{jq}}\quad
     (x\in\Gamma_{jq}),
 \end{array}
\end{equation}
\begin{equation}\label{eqHomog4'}
 \begin{array}{c}
  \|\hat W_{j,q-1}\|_{H_a^{l+2m}(\Omega_{j,q-1})}\le 
  k_2 \sum\limits_{\nu} \|h_{jq\nu}(x)-\\
  -\sum\limits_{k,\,\sigma,\,s}(T_{jq\nu k\sigma s}(D_y,\ D_z)W_k)
     ({\cal G}'_{jqk\sigma s}y,\ z)|_{\Gamma_{jq}}\|_{H_a^{l+2m-\nu+1/2}(\Gamma_{jq})}.
 \end{array}
\end{equation}

Since the functions $\zeta_{jq}$ defined by formula~(\ref{eqZeta}) are the multiplicators in the spaces
${\cal H}_a^{l+2m}(\Omega_j),$ from~(\ref{eqHomog1'})--(\ref{eqHomog4'}), it follows that the functions
$$
V_j(x)=\left\{
 \begin{array}{ll}
  \zeta_{j1}W_{j1}(x)+\zeta_{j2}\hat W_{j1}(x) & \mbox{for } x\in\Omega_{j1},\\
  \zeta_{j,t+1}\hat W_{jt}(x) & \mbox{for } x\in\Omega_{jt}\ (t=2,\ \dots,\ R_j-1),\\
  \zeta_{j,R_j+1}W_{jR_j}(x) &  \mbox{for } x\in\Omega_{jR_j}
 \end{array}
 \right.
$$
satisfy the conditions of the Lemma.
\end{proof}

\subsection{Nonlocal problems with parameter $\theta$ in plane angles.}
Put $v=(v_1,\ \dots,\ v_N)$, $f=(f_1,\ \dots,\ f_N)$. Here the functions $v_j(y)$ ($f_j(y)$) are defined in $K_j$ ($j=1,\ \dots,\ N$).
As before, we shall denote by $v_{jt}$ ($f_{jt}$) the restriction of $v_j$ ($f_j$) to $K_{jt}$. Then we see that 
Green formula~(\ref{eqGrPEta}) (for $\eta=\theta\in S^{n-3}=\{z\in{\mathbb R}^{n-2}:\ |z|=1\}$) generates the problem, which 
is formally adjoint to problem~(\ref{eqPTheta}), (\ref{eqBTheta})
\begin{equation}\label{eqQTheta}
{\cal Q}_j(D_y,\ \theta)v_{jt}=f_{jt}(y) \quad (y\in K_{jt};\ t=1,\ \dots,\ R_j),
\end{equation}
\begin{equation}\label{eqCTheta}
  \begin{array}{c}
  {\cal C}_{j1\mu}(D_y,\ \theta)v=
      C_{j1\mu}(D_y,\ \theta)v_{j1}(y)|_{\gamma_{j1}}=g_{j1\mu}(y) \quad (y\in\gamma_{j1}),\\
  {\cal C}_{j,R_j+1,\mu}(D_y,\ \theta)v=
     C_{j,R_j+1,\mu}(D_y,\ \theta)v_{jR_j}(y)|_{\gamma_{j,R_j+1}}=\\
  =g_{j,R_j+1,\mu}(y) 
      \quad (y\in\gamma_{j,R_j+1}),
  \end{array}
\end{equation}
\begin{equation}\label{eqTTheta}
  \begin{array}{c}
  {\cal T}_{jq\nu}(D_y,\ \theta)v=
     T_{jq\nu}(D_y,\ \theta)v_{j,q-1}(y)|_{\gamma_{jq}}-T_{jq\nu}(D_y,\ \theta)v_{jq}(y)|_{\gamma_{jq}}+\\
   +\sum\limits_{k,\sigma,s}
       (T_{jq\nu k\sigma s}(D_y,\ \theta)v_k)({\cal G}'_{jqk\sigma s}y)|_{\gamma_{jq}}=h_{jq\nu}(y) \quad (y\in\gamma_{jq})
  \end{array}
\end{equation}
$$
 (j=1,\ \dots,\ N;\ \mu=1,\ \dots,\ m;\ q=2,\ \dots,\ R_j;\ \nu=1,\ \dots,\ 2m).
$$
It is easy to see that problem~(\ref{eqQTheta})--(\ref{eqTTheta}) can be also obtained from problem~(\ref{eqQ})--(\ref{eqT})
by substituting $\theta$ for $D_z$.

Consider the space ${\cal H}_a^l(K_j)=\bigoplus\limits_{t=1}^{R_j}H_a^l(K_{jt})$ with the norm
$
 \|v_j\|_{{\cal H}_a^l(K_j)}=\left(\sum\limits_{t=1}^{R_j}\|v_{jt}\|_{H_a^l(K_{jt})}^2\right)^{1/2}
$ 
and the space ${\cal E}_a^l(K_j)=\bigoplus\limits_{t=1}^{R_j}E_a^l(K_{jt})$ with the norm
$
 \|v_j\|_{{\cal E}_a^l(K_j)}=\left(\sum\limits_{t=1}^{R_j}\|v_{jt}\|_{E_a^l(K_{jt})}^2\right)^{1/2}.
$

Introduce the spaces of vector--functions
$$
 {\cal E}_a^{l+2m,\,N}(K)=\prod_{j=1}^N {\cal E}_a^{l+2m}(K_j),\ 
 {\cal E}_a^{l,\,N}(K,\ \gamma)=\prod_{j=1}^N {\cal E}_a^l(K_j,\ \gamma_j),
$$
$$
 \begin{array}{c}
 {\cal E}_a^l(K_j,\ \gamma_j)={\cal E}_a^l(K_j)\times\\
   \times\prod_{\sigma=1,\,R_j+1}\prod_{\mu=1}^m E_a^{l+2m-m'_{j\sigma\mu}-1/2}(\gamma_{j\sigma})
   \times\prod_{q=2}^{R_j}\prod_{\nu=1}^{2m} E_a^{l+2m-\nu+1/2}(\gamma_{jq}).
 \end{array}
$$
We study solutions $v=(v_1,\ \dots,\ v_N)\in{\cal E}_a^{l+2m,\,N}(\Omega)$ for problem~(\ref{eqQTheta})--(\ref{eqTTheta}) 
supposing that $f=\{f_j,\ g_{j\sigma\mu},\ h_{jq\nu}\}\in {\cal E}_a^{l,\,N}(\Omega,\ \Gamma).$ Introduce the bounded operator
${\cal M}(\theta) : {\cal E}_a^{l+2m,\,N}(\Omega)\to {\cal E}_a^{l,\,N}(\Omega,\ \Gamma)$ 
corresponding to problem~(\ref{eqQTheta})--(\ref{eqTTheta}) and given by
$$
{\cal M}v=\{w_j,\ {\cal C}_{j\sigma\mu}(D_y,\ \theta)v,\ {\cal T}_{jq\nu}(D_y,\ \theta)v\}.
$$ 
Here ${\cal C}_{j\sigma\mu}(D_y,\ \theta)v$ and ${\cal T}_{jq\nu}(D_y,\ \theta)v$ are given by~(\ref{eqCTheta}) and~(\ref{eqTTheta}) 
correspondingly; $w_j(y)\equiv {\cal Q}_j(D_y,\ \theta)v_{jt}(y)$ for $y\in K_{jt}$.

Repeating the proof of Lemma~\ref{lHomog'}, from lemma~3.1$'$ \cite{MP}, we get the following statement.
\begin{lemma}\label{lHomogTheta'} 
 For any $g_{j\sigma\mu}\in E_a^{l+2m-m'_{j\sigma\mu}-1/2}(\gamma_{j\sigma}),$ 
 $h_{jq\nu}\in E_a^{l+2m-\nu+1/2}(\gamma_{jq})$
 $(j=1,\ \dots,\ N;\ \sigma=1,\ R_j+1;\ \mu=1,\ \dots,\ m;\ q=2,\ \dots,\ R_j;\ \nu=1,\ \dots,\ 2m)$ there exists
 a vector--function $v\in  {\cal E}_a^{l+2m,\,N}(\Omega)$ such that
  $$
  {\cal C}_{j\sigma\mu}(D_y,\ \theta)v=g_{j\sigma\mu}(y)\ (y\in\gamma_{j\sigma}),\quad
  {\cal T}_{jq\nu}(D_y,\ \theta)v=h_{jq\nu}(y)\ (y\in\gamma_{jq}),
  $$
 $$
 \begin{array}{c}
  \|v\|_{{\cal E}_a^{l+2m,\,N}(\Omega)}\le 
  c\sum\limits_j\Bigl\{\sum_{\sigma,\,\mu}\|g_{j\sigma\mu}\|
       _{E_a^{l+2m-m'_{j\sigma\mu}-1/2}(\gamma_{j\sigma})}+\\
 + \sum_{q,\,\nu}\|h_{jq\nu}\|_{E_a^{l+2m-\nu+1/2}(\gamma_{jq})}\Bigr\}, 
 \end{array}
 $$
 where $c>0$ is independent of $g_{j\sigma\mu},\ h_{jq\nu},\ \theta.$
\end{lemma}
\bigskip

\section{Solvability of nonlocal transmission problems in plane angles}\label{sectTransSolv}
The results of this section are analogous to those of~\S\ref{sectBoundSolv}. We shall need these results 
for obtaining a priori estimates of solutions to nonlocal transmission problems in dihedral angles in~\S\ref{sectTransApr}. 

\subsection{Reduction of nonlocal problems in plane angles to nonlocal problems on arcs.}
Consider the nonlocal transmission problem for a vector--function $V=(V_1,\ \dots,\ V_N)\in {\cal H}_a^{l+2m,\,N}(K)$
\begin{equation}\label{eqQ0}
  {\cal Q}_j(D_y,\ 0)V_{jt}=f_{jt}(y) \quad (y\in K_{jt};\ t=1,\ \dots,\ R_j),
\end{equation}
\begin{equation}\label{eqC0}
  \begin{array}{c}
  {\cal C}_{j1\mu}(D_y,\ 0)V=
      C_{j1\mu}(D_y,\ 0)V_{j1}(y)|_{\gamma_{j1}}=g_{j1\mu}(y) \quad (y\in\gamma_{j1}),\\
  {\cal C}_{j,R_j+1,\mu}(D_y,\ 0)V=
     C_{j,R_j+1,\mu}(D_y,\ 0)V_{jR_j}(y)|_{\gamma_{j,R_j+1}}=\\
  =g_{j,R_j+1,\mu}(y) 
      \quad (y\in\gamma_{j,R_j+1}),
  \end{array}
\end{equation}
\begin{equation}\label{eqT0}
  \begin{array}{c}
  {\cal T}_{jq\nu}(D_y,\ 0)V=
     T_{jq\nu}(D_y,\ 0)V_{j,q-1}(y)|_{\gamma_{jq}}-T_{jq\nu}(D_y,\ 0)V_{jq}(y)|_{\gamma_{jq}}+\\
   +\sum\limits_{k,\sigma,s}
      (T_{jq\nu k\sigma s}(D_y,\ 0)V_k)({\cal G}'_{jqk\sigma s}y)|_{\gamma_{jq}}=h_{jq\nu}(y) \quad (y\in\gamma_{jq})
  \end{array}
\end{equation}
$$
 (j=1,\ \dots,\ N;\ \mu=1,\ \dots,\ m;\ q=2,\ \dots,\ R_j;\ \nu=1,\ \dots,\ 2m),
$$
where $f=\{f_j,\ g_{j\sigma\mu},\ h_{jq\nu}\}\in  {\cal H}_a^{l,\,N}(K,\ \gamma).$

Put formally $D_z=0$ and write the differential operators in the polar coordinates: 
${\cal Q}_j(D_y,\ 0)=r^{-2m}\tilde{\cal Q}_j(\varphi,\ D_\varphi,\ rD_r)$,
$C_{j\sigma\mu}(D_y,\ 0)=r^{-m'_{j\sigma\mu}}\tilde C_{j\sigma\mu}(\varphi,\ D_\varphi,\ rD_r)$,
$T_{jq\nu}(D_y,\ 0)=r^{-\nu+1}\tilde T_{jq\nu}(\varphi,\ D_\varphi,\ rD_r),$
$T_{jq\nu k\sigma s}(D_y,\ 0)=r^{-\nu+1}\tilde T_{jq\nu k\sigma s}(\varphi,\ D_\varphi,\ rD_r)$.

Put $\tau=\ln r$ and do the Fourier transform with respect to~$\tau$; then from~(\ref{eqQ0})--(\ref{eqT0}), we get
\begin{equation}\label{eqQLambda}
  \tilde{\cal Q}_j(\varphi,\ D_\varphi,\ \lambda)\tilde V_{jt}(\varphi,\ \lambda)=\tilde F_{jt}(\varphi,\ \lambda) \quad 
          (\varphi\in(b_{jt},\ b_{j,t+1});\ t=1,\ \dots,\ R_j),
\end{equation}
\begin{equation}\label{eqCLambda}
  \begin{array}{c}
     \tilde {\cal C}_{j1\mu}(\varphi,\ D_\varphi,\ \lambda)\tilde V(\varphi,\ \lambda)=  
     \tilde C_{j1\mu}(\varphi,\ D_\varphi,\ \lambda)\tilde V_{j1}(\varphi,\ \lambda)|_{\varphi=b_{j1}}=\\
             =\tilde G_{j1\mu}(\lambda),\\
     \tilde {\cal C}_{j,R_j+1,\mu}(\varphi,\ D_\varphi,\ \lambda)\tilde V(\varphi,\ \lambda)=
     \tilde C_{j,R_j+1,\mu}(\varphi,\ D_\varphi,\ \lambda)\tilde V_{jR_j}(\varphi,\ \lambda)|_
          {\varphi=b_{j,R_j+1}}=\\
  =\tilde G_{j,R_j+1,\mu}(\lambda),
  \end{array}
\end{equation}
\begin{equation}\label{eqTLambda}
  \begin{array}{c}
     \tilde {\cal T}_{jq\nu}(\varphi,\ D_\varphi,\ \lambda)\tilde V(\varphi,\ \lambda)=
     \tilde T_{jq\nu}(\varphi,\ D_\varphi,\ \lambda)\tilde V_{j,q-1}(\varphi,\ \lambda)|_{\varphi=b_{jq}}-\\
           -\tilde T_{jq\nu}(\varphi,\ D_\varphi,\ \lambda)\tilde V_{jq}(\varphi,\ \lambda)|_{\varphi=b_{jq}}+\\
   +\sum\limits_{k,\sigma,s}
       e^{(i\lambda-(\nu-1))\ln\chi'_{jqk\sigma s}}\tilde T_{jq\nu k\sigma s}(\varphi,\ D_\varphi,\ \lambda)
            \tilde V_k(\varphi+\varphi'_{jqk\sigma},\ \lambda)|_{\varphi=b_{jq}}
   =\tilde H_{jq\nu}(\lambda)
  \end{array}
\end{equation}
$$
 (j=1,\ \dots,\ N;\ \mu=1,\ \dots,\ m;\ q=2,\ \dots,\ R_j;\ \nu=1,\ \dots,\ 2m).
$$
Here $F_{jt}(\varphi,\ \tau)=e^{2m\tau}f_{jt}(\varphi,\ \tau)$, 
$G_{j\sigma\mu}(\tau)=e^{m'_{j\sigma\mu}\tau}g_{j\sigma\mu}(\tau)$;
$H_{jq\nu}(\tau)=e^{(\nu-1)\tau}h_{jq\nu}(\tau)$; $\tilde V_{jt}$, $\tilde F_{jt}$, $\tilde G_{j\sigma\mu}$, and $\tilde H_{jq\nu}$ are
the Fourier transforms of $V_{jt}$, $F_{jt}$, $G_{j\sigma\mu}$, and $H_{jq\nu}$ correspondingly.

This problem is a system of $R_1+\dots+ R_N$ ordinary differential equations~(\ref{eqQLambda}) for the functions
$\tilde V_{jt}\in W^{l+2m}(b_{jt},\ b_{j,t+1})$ with boundary conditions~(\ref{eqCLambda}) and nonlocal transmission 
conditions~(\ref{eqTLambda}) connecting jumps of the functions $\tilde V_j$ and their derivatives at the points of the intervals
$(b_{j1},\ b_{j,R_j+1})$ with values of the functions $\tilde V_{k1}$ and $\tilde V_{k,R_k+1}$ and their derivatives
at the points $\varphi=b_{k1}$ and $\varphi=b_{k,R_k+1}$ correspondingly.

Notice that problem~(\ref{eqQLambda})--(\ref{eqTLambda}) is formally adjoint to problem~(\ref{eqPLambda}), (\ref{eqBLambda})
with respect to Green formula~(\ref{eqGrPLambda}).

\subsection{Solvability of nonlocal problems with parameter $\lambda$ on arcs.}
Consider the space ${\cal W}^l(b_{j1},\ b_{j,R_j+1})=\bigoplus\limits_{t=1}^{R_j}W^l(b_{jt},\ b_{j,t+1})$ with the norm
$
 \|\tilde V_j\|_{{\cal W}^l(b_{j1},\ b_{j,R_j+1})}=
   \left(\sum\limits_{t=1}^{R_j}\|\tilde V_{jt}\|_{W^l(b_{jt},\ b_{j,t+1})}^2\right)^{1/2}.
$
Introduce the spaces of vector--functions
$$
 {\cal W}^{l+2m,\, N}(b_1,\ b_2)=\prod_{j=1}^N {\cal W}^{l+2m}(b_{j1},\ b_{j,R_j+1}),
$$
$$
 {\cal W}^{l,\, N}[b_1,\ b_2]=\prod_{j=1}^N {\cal W}^{l}[b_{j1},\ b_{j,R_j+1}], 
$$
$$
  {\cal W}^{l}[b_{j1},\ b_{j,R_j+1}]={\cal W}^l(b_{j1},\ b_{j,R_j+1}) \times{\mathbb C}^m \times{\mathbb C}^m\times
  \prod_{q=2}^{R_j}{\mathbb C}^{2m}.
$$
Introduce the equivalent norms depending on the parameter  $\lambda$ ($|\lambda|\ge1$) in the Hilbert spaces
${\cal W}^{l}(b_{j1},\ b_{j,R_j+1})$ and ${\cal W}^{l}[b_{j1},\ b_{j,R_j+1}]$:
$$
 |||\tilde V_j|||_{{\cal W}^{l}(b_{j1},\ b_{j,R_j+1})}=
 \bigl(\|\tilde V_j\|_{{\cal W}^{l}(b_{j1},\ b_{j,R_j+1})}^2+
    |\lambda|^{2l}\,\|\tilde V_j\|_{L_2(b_{j1},\ b_{j,R_j+1})}^2\bigr)^{1/2},
$$
$$
 \begin{array}{c}
 |||\{\tilde F_j,\ \tilde G_{j\sigma\mu},\ \tilde H_{jq\nu}\}|||_{{\cal W}^{l}[b_{j1},\ b_{j,R_j+1}]}=
 \Bigl(|||\tilde F_j|||_{{\cal W}^l(b_{j1},\ b_{j,R_j+1})}^2+\\
 +\sum\limits_{\sigma,\,\mu}(1+|\lambda|^{2(l+2m-m'_{j\sigma\mu}-1/2)}) |\tilde G_{j\sigma\mu}|^2+
 \sum\limits_{q,\,\nu}(1+|\lambda|^{2(l+2m-\nu+1/2)})|\tilde H_{jq\nu}|^2\Bigr)^{1/2},
 \end{array}
$$
where $\tilde V_j\in {\cal W}^l(b_{j1},\ b_{j,R_j+1}),\ \{\tilde F_j,\ \tilde G_{j\sigma\mu},\ \tilde H_{jq\nu}\}\in 
  {\cal W}^l[b_{j1},\ b_{j,R_j+1}].$ And therefore we have
$$
  |||\tilde V|||_{W^{l+2m,\,N}(b_1,\ b_2)}=
  \Bigl(\sum_j |||\tilde V_j|||^2_{{\cal W}^{l+2m}(b_{j1},\ b_{j,R_j+1})}\Bigr)^{1/2},
$$
$$
  |||\tilde \Phi|||_{{\cal W}^{l,\,N}[b_1,\ b_2]}=
  \Bigl(\sum_j |||\tilde \Phi_j|||^2_{{\cal W}^l[b_{j1},\ b_{j,R_j+1}]}\Bigr)^{1/2},
$$
where $\tilde V=(\tilde V_1,\ \dots,\ \tilde V_N)\in {\cal W}^{l+2m,\,N}(b_1,\ b_2),$
$\tilde \Phi=(\tilde \Phi_1,\ \dots,\ \tilde\Phi_N)\in {\cal W}^{l,\,N}[b_1,\ b_2].$

Consider the operator--valued function $\tilde{\cal M}(\lambda): {\cal W}^{l+2m,\, N}(b_1,\ b_2)\to {\cal W}^{l,\, N}[b_1,\ b_2]$
corresponding to problem~(\ref{eqQLambda})--(\ref{eqTLambda}) and given by
$$
\tilde{\cal M}(\lambda)\tilde V=\{\tilde W_j,\ 
  \tilde{\cal C}_{j\sigma\mu}(\varphi,\ D_\varphi,\ \lambda)\tilde V,\ 
  \tilde{\cal T}_{jq\nu}(\varphi,\ D_\varphi,\ \lambda)\tilde V\}.
$$
Here $\tilde{\cal C}_{j\sigma\mu}(\varphi,\ D_\varphi,\ \lambda)\tilde V$ and 
$\tilde{\cal T}_{jq\nu}(\varphi,\ D_\varphi,\ \lambda)\tilde V$ are given by~(\ref{eqCLambda}) and~(\ref{eqTLambda}) 
correspondingly; $\tilde W_j(\varphi)=\tilde{\cal Q}_j(\varphi,\ D_\varphi,\ \lambda)\tilde V_{jt}(\varphi)$ for $\varphi\in(b_{jt},\ b_{j,t+1})$.

\begin{lemma}\label{lSolvLambda1'} 
 For all $\lambda\in{\mathbb C}$, the operator
 $\tilde{\cal M}(\lambda):{\cal W}^{l+2m,\, N}(b_1,\ b_2)\to {\cal W}^{l,\, N}[b_1,\ b_2]$ is Fredholm,
 $\ind\tilde{\cal M}(\lambda)=0;$ for any $h\in\mathbb R$, there exists a $q_0>0$ such that for
 $\lambda\in J_{h,\,q_0}=\{\lambda\in{\mathbb C}:\ \Im\,\lambda=h,\ |\Re\,\lambda|\ge q_0\}$, the operator
 $\tilde{\cal M}(\lambda)$ has the bounded inverse 
 $\tilde{\cal M}^{-1}(\lambda):{\cal W}^{l,\, N}[b_1,\ b_2]\to {\cal W}^{l+2m,\, N}(b_1,\ b_2)$ and
 \begin{equation}\label{eqSolvLambda1'}
  |||\tilde{\cal M}^{-1}(\lambda)\tilde \Phi|||_{{\cal W}^{l+2m,\, N}(b_1,\ b_2)}\le 
          c|||\tilde \Phi|||_{{\cal W}^{l,\, N}[b_1,\ b_2]}
 \end{equation}
 for all $\tilde \Phi\in {\cal W}^{l,\, N}[b_1,\ b_2],$ where $c>0$ is independent of $\lambda$ and $\tilde\Phi;$ 
 the operator--valued function
 $\tilde{\cal M}^{-1}(\lambda):{\cal W}^{l,\, N}[b_1,\ b_2]\to {\cal W}^{l+2m,\, N}(b_1,\ b_2)$ 
 is finitely meromorphic.
\end{lemma}
\begin{proof}
If
$$
 \begin{array}{c}
\tilde{\cal T}_{jq\nu}(\varphi,\ D_\varphi,\ \lambda)\tilde V(\varphi,\ \lambda)=
 \tilde T_{jq\nu}(\varphi,\ D_\varphi,\ \lambda)\tilde V_{j,q-1}(\varphi,\ \lambda)|_{\varphi=b_{jq}}-\\
           -\tilde T_{jq\nu}(\varphi,\ D_\varphi,\ \lambda)\tilde V_{jq}(\varphi,\ \lambda)|_{\varphi=b_{jq}}
 \end{array}
$$
(i.e., if the operators $T_{jq\nu k\sigma s}(\varphi,\ D_\varphi,\ rD_r)$ corresponding to the nonlocal terms are absent), then 
we denote by $\tilde{\cal M}_0(\lambda)$ the operator $\tilde{\cal M}(\lambda)$.
Following the scheme developed by M.S.~Agranovich and M.I.~Vishik in~\cite{AV} (see also~\cite[\S5]{RSh}), 
one can show that there exist $0<\varepsilon_1<\pi/2$ and $q_1>0$ such that for
$$
 \lambda\in Q_{\varepsilon_1,q_1}=\{\lambda: |\lambda|\ge q_1,\ |\arg\lambda|\le\varepsilon_1\}\cup
 \{\lambda: |\lambda|\ge q_1,\ |\arg\lambda-\pi|\le\varepsilon_1\},
$$
there exists the bounded inverse operator $\tilde{\cal M}^{-1}_0(\lambda);$ moreover, for all $\tilde \Phi\in W^{l,\, N}[b_1,\ b_2]$,
\begin{equation}\label{eqSolvLambda11'}
  |||\tilde{\cal M}_0^{-1}(\lambda)\tilde \Phi|||_{{\cal W}^{l+2m,\, N}(b_1,\ b_2)}\le 
      k_1|||\tilde \Phi|||_{{\cal W}^{l,\, N}[b_1,\ b_2]}.
\end{equation}
Here $k_1>0$ is independent of $\lambda$ and $\tilde\Phi.$ 

Consider the operator
$\tilde{\cal M}_t(\lambda)=\tilde{\cal M}_0(\lambda)+t(\tilde{\cal M}(\lambda)-\tilde{\cal M}_0(\lambda)),$
$0\le t\le 1.$ We shall prove that for any $h\in{\mathbb R}$, there exists a $q_0>0$ such that if
$\lambda\in J_{h,q_0}$ and $0\le t\le 1$, then we have 
\begin{equation}\label{eqSolvLambda12'}
 k_2|||\tilde{\cal M}_t(\lambda)\tilde V|||_{{\cal W}^{l,\, N}[b_1,\ b_2]}\le
 |||\tilde V|||_{{\cal W}^{l+2m,\, N}(b_1,\ b_2)}\le
 k_3|||\tilde{\cal M}_t(\lambda)\tilde V|||_{{\cal W}^{l,\, N}[b_1,\ b_2]}
\end{equation}
for all $\tilde V\in {\cal W}^{l+2m,\, N}(b_1,\ b_2)$. Here $k_2,\ k_3>0$ are independent of $\lambda,\ t$ and $V.$

Denote $\tilde{\cal M}_t(\lambda)\tilde V=\tilde\Phi$; then we have
$$
 \tilde{\cal M}_0(\lambda)\tilde V=\tilde\Phi+\tilde\Psi,
$$
where
$$
  \tilde\Psi=\{0,\ 0,\ -t\sum\limits_{k,\,\sigma,\,s}
     e^{(i\lambda-(\nu-1))\ln\chi'_{jqk\sigma s}}\tilde T_{jq\nu k\sigma s}(\varphi,\ D_\varphi,\ \lambda)
            \tilde V_k(\varphi+\varphi'_{jqk\sigma},\ \lambda)|_{\varphi=b_{jq}}\}.
$$
By virtue of (\ref{eqSolvLambda11'}), we have
\begin{equation}\label{eqSolvLambda13'}
  |||\tilde V|||_{{\cal W}^{l+2m,\, N}(b_1,\ b_2)}\le
 k_1|||\tilde\Phi+\tilde\Psi|||_{{\cal W}^{l,\, N}[b_1,\ b_2]}.
\end{equation}
Take $\varepsilon>0$ from formula~(\ref{eqEpsilon}) and a $q_0\ge q_1$ such that
$J_{h,q_0}\subset Q_{\varepsilon_1,q_1}.$ Then using inequalities~(\ref{eqInterp1}), (\ref{eqInterp2}), we get
\begin{equation}\label{eqSolvLambda14'}
 \begin{array}{c}
  I_{jq\nu k1s}=(1+|\lambda|^{l+2m-\nu+1/2})
     \Big| e^{(i\lambda-(\nu-1))\ln\chi'_{jqk1s}}\times\\
  \times\tilde T_{jq\nu k1s}(\varphi,\ D_\varphi,\ \lambda)
            \tilde V_k(\varphi+\varphi'_{jqk1})|_{\varphi=b_{jq}}\Big|\le\\
  k_4|\lambda|^{l+2m-\nu}\bigl\{\|\tilde T_{jq\nu k1s}(\varphi,\ D_\varphi,\ \lambda)
            \tilde V_{k1}\|_
           {W^1(b_{k1},\ b_{k1}+\varepsilon/2)}+\\
  |\lambda|\,\|\tilde T_{jq\nu k1s}(\varphi,\ D_\varphi,\ \lambda)
            \tilde V_{k1}\|_
           {L^2(b_{k1},\ b_{k1}+\varepsilon/2)}\bigr\}\le 
            k_5 |||\tilde V_{k1}|||_{W^{l+2m}(b_{k1},\ b_{k1}+\varepsilon/2)}.
  \end{array}
\end{equation}
If $\varepsilon_1$ is sufficiently small and $q_1$ is sufficiently large, 
then from inequality~(\ref{eqSolvLambda14'}), theorem~4.1 \cite[Chapter 1, \S4]{AV}, 
Leibniz' formula, and interpolation inequality~(\ref{eqInterp1}), we obtain
\begin{equation}\label{eqSolvLambda15'}
 \begin{array}{c}
  I_{jq\nu k1s}\le k_5 |||\zeta_{k1}\tilde V_{k1}|||_{W^{l+2m}(b_{k1},\ b_{k1}+\varepsilon/2)}\le
    k_6\bigl( |||{\tilde {\cal Q}_k}(\zeta_{k1}\tilde V_{k1})|||_{W^l(b_{k1},\ b_{k2})}+\\
    +\sum\limits_{\mu=1}^m (1+|\lambda|^{l+2m-m'_{k1\mu}-1/2})\Big|\tilde C_{k1\mu}(\varphi,\ D_\varphi,\ \lambda)
            \tilde V_{k1}(\varphi)|_{\varphi=b_{k1}}\Big|\bigr)\le \\
    \le  k_7\bigl( |||{\tilde {\cal Q}_k}\tilde V_{k1}|||_{W^l(b_{k1},\ b_{k2})}
    +|\lambda|^{-1}|||\tilde V_{k1}|||_{W^{l+2m}(b_{k1},\ b_{k2})}+\\
          +\sum\limits_{\mu=1}^m (1+|\lambda|^{l+2m-m'_{k1\mu}-1/2})|\tilde C_{k1\mu}(\varphi,\ D_\varphi,\ \lambda)
            \tilde V_{k1}(\varphi)|_{\varphi=b_{k1}}|\bigr).
 \end{array}
\end{equation}
Similarly to~(\ref{eqSolvLambda14'}), (\ref{eqSolvLambda15'}), one can estimate
$$
 \begin{array}{c}
 I_{jq\nu k,R_k+1,s}=(1+|\lambda|^{l+2m-\nu+1/2})\times\\
     \times\Big| e^{(i\lambda-(\nu-1))\ln\chi'_{jqk,R_k+1,s}}\tilde T_{jq\nu k,R_k+1,s}(\varphi,\ D_\varphi,\ \lambda)
            \tilde V_k(\varphi+\varphi'_{jqk,R_k+1})|_{\varphi=b_{jq}}\Big|:
 \end{array}
$$
\begin{equation}\label{eqSolvLambda16'}
 \begin{array}{c}
 I_{jq\nu k,R_k+1,s}\le k_8 \bigl( |||{\tilde {\cal Q}_k}\tilde V_{kR_k}|||_{W^l(b_{kR_k},\ b_{k,R_k+1})}+\\
    +|\lambda|^{-1}|||\tilde V_{kR_k}|||_{W^{l+2m}(b_{kR_k},\ b_{kR_k+1})}+\\
     +\sum\limits_{\mu=1}^m (1+|\lambda|^{l+2m-m'_{k,R_k+1,\mu}-1/2})
      |\tilde C_{k,R_k+1,\mu}(\varphi,\ D_\varphi,\ \lambda)\tilde V_{kR_k}(\varphi)|_{\varphi=b_{k,R_k+1}}|\bigr).
 \end{array}
\end{equation}
Now if $q_0$ is sufficiently large, then~(\ref{eqSolvLambda13'}),
(\ref{eqSolvLambda15'}), and (\ref{eqSolvLambda16'}) imply right-hand side of inequality~(\ref {eqSolvLambda12'}). 
Left-hand side of inequality~(\ref {eqSolvLambda12'}) is obvious.
Using a standard method of continuation with respect to parameter $t$ (see the proof of theorem~7.1 
\cite[Chapter 2, \S7]{Lad}), inequality~(\ref {eqSolvLambda12'}) and existence of a bounded inverse operator
$\tilde{\cal M}^{-1}_0(\lambda)$ for $\lambda\in Q_{\varepsilon_1,q_1}$, one can easily see that for
$\lambda\in J_{h,q_0}$, the operator $\tilde{\cal M}(\lambda)$ also has a bounded inverse
and~(\ref{eqSolvLambda1'}) holds.

Let us prove that the operator $\tilde{\cal M}(\lambda)$ is Fredholm.
For $\lambda_0\in Q_{\varepsilon_1,q_1}$, we have
$$
 \tilde{\cal M}(\lambda)\tilde{\cal M}_0^{-1}(\lambda_0)=
 I+(\tilde{\cal M}(\lambda)-\tilde{\cal M}_0(\lambda_0))\tilde{\cal M}_0^{-1}(\lambda_0),
$$
where $I$ is the identity operator in ${\cal W}^{l,\, N}[b_1,\ b_2].$ Since the 
operators $\tilde{\cal Q}_j(\varphi,\ D_\varphi,\ \lambda)$ contain the parameter $\lambda$ only in junior terms, the operator
$$
 \tilde{\cal M}(\lambda)-\tilde{\cal M}_0(\lambda_0): {\cal W}^{l+2m,\, N}(b_1,\ b_2)\to
 {\cal W}^{l+1,\, N}[b_1,\ b_2]
$$
is bounded for every fixed $\lambda\in{\mathbb C}$. 
Hence from the compactness of the imbedding
operator of $W^{l+1}(b_{jt},\ b_{j,t+1})$ into
$W^{l}(b_{jt},\ b_{j,t+1}),$ it follows that the operator
$$
 (\tilde{\cal M}(\lambda)-\tilde{\cal M}_0(\lambda_0))\tilde{\cal M}_0^{-1}(\lambda_0): 
   {\cal W}^{l,\, N}[b_1,\ b_2]\to {\cal W}^{l,\, N}[b_1,\ b_2]
$$
is compact. Thus by theorem~15.1 \cite[\S15]{Kr}, the operator $\tilde{\cal M}(\lambda)$ is Fredholm and
$\ind\tilde{\cal M}(\lambda)=0$ for all $\lambda\in{\mathbb C}$.

From this, from existence of the bounded inverse operator $\tilde{\cal M}^{-1}(\lambda)$ for $\lambda\in J_{h,q_0}$, and from
theorem~1 \cite{Bl}, it follows that the operator--valued function $\tilde{\cal M}^{-1}(\lambda)$ is finitely meromorphic.
\end{proof}

Repeating the proof of lemma~2.2 \cite[\S2]{SkDu90}, from~(\ref{eqSolvLambda13'})--(\ref{eqSolvLambda16'}), we obtain
the following statement.
\begin{lemma}\label{lSolvLambda2'}
 For any $0<\varepsilon'<1/d'$, there exists a $q>1$ such that the set
 $\{\lambda\in{\mathbb C}:\ |\Im\lambda|\le \varepsilon'\ln|\Re\lambda|,\ |\Re\lambda|\ge q\}$ contains no poles 
of the operator--valued function $\tilde{\cal M}^{-1}(\lambda),$ where $d'=\max\,|\ln\,\chi'_{jqk\sigma s}|;$ for every pole
 $\lambda_0$ of the operator--valued function $\tilde{\cal M}^{-1}(\lambda)$, there exists
 a $\delta>0$ such that the set $\{\lambda\in{\mathbb C}:\ 0<|\Im\lambda - \Im\lambda_0|<\delta\}$ contains no poles of
 the operator--valued function $\tilde{\cal M}^{-1}(\lambda).$
\end{lemma}

\subsection{One--valued solvability of nonlocal problems in plane angles.}
Replacing in the proof of theorem~2.1 \cite[\S2]{SkDu90}  Sobolev spaces $W^l(\cdot)$ by 
${\cal W}^l(\cdot)$ and weighted spaces $H_a^l(\cdot)$ by ${\cal H}_a^l(\cdot),$ from
Lemma~\ref{lSolvLambda1'}, we obtain the following result.
\begin{theorem}\label{thSolvFlat'} 
 Suppose the line $\Im\,\lambda=a+1-l-2m$ contains no poles of the operator--valued function 
 $\tilde{\cal M}^{-1}(\lambda)$; then nonlocal transmission problem~(\ref{eqQ0})--(\ref{eqT0}) 
 has a unique solution $V\in {\cal H}_a^{l+2m,\,N}(K)$ for every
 right-hand side $f\in {\cal H}_a^{l,\,N}(K,\ \gamma)$ and
 $$
  \|V\|_{{\cal H}_a^{l+2m,\,N}(K)}\le c\|f\|_{{\cal H}_a^{l,\,N}(K,\ \gamma)},
 $$
 where $c>0$ does not depend on $f.$
\end{theorem}
\section{A priori estimates of solutions for nonlocal transmission problems}\label{sectTransApr}
In this section, we prove a priori estimates for solutions to nonlocal transmission
problems analogous to those of~\S\ref{sectBoundApr}.

\subsection{A priori estimates in dihedral angles.}
Denote $d'_1=\min\{1,\ \chi'_{jqk\sigma s}\}/2,$ $d'_2=2 \max\{1,\ \chi'_{jqk\sigma s}\},$
$
 \Omega_{j}^p=\Omega_j\cap
     \{r_1(d'_1)^{6-p}<r<r_2(d'_2)^{6-p},\ |z|<2^{-p-1}\},
$ 
$
 \Omega_{jt}^p=\Omega_{jt}\cap
     \{r_1(d'_1)^{6-p}<r<r_2(d'_2)^{6-p},\ |z|<2^{-p-1}\},
$
where $j=1,\ \dots,\ N;\ t=1,\ \dots,\ R_j;\ p=0,\ \dots,\ 6;\ 0<r_1<r_2.$

Introduce the space ${\cal W}^l(\Omega_j^p)=\bigoplus\limits_{t=1}^{R_j}W^l(\Omega_{jt}^p)$ with the norm
$
 \|V_j\|_{{\cal W}^l(\Omega_j^p)}=\left(\sum\limits_{t=1}^{R_j}\|V_{jt}\|_{W^l(\Omega_{jt}^p)}^2\right)^{1/2}.
$

\begin{lemma}\label{lAprW'} 
 Suppose $V_j\in {\cal W}^{2m}(\Omega_j^0),$
 \begin{equation}\label{eqMUl}
  \begin{array}{c}
   {\cal Q}_j(D_y,\ D_z)V_{jt}\in W^l(\Omega_{jt}^0),\\
   {\cal C}_{j\sigma\mu}(D_y,\ D_z)V\in 
      W^{l+2m-m'_{j\sigma\mu}-1/2}(\Gamma_{j\sigma}\cap\bar\Omega_j^0),\\
   {\cal T}_{jq\nu}(D_y,\ D_z)V\in 
      W^{l+2m-\nu+1/2}(\Gamma_{jq}\cap\bar\Omega_j^0)
  \end{array}
 \end{equation}
 $$
 \begin{array}{c}
 ( j=1,\ \dots,\ N;\ \sigma=1,\ R_j+1;\ \mu=1,\ \dots,\ m;\\
  q=2,\ \dots,\ R_j;\ \nu=1,\ \dots,\ 2m);
 \end{array}
 $$
 then we have $V\in \prod\limits_j {\cal W}^{l+2m}(\Omega_j^3)$ and for $|\lambda|\ge 1$,
 \begin{equation}\label{eqAprW'}
  \begin{array}{c}
   \sum\limits_j \|V_j\|_{{\cal W}^{l+2m}(\Omega_j^6)}\le
           c \sum\limits_j \bigl\{\sum\limits_t\|{\cal Q}_j(D_y,\ D_z)V_{jt}\|_{W^l(\Omega_{jt}^3)}+\\
   +\sum\limits_{\sigma,\,\mu}\|{\cal C}_{j\sigma\mu}(D_y,\ D_z)V\|
           _{W^{l+2m-m'_{j\sigma\mu}-1/2}(\Gamma_{j\sigma}\cap\bar \Omega_j^3)}+\\
   +\sum\limits_{q,\,\nu}\|{\cal T}_{jq\nu}(D_y,\ D_z)V\|
           _{W^{l+2m-\nu+1/2}(\Gamma_{jq}\cap\bar \Omega_j^3)}+\\
        + |\lambda|^{-1}\|V_j\|_{{\cal W}^{l+2m}(\Omega_j^3)}+
             |\lambda|^{l+2m-1}\|V_j\|_{L_2(\Omega_j^3)}\bigr\},
  \end{array}
 \end{equation}
 where $c>0$ is independent of $\lambda$ and $V.$
\end{lemma}
\begin{proof}
Since the functions $\zeta_{jq}$ $(q=1,\ \dots,\ R_j+1)$ given by~(\ref{eqZeta}) are the multiplicators in the spaces
$W^l(\Omega_{jt}^p)$ $(t=1,\ \dots,\ R_j),$ we have
$\zeta_{j\sigma}V_{j}\in W^{2m}(\Omega_j^0)$ $(\sigma=1,\ R_j+1).$
Apply theorem~5.1 \cite[Chapter 2, \S 5.1]{LM} to the functions $\zeta_{j\sigma}V_{j}$ and to the operator
$\{{\cal Q}_j(D_y,\ D_z),\ {\cal C}_{j\sigma\mu}(D_y,\ D_z)\}$; then from~(\ref{eqMUl}) and Leibniz' formula, we get
\begin{equation}\label{eqAprW2'}
  \zeta_{j\sigma}V_{j}\in\ W^{l+2m}(\Omega_j^1).
\end{equation}
Denote $W_{jq\nu}=
  \sum\limits_{k,\,\sigma,\,s}(T_{jq\nu k\sigma s}(D_y,\ D_z)(\zeta_{k\sigma}V_k))({\cal G}'_{jqk\sigma s}y,\ z).$
Clearly,
\begin{equation}\label{eqAprW3'}
 W_{jq\nu}|_{\Gamma_{jq}\cap\bar\Omega_j^2}=
  \sum\limits_{k,\ \sigma,\,s}(T_{jq\nu k\sigma s}(D_y,\ D_z)V_k))({\cal G}'_{jqk\sigma s}y,\ z)|
         _{\Gamma_{jq}\cap\bar\Omega_j^2}.
\end{equation}
From equality~(\ref{eqAprW3'}) and relations~(\ref{eqMUl}), (\ref{eqAprW2'}), it follows that
\begin{equation}\label{eqAprW4'}
 \begin{array}{c}
  T_{jq\nu}(D_y,\ D_z)V_{j,q-1}|_{\Gamma_{jq}\cap\bar\Omega_j^2}
      -T_{jq\nu}(D_y,\ D_z)V_{jq}|_{\Gamma_{jq}\cap\bar\Omega_j^2}=\\
  ={\cal T}_{jq\nu}(D_y,\ D_z)V - W_{jq\nu}|_{\Gamma_{jq}\cap\bar\Omega_j^2}
   \in W^{l+2m-\nu+1/2}(\Gamma_{jq}\cap\bar\Omega_j^2).
 \end{array}
\end{equation}
Now~(\ref{eqMUl}), (\ref{eqAprW4'}), and theorem~1 \cite[\S2]{Sh65} imply that
$V_j\in {\cal W}^{l+2m}(\Omega_j^3)$ and
\begin{equation}\label{eqAprW5'}
 \begin{array}{c}
    \sum\limits_j \|V_j\|_{{\cal W}^{l+2m}(\Omega_j^6)}
           \le k_1 \sum\limits_j \bigl\{\sum\limits_t\|{\cal Q}_j(D_y,\ D_z)V_{jt}\|_{W^l(\Omega_{jt}^5)}+\\
   +\sum\limits_{\sigma,\,\mu}\|{\cal C}_{j\sigma\mu}(D_y,\ D_z)V\|
           _{W^{l+2m-m'_{j\sigma\mu}-1/2}(\Gamma_{j\sigma}\cap\bar \Omega_j^5)}+\\
   +\sum\limits_{q,\,\nu}\|T_{jq\nu}(D_y,\ D_z)V_{j,q-1}|_{\Gamma_{jq}\cap\bar \Omega_j^5}-\\
     -T_{jq\nu}(D_y,\ D_z)V_{jq}|_{\Gamma_{jq}\cap\bar \Omega_j^5}\|
           _{W^{l+2m-\nu+1/2}(\Gamma_{jq}\cap\bar \Omega_j^5)}+ \|V_j\|_{L_2(\Omega_j^5)}
    \bigr\}.
 \end{array}
\end{equation}
Again using theorem~5.1 \cite[Chapter 2, \S 5.1]{LM}, Leibniz' formula, and inequality~(\ref{eqInterp1}), we get
\begin{equation}\label{eqAprW6'}
 \begin{array}{c}
  \|W_{jq\nu}|_{\Gamma_{jq}\cap\bar\Omega_j^5}\|_{W^{l+2m-\nu+1/2}(\Gamma_{jq}\cap\bar\Omega_j^5)}\le
  k_2\sum\limits_{k,\,\sigma}\|\zeta_{k\sigma}V_k\|_{W^{l+2m}(\Omega_k^4)}\le\\
 \le k_3\sum\limits_k\bigl\{\sum\limits_t\|{\cal Q}_k(D_y,\ D_z)V_{kt}\|_{W^l(\Omega_{kt}^3)}+\\
   +\sum\limits_{\sigma,\,\mu}\|{\cal C}_{k\sigma\mu}(D_y,\ D_z)V\|
           _{W^{l+2m-m'_{k\sigma\mu}-1/2}(\Gamma_{k\sigma}\cap\bar \Omega_j^3)}+\\
    +|\lambda|^{-1}\|V_k\|_{{\cal W}^{l+2m}(\Omega_k^3)}+
             |\lambda|^{l+2m-1}\|V_k\|_{L_2(\Omega_k^3)}\bigr\}.
 \end{array}
\end{equation}
From~(\ref{eqAprW5'}), (\ref{eqAprW3'}), and (\ref{eqAprW6'}), it follows inequality~(\ref{eqAprW'}).
\end{proof}

Denote ${\cal W}_\loc^l(\bar\Omega_j\backslash M)=\bigoplus\limits_{t=1}^{R_j}W_\loc^l(\bar\Omega_{jt}\backslash M).$

\begin{theorem}\label{thAprH'} 
 Let $V\in\prod\limits_j {\cal W}_\loc^{2m}(\bar\Omega_j\backslash M)$ be a solution for nonlocal transmission 
problem~(\ref{eqQ})--(\ref{eqT}) such that $V\in H_{a-l-2m}^{0,\,N}(\Omega)$ and 
 $f\in {\cal H}_a^{l,\,N}(\Omega,\ \Gamma)$; then $V\in {\cal H}_a^{l+2m,\,N}(\Omega)$ and
 \begin{equation}\label{eqAprH'}
  \|V\|_{{\cal H}_a^{l+2m,\,N}(\Omega)}\le c \bigl(\|f\|_{{\cal H}_a^{l,\,N}(\Omega,\ \Gamma)}+
      \|V\|_{H_{a-l-2m}^{0,\,N}(\Omega)}\bigr),
 \end{equation}
 where $c>0$ is independent of $V.$
\end{theorem}
\begin{proof}
From Lemma~\ref{lAprW'}, it follows that $V\in\prod\limits_j {\cal W}_\loc^{l+2m}(\bar\Omega_j\backslash M).$
Now repeating the proof of lemma~3.2 \cite[\S3]{SkDu90} and 
replacing there $W^l(\cdot)$ by ${\cal W}^l(\cdot)$ and weighted spaces $H_a^l(\cdot)$ by ${\cal H}_a^l(\cdot)$, from
Lemmas~\ref{lHomog'} and~\ref{lAprW'}, we derive that $V\in{\cal H}_a^{l+2m,\,N}(\Omega)$ and a priori 
estimate~(\ref{eqAprH'}) holds.
\end{proof}

\subsection{A priori estimates in plane angles.}
Put $K_j^{ps}=K_j\cap\{r_1(d'_1)^{6-p}\cdot 2^s<r<r_2(d'_2)^{6-p}\cdot 2^s\},$
$K_{jt}^{ps}=K_{jt}\cap\{r_1(d'_1)^{6-p}\cdot 2^s<r<r_2(d'_2)^{6-p}\cdot 2^s\},$ where $0<r_1<r_2;$ 
$s\ge 1;\   j=1,\ \dots,\ N;\ p=0,\ \dots,\ 6.$

Introduce the space ${\cal W}^l(K_j^{ps})=\bigoplus\limits_{t=1}^{R_j}W^l(K_{jt}^{ps})$ with the norm
$
 \|v_j\|_{{\cal W}^l(K_j^{ps})}=\left(\sum\limits_{t=1}^{R_j}\|v_{jt}\|_{W^l(\Omega_{jt}^{ps})}^2\right)^{1/2}.
$

\begin{lemma}\label{lAprWTheta'}
 Suppose $s\ge 1,\ \theta\in S^{n-3}.$ Assume that $v_j\in {\cal W}^{2m}(K_j^{0s}),$
$$
  \begin{array}{c}
   {\cal Q}_j(D_y,\ \theta)v_{jt}\in W^l(K_{jt}^{0s}),\\
    \\
   {\cal C}_{j\sigma\mu}(D_y,\ \theta)v = 0 \quad (y\in \gamma_{j\sigma}\cap \bar K_j^{0s}),\
   {\cal T}_{jq\nu}(D_y,\ \theta)v = 0 \quad (y\in \gamma_{jq}\cap \bar K_j^{0s})
  \end{array}
$$
 $$ 
  \begin{array}{c}
 ( j=1,\ \dots,\ N,\ \sigma=1,\ R_j+1,\ \mu=1,\ \dots,\ m,\\
  q=2,\ \dots,\ R_j,\ \nu=1,\ \dots,\ 2m);
 \end{array}
 $$
 then $v\in \prod\limits_j {\cal W}^{l+2m}(K_j^{3s})$ and for all $|\lambda|\ge 1$,
 \begin{equation}\label{eqAprWTheta'}
  \begin{array}{c}
   \sum\limits_j  2^{sa}\|v_j\|_{{\cal W}^{l+2m}(K_j^{6s})}\le
           c \sum\limits_j \bigl\{2^{sa}\sum\limits_t\|{\cal Q}_j(D_y,\ \theta)v_{jt}\|_{W^l(K_{jt}^{3s})}+\\
   +|\lambda|^{-1}2^{sa}\|v_j\|_{{\cal W}^{l+2m}(K_j^{3s})}+
             |\lambda|^{l+2m-1}2^{s(a-l-2m)}\|v_j\|_{L_2(K_j^{3s})} \bigr\},
  \end{array}
 \end{equation}
 where $c>0$ is independent of $v,$ $\theta,\ \lambda$, and $s.$
\end{lemma}
\begin{proof}
Repeating the proof of Lemma~\ref{lAprW'} and replacing $\Omega_j^p$ by $K_j^{ps}$ and $D_z$ by
$\theta$, we get $v\in \prod\limits_j {\cal W}^{l+2m}(K_j^{3s}).$
Now repeating the proof of lemma~3.3 \cite[\S3]{SkDu90} and replacing there
$W^l(\cdot)$ by ${\cal W}^l(\cdot)$ and $H_a^l(\cdot)$ by ${\cal H}_a^l(\cdot),$ from a priori estimate~(\ref{eqAprW'}),
we derive estimate~(\ref{eqAprWTheta'}).
\end{proof}

\begin{theorem}\label{thAprE1'}
 Let $v\in\prod\limits_j {\cal W}_\loc^{2m}(\bar K_j\backslash\{0\})$ be a solution for problem~(\ref{eqQTheta})--(\ref{eqTTheta})
 such that $v\in E_{a-l-2m}^{0,\,N}(K)$ and
  $f\in {\cal E}_a^{l,\,N}(K,\ \gamma)$; then $v\in {\cal E}_a^{l+2m,\,N}(K)$ and
  \begin{equation}\label{eqAprE1'}
   \|v\|_{{\cal E}_a^{l+2m,\,N}(K)}\le c \bigl(\|f\|_{{\cal E}_a^{l,\,N}(K,\ \gamma)}+
       \|v\|_{E_{a-l-2m}^{0,\,N}(K)}\bigr),
  \end{equation}
  where $c>0$ is independent of $v$ and $\theta\in S^{n-3}.$
\end{theorem}

\begin{proof}
The proof is analogous to the proof of Theorem~\ref{thAprE1}, where one must replace
$W^l(\cdot)$, $H_a^l(\cdot)$, $E_a^l(\cdot)$ by ${\cal W}^l(\cdot),$  ${\cal H}_a^l(\cdot),$  ${\cal E}_a^l(\cdot)$;
Lemmas~\ref{lHomogTheta}, \ref{lAprW}, \ref{lAprWTheta} by 
Lemmas~\ref{lHomogTheta'}, \ref{lAprW'}, \ref{lAprWTheta'} correspondingly; Theorem~\ref{thAprH} by 
Theorem~\ref{thAprH'}.
\end{proof}

From Theorem~\ref{thSolvFlat'} and Lemma~\ref{lAprWTheta'}, we obtain the following result
(see theorem~3.1 \cite[\S3]{SkDu90} with $E_a^l(\cdot)$ replaced by ${\cal E}_a^l(\cdot)$).
\begin{theorem}\label{thAprE2'} 
 Suppose the line $\Im\,\lambda=a+1-l-2m$ contains no poles of the operator--valued function $\tilde{\cal M}^{-1}(\lambda)$;
 then for all solutions $v\in {\cal E}_a^{l+2m,\,N}(K)$ to nonlocal transmission 
 problem~(\ref{eqQTheta})--(\ref{eqTTheta}) and all $\theta\in S^{n-3}$, we have
 \begin{equation}\label{eqAprE2'}
   \|v\|_{{\cal E}_a^{l+2m,\,N}(K)}\le c \bigl(\|f\|_{{\cal E}_a^{l,\,N}(K,\ \gamma)}+
      \sum\limits_j \|v_j\|_{L_2(K_j\cap S')}\bigr),
 \end{equation}
 where $S'=\{y\in{\mathbb R}^2:\ 0<R'_1<r<R'_2\};$ $c>0$ is independent of  $\theta$ and $v.$

 If for any $\theta\in S^{n-3}$, estimate~(\ref{eqAprE2'}) holds for all solutions to nonlocal transmission 
 problem~(\ref{eqQTheta})--(\ref{eqTTheta}), then the line $\Im\,\lambda=a+1-l-2m$ contains no poles
 of the operator--valued function $\tilde{\cal M}^{-1}(\lambda).$
\end{theorem}

Theorem~\ref{thAprE2'} implies that kernel of ${\cal M}(\theta)$ is of finite dimension and range of ${\cal L}(\theta)$
is closed.
\section{Adjoint nonlocal problems}\label{sectLMAdj}
In this section, we study operators that are adjoint to the operators of the nonlocal boundary value problems with 
parameter $\theta\in S^{n-3}$.

\subsection{Operators ${\cal L}(\theta)^*.$}
Let ${\cal L}(\theta)=\{{\cal P}_j(D_y,\ \theta),\ {\cal B}_{j\sigma\mu}(D_y,\ \theta)\}:
E_a^{2m,\,N}(K)\to  E_a^{0,\,N}(K,\ \gamma)$ be the operator corresponding to problem~(\ref{eqPTheta}), (\ref{eqBTheta}).
Consider the adjoint operator ${\cal L}(\theta)^*: (E_a^{0,\,N}(K,\ \gamma))^*\to (E_a^{2m,\,N}(K))^*,$ where
$$
 (E_a^{0,\,N}(K,\ \gamma))^*=
\prod\limits_{j=1}^N\bigl\{E_{-a}^0(K_j)\times
\prod\limits_{\sigma=1,\,R_j+1}\prod\limits_{\mu=1}^m(E_a^{2m-m_{j\sigma\mu}-1/2}(\gamma_{j\sigma}))^*
\bigr\},
$$ 
$$
(E_a^{2m,\,N}(K))^*=\prod\limits_{j=1}^N (E_a^{2m}(K_j))^*.
$$  
${\cal L}(\theta)^*$ takes $f=\{f_j,\ g_{j\sigma\mu}\}\in(E_a^{0,\,N}(K,\ \gamma))^*$ to ${\cal L}(\theta)^*f$ by the rule
$$
 \begin{array}{c}
 <u,\ {\cal L}(\theta)^*f>=\sum\limits_j\bigl\{<{\cal P}_j(D_y,\ \theta)u_j,\ f_j>_{K_j}+\\
 +\sum\limits_{\sigma,\mu}
 <{\cal B}_{j\sigma\mu}(D_y,\ \theta)u,\ g_{j\sigma\mu}>_{\gamma_{j\sigma}}\bigr\} 
 \end{array}
$$
for all $u\in E_a^{2m,\,N}(K)$. Here $<\cdot,\ \cdot>,$ 
$<\cdot,\ \cdot>_{K_j},$ $<\cdot,\ \cdot>_{\gamma_{j\sigma}}$ are the sesquilinear forms on the corresponding dual pairs of the spaces. 

Introduce the space ${\cal W}^l(K_j)=\bigoplus\limits_{t=1}^{R_j}W^l(K_{jt})$ with the norm
$
 \|v_j\|_{{\cal W}^l(K_j)}=\left(\sum\limits_{t=1}^{R_j}\|v_{jt}\|_{W^l(K_{jt})}^2\right)^{1/2}.
$
Further (see Theorem~\ref{thLThetaAdjApr}),
we shall see that if the $j$-th component of ${\cal L}(\theta)^*f$ is smooth in $ K_j$ ($j=1,\ \dots,\ N$), then $f_j$ is
smooth only in $ K_{jt}$ and, generally, may have discontinuity on $\gamma_{jq}$ ($q=2,\ \cdots,\ R_j$). This happens
because of nonlocal terms with supports on $\gamma_{jq}$ in the operator ${\cal L}(\theta)$ and therefore in the operator ${\cal L}(\theta)^*$. 
Hence it is natural to consider spaces ${\cal W}^l(\cdot)$ (but not $W^l(\cdot)$) when studying smoothness of $f$.

Consider the functions $\psi_p\in C_0^\infty({\mathbb R}^1)$ such that 
$$
 \begin{array}{c}
 \psi_p(r)=1\ \mbox{for } r_1d_1^{3-p}<r<r_2d_2^{3-p},\\
 \psi_p(r)=0\ \mbox{for }
 r<\frac{\displaystyle 2}{\displaystyle 3}r_1d_1^{3-p}\ \mbox{and }
 r>\frac{\displaystyle 3}{\displaystyle 2}r_2d_2^{3-p},
 \end{array}
$$
where $0<r_1<r_2;$ $p=0,\ \dots,\ 3.$
Put $\hat\gamma_{jq}=\{y:\ \varphi=b_{jq}\ \mbox{or } \varphi=b_{jq}+\pi\}$
 $(j=1,\ \dots,\ N;\ q=1,\ \dots,\ R_j+1).$ Clearly, $\gamma_{jq}\subset\hat\gamma_{jq}$.

\begin{theorem}\label{thLThetaAdjApr}
 Suppose $f=\{f_j,\ g_{j\sigma\mu}\}\in(E_a^{0,\,N}(K,\ \gamma))^*,$
 ${\cal L}(\theta)^*f\in (E_a^{2m,\,N}(K))^*,$ 
 $$
  \psi_0{\cal L}(\theta)^*f\in\left\{
  \begin{array}{l}
    \prod\limits_j W_{\bar K_j}^{-2m+l}({\mathbb R}^n)\footnotemark\ \mbox{for } l<2m,\\
    \prod\limits_j {\cal W}^{-2m+l}(K_j)\ \mbox{for } l\ge2m;
  \end{array}
  \right.
 $$
 then $\psi_3f\in\prod\limits_j\bigl\{{\cal W}^l(K_j)\times
 \prod\limits_{\sigma,\mu}
  W^{-2m+l+m_{j\sigma\mu}+1/2}(\hat\gamma_{j\sigma})
 \bigr\}$ and
 \begin{equation}\label{eqLThetaAdjApr}
  \begin{array}{c}
    \|\psi_3f\|_{\prod\limits_j\{{\cal W}^l( K_j)\times
    \prod\limits_{\sigma,\mu} 
     W^{-2m+l+m_{j\sigma\mu}+1/2}(\hat\gamma_{j\sigma})\}}\le\\
    c_{l}\bigl(\|\psi_0{\cal L}(\theta)^*f\|_{-2m+l}+
    \|\psi_0f\|_{\prod\limits_j\bigl\{W^{-1}_{\bar K_j}({\mathbb R}^n)\times
    \prod\limits_{\sigma,\mu} 
    W^{-2m-1+m_{j\sigma\mu}+1/2}(\hat\gamma_{j\sigma})\bigr\}}\bigr),
  \end{array}
 \end{equation}
 where
 $$
  \|\cdot\|_{-2m+l}=\left\{
  \begin{array}{l}
    \|\cdot\|_{\prod\limits_j W_{\bar K_j}^{-2m+l}({\mathbb R}^n)}\ \mbox{for } l<2m,\\
    \|\cdot\|_{\prod\limits_j {\cal W}^{-2m+l}(K_j)}\ \mbox{for } l\ge2m,
  \end{array}
  \right.
 $$
 \footnotetext{
  $W_{\bar K_j}^{-l}({\mathbb R}^n)$ ($l>0$) is the space that is adjoint to
  $W^l( K_j).$  One can identify the space $W_{\bar K_j}^{-l}({\mathbb R}^n)$ 
  with the subspace of the space $W^{-l}({\mathbb R}^n)$ consisting of distributions with supports from $\bar K_j$
 (see remark~12.4 \cite[Chapter 1, \S 12.6]{LM}).
 }
 $c_{l}>0$ depends on $l\ge0$ and does not depend on $f.$
\end{theorem}

\begin{proof}
1) For any $g\in(E_a^{l-1/2}(\gamma_{jq}))^*$ and $\psi_p$,
denote by $\psi g\otimes \delta(\gamma_{jq})$ the distribution from
$W^{-l}_{\bar K_j}({\mathbb R}^n)$ given by
$$ 
 <u_j,\ \psi g\otimes \delta(\gamma_{jq})>_{ K_j}=<\psi u_j|_{\gamma_{jq}},\ g>_{\gamma_{jq}}\ 
 \mbox{for all } u_j\in W^l( K_j),
$$
$j=1,\ \dots,\ N;\ q=1,\ \dots,\ R_j+1.$ 

Introduce the auxiliary operator 
$$
 \begin{array}{c}
 {\cal L}_{\cal G}(\theta)^*: \prod\limits_{j=1}^N\bigl\{E_{-a}^0( K_j)\times
 \prod\limits_{\sigma=1,\,R_j+1}\prod\limits_{\mu=1}^m
  \bigl(E_a^{2m-m_{j\sigma\mu}-1/2}(\gamma_{j\sigma}))^*\times\\
 \times \prod\limits_{k=1}^N\prod\limits_{q=2}^{R_k}\prod\limits_{s=1}^{S_{j\sigma kq}}
 (E_a^{2m-m_{j\sigma\mu}-1/2}(\gamma_{kq}))^*\bigr)\bigr\}
 \to (E_a^{2m,\,N}( K))^*
 \end{array}
$$
that takes $f'=\{f_j,\ g_{j\sigma\mu},\ g'_{j\sigma\mu kqs}\}\in \prod\limits_j\bigl\{E_{-a}^0( K_j)\times
 \prod\limits_{\sigma,\mu}
         \bigl( E_a^{2m-m_{j\sigma\mu}-1/2}(\gamma_{j\sigma}))^*\times\prod\limits_{k,q,s}
 (E_a^{2m-m_{j\sigma\mu}-1/2}(\gamma_{kq}))^*\bigr)\bigr\}$ to ${\cal L}_{\cal G}(\theta)^*f'$ by the rule
$$
 \begin{array}{c}
 <u,\ {\cal L}_{\cal G}(\theta)^*f'>=\sum\limits_j\bigl\{<{\cal P}_j(D_y,\ \theta)u_j,\ f_j>_{ K_j}+\\
 +\sum\limits_{\sigma,\mu}
 \bigl(<B_{j\sigma\mu}(D_y,\ \theta)u_j|_{\gamma_{j\sigma}},\ g_{j\sigma\mu}>_{\gamma_{j\sigma}}+\\
 \sum\limits_{k,q,s} 
<B_{j\sigma\mu kqs}(D_y,\ \theta)u_k|_{\gamma_{kq}},\ g'_{j\sigma\mu kqs}>_{\gamma_{kq}}\bigr)\bigr\}
\quad \mbox{for all } u\in E_a^{2m,\,N}( K).
\end{array}
$$

Now for every $g_{j\sigma\mu}\in(E_a^{2m-m_{j\sigma\mu}-1/2}(\gamma_{j\sigma}))^*$ and
$\psi_p$ we introduce the distributions
$g^{\cal G}_{j\sigma\mu kqs}\in (E_a^{2m-m_{j\sigma\mu}-1/2}(\gamma_{kq}))^*$ and
$\psi_p g^{\cal G}_{j\sigma\mu kqs}\in W^{-2m+m_{j\sigma\mu}+1/2}(\hat\gamma_{kq})$ given by
$$
 \begin{array}{c}
 <u_{\gamma_{kq}},\ g^{\cal G}_{j\sigma\mu kqs}>_{\gamma_{kq}}=
 <u_{\gamma_{kq}}({\cal G}_{j\sigma kqs}\cdot),\ g_{j\sigma\mu}>_{\gamma_{j\sigma}}\\ 
 \mbox{for all } u_{\gamma_{kq}}\in E_a^{2m-m_{j\sigma\mu}-1/2}(\gamma_{kq})
 \end{array}
$$
and
$$
 \begin{array}{c}
 <W_{\gamma_{kq}},\ \psi_p g^{\cal G}_{j\sigma\mu kqs}>_{\hat\gamma_{kq}}=
 <(\psi_p W_{\gamma_{kq}})({\cal G}_{j\sigma kqs}\cdot),\ g_{j\sigma\mu}>_{\gamma_{j\sigma}}\\
 \mbox{for all } W_{\gamma_{kq}}\in W^{2m-m_{j\sigma\mu}-1/2}(\hat\gamma_{kq}).
 \end{array}
$$
From this, it follows in particular that
$\psi_p g^{\cal G}_{j\sigma\mu kqs}\in W^{-2m+l+m_{j\sigma\mu}+1/2}(\hat\gamma_{kq})$ iff
$\psi_p({\cal G}_{j\sigma kqs}\cdot) g_{j\sigma\mu}\in W^{-2m+l+m_{j\sigma\mu}+1/2}(\hat\gamma_{j\sigma});$ moreover, there are
constants $k_1,\ k_2>0$ (depending on $l$) such that
\begin{equation}\label{eqLThetaAdjApr1}
 \begin{array}{c}
 k_1\|\psi_p({\cal G}_{j\sigma kqs}\cdot) g_{j\sigma\mu}\|_{W^{-2m+l+m_{j\sigma\mu}+1/2}(\hat\gamma_{j\sigma})}\le
 \|\psi_p g^{\cal G}_{j\sigma\mu kqs}\|_{W^{-2m+l+m_{j\sigma\mu}+1/2}(\hat\gamma_{kq})}\le\\
 \le k_2\|\psi_p({\cal G}_{j\sigma kqs}\cdot) g_{j\sigma\mu}\|_{W^{-2m+l+m_{j\sigma\mu}+1/2}(\hat\gamma_{j\sigma})}.
 \end{array}
\end{equation}

Put $f^{\cal G}=\{f_j,\ g_{j\sigma\mu},\ g^{\cal G}_{j\sigma\mu kqs}\}.$ From the definitions of the 
operators ${\cal L}(\theta)^*$ and ${\cal L}_{\cal G}(\theta)^*$, it follows that
\begin{equation}\label{eqLThetaAdjApr2} 
 {\cal L}_{\cal G}(\theta)^*f^{\cal G}={\cal L}(\theta)^*f.
\end{equation}

Denote $\Xi f=\{\Xi_jf_j,\ \Xi_j g_{j\sigma\mu}\},$ 
$\Xi f^{\cal G}=\{\Xi_jf_j,\ \Xi_jg_{j\sigma\mu},\ \Xi_k g^{\cal G}_{j\sigma\mu kqs}\},$ where $\Xi=(\Xi_1,\ \dots,\ \Xi_N),$
$\Xi_j=\Xi_j(\varphi)$ are arbitrary infinitely differentiable on $[b_{j1},\ b_{j,R_j+1}]$ functions. Notice that in the formula
$\Xi f^{\cal G}=\{\Xi_jf_j,\ \Xi_jg_{j\sigma\mu},\ \Xi_k g^{\cal G}_{j\sigma\mu kqs}\},$ a distribution $g^{\cal G}_{j\sigma\mu kqs}$
is multiplied by $\Xi_k$, but not by $\Xi_j$. This will be important further.

2) Let $\zeta_{jq}$ be the functions given by formula~(\ref{eqZeta}). We also consider the functions
\begin{equation}\label{eqHatZeta}
 \hat\zeta_{jq}\in C^\infty({\mathbb R}),\ \hat\zeta_{jq}(\varphi)=1\ \mbox{for } |b_{jq}-\varphi|<3\varepsilon/2,\ 
     \hat\zeta_{jq}(\varphi)=0\ \mbox{for } |b_{jq}-\varphi|>2\varepsilon;
\end{equation}
\begin{equation}\label{eqBarZeta}
 \bar\zeta_{jq}\in C^\infty({\mathbb R}),\ \bar\zeta_{jq}(\varphi)=1\ \mbox{for } |b_{jq}-\varphi|<\varepsilon/8,\ 
     \bar\zeta_{jq}(\varphi)=0\ \mbox{for } |b_{jq}-\varphi|>\varepsilon/4
\end{equation}
$(j=1,\ \dots,\ N;\ q=1,\ \dots,\ R_j+1),$ where $\varepsilon$ is given by formula~(\ref{eqEpsilon}).

Introduce the $N$-dimensional vector--function 
$$
 \Xi^{j'\sigma'}=(0,\ \dots,\ \zeta_{j'\sigma'},\ \dots,\ 0).
$$
Here ``zeroes" are
everywhere, except the $j'$-th position, $j'=1,\ \dots,\ N;$ $\sigma'=1,\ R_{j'}+1.$ If $j\ne j'$, then we have 
$\Xi^{j'\sigma'}_j=0$. If $j=j'$, then we see that the support of $\Xi^{j'\sigma'}_{j'}=\zeta_{j'\sigma'}$ does not intersect 
with $\gamma_{j'q}$, 
but the support of $g^{\cal G}_{j\sigma\mu j'qs}$ is contained in $\gamma_{j'q}$ ($q=2,\ \dots,\ R_{j'}$); therefore,
$\zeta_{j'\sigma'}g^{\cal G}_{j\sigma\mu j'qs}=0$. Thus we have
$$
 \begin{array}{c} 
 {\cal L}_{\cal G}(\theta)^*(\psi_p\Xi^{j'\sigma'}f^{\cal G})=(0,\ \dots,\ {\cal Q}_{j'}(D_y,\ \theta) (\psi_p\zeta_{j'\sigma'}f_{j'})+\\
 +\sum\limits_{\mu=1}^m
 B^*_{j'\sigma'\mu}(D_y,\ \theta)(\psi_p\zeta_{j'\sigma'}g_{j'\sigma'\mu}\otimes\delta(\gamma_{j'\sigma'})),\ \dots,\ 0)
 \end{array}
$$
$(p=0,\ \dots,\ 3).$
Here ``zeroes" are everywhere, except the $j'$-th position, ${\cal Q}_{j'}(D_y,\ \theta)$ and $B^*_{j'\sigma'\mu}(D_y,\ \theta)$ are
formally adjoint to ${\cal P}_{j'}(D_y,\ \theta)$ and $B_{j'\sigma'\mu}(D_y,\ \theta)$ correspondingly.

Notice that the operator 
$$
 {\cal Q}_{j'}(D_y,\ \theta) (\psi_p\zeta_{j'\sigma'}f_{j'})+
 \sum\limits_{\mu=1}^m
 B^*_{j'\sigma'\mu}(D_y,\ \theta)(\psi_p\zeta_{j'\sigma'}g_{j'\sigma'\mu}\otimes\delta(\gamma_{j'\sigma'}))
$$ can be identified with
the adjoint to the operator
$$
 \{{\cal P}_{j'}(D_y,\ \theta)u_{j'},\ B_{j'\sigma'\mu}(D_y,\ \theta)u_{j'}|_{\hat\gamma_{j'\sigma'}}\}_{\mu=1}^m.
$$
Therefore we can use theorem~4.3 \cite[Chapter 2, \S 4.5]{LM}\footnote{
 Theorem~4.3 \cite[Chapter 2, \S 4.5]{LM} deals with operators having variable coefficients; therefore some additional
 restrictions are imposed on supports of considered functions. It is easy to see that these restrictions may be omitted if the 
 coefficients are constant.
}. Thus from relation~(\ref{eqLThetaAdjApr2}) and Leibniz' formula, it follows that
$$
 \begin{array}{c}
\psi_1\Xi^{j'\sigma'}f^{\cal G}\in \prod\limits_j\bigl\{W^l( K_j)\times
 \prod\limits_{\sigma,\mu}
         \bigl( W^{-2m+l+m_{j\sigma\mu}+1/2}(\hat\gamma_{j\sigma})\times\\
  \times \prod\limits_{k,q,s}
  W^{-2m+l+m_{j\sigma\mu}+1/2}(\hat\gamma_{kq})\bigr)\bigr\}
 \end{array}
$$
and
\begin{equation}\label{eqLThetaAdjApr3}
 \begin{array}{c}
 \|\psi_1\Xi^{j'\sigma'}f^{\cal G}\|_{\prod\limits_j\bigl\{W^l( K_j)\times
 \prod\limits_{\sigma,\mu}
         \bigl( W^{-2m+l+m_{j\sigma\mu}+1/2}(\hat\gamma_{j\sigma})\times
 \prod\limits_{k,q,s}
 W^{-2m+l+m_{j\sigma\mu}+1/2}(\hat\gamma_{kq})\bigr)\bigr\}}\le \\
 \le k_3\bigl(\|\psi_0{\cal L}(\theta)^*f\|_{-2m+l}
+\|\psi_0\hat\zeta_{j'\sigma'}f_{j'}\|_{W_{\bar K_{j'}}^{-1}({\mathbb R}^n)}+\\
 +\sum\limits_{\mu=1}^m\|\psi_0 g_{j'\sigma'\mu}\|_
 {W^{-2m-1+m_{j'\sigma'\mu}+1/2}(\hat\gamma_{j'\sigma'})}\bigr).
 \end{array}
\end{equation}
From~(\ref{eqLThetaAdjApr3}) and~(\ref{eqLThetaAdjApr1}), it follows in particular that
$\psi_2g^{\cal G}_{j'\sigma'\mu kqs}\in W^{-2m+l+m_{j'\sigma'\mu}+1/2}(\hat\gamma_{kq})$ and
\begin{equation}\label{eqLThetaAdjApr4}
 \begin{array}{c}
 \|\psi_2g^{\cal G}_{j'\sigma'\mu kqs}\|_{W^{-2m+l+m_{j\sigma\mu}+1/2}(\hat\gamma_{kq})}
 \le k_4\bigl(\|\psi_0{\cal L}(\theta)^*f\|_{-2m+l}+\\
+\|\psi_0\hat\zeta_{j'\sigma'}f_{j'}\|_{W_{\bar K_{j'}}^{-1}({\mathbb R}^n)}+
 \sum\limits_{\mu=1}^m\|\psi_0 g_{j'\sigma'\mu}\|_
 {W^{-2m-1+m_{j'\sigma'\mu}+1/2}(\hat\gamma_{j'\sigma'})}\bigr).
 \end{array}
\end{equation}

3) Put  $\Xi^{k'q'}=(0,\ \dots,\ \zeta_{k'q'},\ \dots,\ 0).$ Here ``zeroes" are everywhere, except the 
$k'$-th position, $k'=1,\ \dots,\ N;$ $q'=2,\ \dots,\ R_{k'}.$ 
If $k\ne k'$, then we have $\Xi^{k'q'}_k=0$. If $k=k'$, then we see that the support of $\Xi^{k'q'}_{k'}=\zeta_{k'q'}$ does not intersect
with the supports of $g_{k'\sigma\mu}$ and $g^{\cal G}_{j\sigma\mu k'qs}$ for $q\ne q'$; therefore, 
$\zeta_{k'q'}g_{k'\sigma\mu}=0$ and $\zeta_{k'q'}g^{\cal G}_{j\sigma\mu k'qs}=0$ for $q\ne q'$. Thus we have
$$ 
 \begin{array}{c}
 {\cal L}_{\cal G}(\theta)^*(\psi_p\Xi^{k'q'}f^{\cal G})=(0,\ \dots,\ {\cal Q}_{k'}(D_y,\ \theta) (\psi_p\zeta_{k'q'}f_{k'})+\\
 +\sum\limits_{j,\sigma,\mu,s}
 B^*_{j\sigma\mu k'q's}(D_y,\ \theta)(\psi_p\zeta_{k'q'}g^{\cal G}_{j\sigma\mu k'q's}\otimes\delta(\gamma_{k'q'})),\ \dots,\ 0)
 \end{array}
$$
$(p=0,\ \dots,\ 3),$
where ``zeroes" are everywhere, except the $k'$-th position, $B^*_{j\sigma\mu k'q's}(D_y,\ \theta)$ is formally adjoint to
$B_{j\sigma\mu k'q's}(D_y,\ \theta).$

 Notice that the operator 
$$
 {\cal Q}_{k'}(D_y,\ \theta) (\psi_p\zeta_{k'q'}f_{k'})+
 \sum\limits_{j,\sigma,\mu,s}
 B^*_{j\sigma\mu k'q's}(D_y,\ \theta)(\psi_p\zeta_{k'q'}g^{\cal G}_{j\sigma\mu k'q's}\otimes\delta(\gamma_{k'q'}))
$$
 can be 
 identified with the adjoint to the operator of the problem
$$ 
 \begin{array}{c}
 {\cal P}_{k'}(D_y,\ \theta)u_{k'}=\hat f_{k'}(y)\ (y\in{\mathbb R}^2),\\
 B_{j\sigma\mu k'q's}(D_y,\ \theta)u_{k'}|_{\hat\gamma_{k'q'}}=\hat g_{j\sigma\mu s}(y)\ (y\in\hat\gamma_{k'q'})\\
 (j=1,\ \dots,\ N;\ \sigma=1,\ R_j+1;\ \mu=1,\ \dots,\ m;\ s=1,\ \dots,\ S_{j\sigma k'q'}).
 \end{array}
$$
This problem differs from the problem studied in Appendix~\ref{appendL*Rn} only in junior terms. 

In 1), we showed 
that $\psi_2g^{\cal G}_{j\sigma\mu k'q's}\in W^{-2m+l+m_{j\sigma\mu}+1/2}(\hat\gamma_{k'q'})$; hence
we can apply theorem~\ref{thLAdjAprRn}. Thus from relation~(\ref{eqLThetaAdjApr2}) and Leibniz' formula,
we obtain
$$
 \begin{array}{c}
\psi_3\Xi^{k'q'}f^{\cal G}\in \prod\limits_j\bigl\{{\cal W}^l( K_j)\times
 \prod\limits_{\sigma,\mu}
         \bigl( W^{-2m+l+m_{j\sigma\mu}+1/2}(\hat\gamma_{j\sigma})\times\\
 \times  \prod\limits_{k,q,s}
  W^{-2m+l+m_{j\sigma\mu}+1/2}(\hat\gamma_{kq})\bigr)\bigr\}
 \end{array}
$$
and
\begin{equation}\label{eqLThetaAdjApr5}
 \begin{array}{c}
 \|\psi_3\Xi^{k'q'}f^{\cal G}\|_{\prod\limits_j\bigl\{{\cal W}^l( K_j)\times
 \prod\limits_{\sigma,\mu}
         \bigl( W^{-2m+l+m_{j\sigma\mu}+1/2}(\hat\gamma_{j\sigma})\times
 \prod\limits_{k,q,s}
  W^{-2m+l+m_{j\sigma\mu}+1/2}(\hat\gamma_{kq})\bigr)\bigr\}}\le \\
 \le k_5\bigl(\|\psi_2{\cal L}(\theta)^*f\|_{-2m+l}+\|\psi_2\hat\zeta_{k'q'}f_{k'}\|_{W_{\bar K_{k'}}^{-1}({\mathbb R}^n)}+\\
 +\sum\limits_{j,\sigma,\mu,s}
 \|\psi_2 g^{\cal G}_{j\sigma\mu k'q's}\|_
 {W^{-2m+l+m_{j\sigma\mu}+1/2}(\hat\gamma_{k'q'})}\bigr).
 \end{array}
\end{equation}
Notice that the space ${\cal W}^l(\cdot)$ appeared just here. As we noted earlier, this is connected with the nonlocal
terms $g^{\cal G}_{j\sigma\mu k'q's}$, which have supports on $\gamma_{k'q'}$ ($q'=2,\ \dots,\ R_{k'}$).

From inequalities~(\ref{eqLThetaAdjApr5}) and~(\ref{eqLThetaAdjApr4}), we get
\begin{equation}\label{eqLThetaAdjApr6}
 \begin{array}{c}
 \|\psi_3\Xi^{k'q'}f^{\cal G}\|_{\prod\limits_j\bigl\{{\cal W}^l( K_j)\times
 \prod\limits_{\sigma,\mu}
         \bigl( W^{-2m+l+m_{j\sigma\mu}+1/2}(\hat\gamma_{j\sigma})\times
  \prod\limits_{k,q,s}
 W^{-2m+l+m_{j\sigma\mu}+1/2}(\hat\gamma_{kq})\bigr)\bigr\}}\le \\
 \le k_6\bigl(\|\psi_0{\cal L}(\theta)^*f\|_{-2m+l}+
 \sum\limits_{j=1}^N \sum\limits_{\sigma=1,\,R_j+1}
 \{\|\psi_0\hat\zeta_{j\sigma}f_j\|_{W_{\bar K_j}^{-1}({\mathbb R}^n)}+\\
 +\sum\limits_{\mu=1}^m\|\psi_0 g_{j\sigma\mu}\|_
 {W^{-2m-1+m_{j\sigma\mu}+1/2}(\hat\gamma_{j\sigma})}\}\bigr).
 \end{array}
\end{equation}

4) Finally, we put $\zeta_{i0}=1-\sum\limits_{q=1}^{R_i+1}\zeta_{iq}$,
$\Xi^{i0}=(0,\ \dots,\ \zeta_{i0},\ \dots,\ 0).$ Here ``zeroes" are everywhere, except the $i$-th position, $i=1,\ \dots,\ N.$

Since the support of $\zeta_{i0}$ does not intersect with $\gamma_{iq}$ 
$(q=1,\ \dots, R_i+1),$ we have
$$ 
 {\cal L}_{\cal G}(\theta)^*(\psi_p\Xi^{i0}f^{\cal G})=(0,\ \dots,\ {\cal Q}_i(D_y,\ \theta)
 (\psi_p\zeta_{i0}f_i),\ \dots,\ 0)
$$
$(p=0,\ \dots,\ 3).$ Here ``zeroes" are everywhere, except the $i$-th position.

 The operator ${\cal Q}_{i}(D_y,\ \theta) (\psi_p\zeta_{j'\sigma'}f_{j'})$ can be identified with
 the adjoint one to the operator of the problem
$$
 {\cal P}_{i}(D_y,\ \theta)u_{i}=\hat f_i(x)\ (y\in{\mathbb R}^2).
$$
Therefore applying theorem~3.1 \cite[Chapter 2, \S 3.2]{LM}, from~(\ref{eqLThetaAdjApr2}) and Leibniz' formula, we get
$$
 \begin{array}{c}
\psi_1\Xi^{i0}f^{\cal G}\in \prod\limits_j\bigl\{W^l( K_j)\times
 \prod\limits_{\sigma,\mu}
         \bigl( W^{-2m+l+m_{j\sigma\mu}+1/2}(\hat\gamma_{j\sigma})\times\\
  \times  \prod\limits_{k,q,s}
 W^{-2m+l+m_{j\sigma\mu}+1/2}(\hat\gamma_{kq})\bigr)\bigr\}
 \end{array}
$$
and
\begin{equation}\label{eqLThetaAdjApr7}
 \begin{array}{c}
 \|\psi_1\Xi^{i0}f^{\cal G}\|_{\prod\limits_j\bigl\{W^l( K_j)\times
 \prod\limits_{\sigma,\mu}
         \bigl( W^{-2m+l+m_{j\sigma\mu}+1/2}(\hat\gamma_{j\sigma})\times
 \prod\limits_{k,q,s}
 W^{-2m+l+m_{j\sigma\mu}+1/2}(\hat\gamma_{kq})\bigr)\bigr\}}\le \\
 \le k_7\bigl(\|\psi_0{\cal L}(\theta)^*f\|_{-2m+l}
+\|\psi_0\bar\zeta_{i0}f_i\|_{W_{\bar K_i}^{-1}({\mathbb R}^n)}\bigr),
 \end{array}
\end{equation}
where $\zeta_{i0}=1-\sum\limits_{q=1}^{R_i+1}\bar\zeta_{iq}$.

Now a priori estimate~(\ref{eqLThetaAdjApr}) follows from
inequalities~(\ref{eqLThetaAdjApr3}), (\ref{eqLThetaAdjApr6}), and (\ref{eqLThetaAdjApr7}).
\end{proof}

\subsection{Connection between kernel of ${\cal L}(\theta)^*$ and kernel of ${\cal M}(\theta)$.}

\begin{lemma}\label{lKerLThetaAdj}
 The kernel  $\ker({\cal L}(\theta)^*)$ of the operator ${\cal L}(\theta)^*$ coincides with the set
$\{v_j,\ F_{j\sigma\mu}(D_y,\ \theta)v|_{\gamma_{j\sigma}}\},$ where $v_j\in {\cal E}_{-a+2m}^{2m}( K_j)$,
$v_{jt}\in C^\infty(\bar K_{jt}\backslash \{0\})$ ($j=1,\dots,\ N$; $t=1,\ \dots,\ R_j$), and $v$ is a solution to 
problem~(\ref{eqQTheta})--(\ref{eqTTheta}) for $\{f_j,\ g_{j\sigma\mu},\ h_{jq\nu}\}=0.$
\end{lemma}
\begin{proof}
1)  In this proof, we shall omit the arguments $(D_y,\ \theta)$ in differential operators; so we shall write ${\cal P}_j$
instead of ${\cal P}_j(D_y,\ \theta)$ and so on.

Suppose $v_j\in {\cal E}_{-a+2m}^{2m}( K_j)$, $v_{jt}\in C^\infty(\bar K_{jt}\backslash \{0\})$ and
$v$ is a solution to problem~(\ref{eqQTheta})--(\ref{eqTTheta}) for $\{f_j,\ g_{j\sigma\mu},\ h_{jq\nu}\}=0.$
Then for any functions $u_j\in C_0^\infty(\bar K_j\backslash \{0\})$, by virtue of Theorem~\ref{thGrP}, we have
\begin{equation}\label{eqKerLThetaAdj1}
 \begin{array}{c}
  \sum\limits_j\bigl\{\sum\limits_t({\cal P}_ju_j,\ v_{jt})_{ K_{jt}}+
   \sum\limits_{\sigma,\mu}
    ({\cal B}_{j\sigma\mu}u,\ F_{j\sigma\mu}v_j|_{\gamma_{j\sigma}})_{\gamma_{j\sigma}}=0.
 \end{array} 
\end{equation}

Since the imbedding operator of 
${\cal E}_{-a+2m}^{2m}( K_j)$ into $E_{-a}^{0}( K_j)$ is bounded, we have 
$v_j\in{\cal E}_{-a}^{0}( K_j).$ Besides, the operator $F_{j\sigma\mu}(D_y,\ \theta)$ is of order $2m-1-m_{j\sigma\mu}$; 
hence, from the Schwarz inequality and Theorem~\ref{thTraceE}, for all 
$u_{\gamma_{j\sigma}}\in E_a^{2m-m_{j\sigma\mu-1/2}}(\gamma_{j\sigma})$, we obtain
$$
 \begin{array}{c}
  |(u_{\gamma_{j\sigma}},\ F_{j\sigma\mu}v_j|_{\gamma_{j\sigma}})_{\gamma_{j\sigma}}|^2\le 
 \int\limits_{\gamma_{j\sigma}}r^{2(a-(2m-m_{j\sigma\mu}-1/2))}|u_{\gamma_{j\sigma}}|^2\,d\gamma\times\\
 \times \int\limits_{\gamma_{j\sigma}}r^{2(-a+2m-(m_{j\sigma\mu}+1/2))}
           |F_{j\sigma\mu}v_j|_{\gamma_{j\sigma}}|^2\,d\gamma\le\\
\le k_1\|u_{\gamma_{j\sigma}}\|_{E_a^{2m-m_{j\sigma\mu-1/2}}(\gamma_{j\sigma})}^2\cdot
 \|F_{j\sigma\mu}v_j|_{\gamma_{j\sigma}}\|_{E_{-a+2m}^{m_{j\sigma\mu}+1/2}(\gamma_{j\sigma})}^2
 \end{array}
$$
Therefore, 
$F_{j\sigma\mu}v_j|_{\gamma_{j\sigma}}\in (E_a^{2m-m_{j\sigma\mu}-1/2}(\gamma_{j\sigma}))^*.$

Thus, $\{v_j,\ F_{j\sigma\mu}(D_y,\ \theta)v|_{\gamma_{j\sigma}}\}\in 
\prod\limits_j\bigl\{E_{-a}^{0}( K_j)\times\prod\limits_{\sigma,\mu}
 (E_a^{2m-m_{j\sigma\mu}-1/2}(\gamma_{j\sigma}))^*\bigr\}$ and from the definition of the operator ${\cal L}(\theta)^*$ 
and identity~(\ref{eqKerLThetaAdj1}), we get
$$
 <u,\ {\cal L}(\theta)^*\{v_j,\ F_{j\sigma\mu}(D_y,\ \theta)v|_{\gamma_{j\sigma}}\}>=0\ \mbox{for all }
 u\in\prod\limits_j C_0^\infty(\bar K_j\backslash \{0\}).
$$
But $\prod\limits_j C_0^\infty(\bar K_j\backslash \{0\})$ is dense in $E_a^{2m,\,N}( K)$; hence,
$\{v_j,\ F_{j\sigma\mu}(D_y,\ \theta)v|_{\gamma_{j\sigma}}\}\in\ker({\cal L}(\theta)^*).$

2) Now suppose $\{v_j,\ \psi_{j\sigma\mu}\}\in\ker({\cal L}(\theta)^*).$ From Theorem~\ref{thLThetaAdjApr}, it follows that
$v_{jt}\in C^\infty(\bar K_{jt}\backslash \{0\})$, $\psi_{j\sigma\mu}\in C^\infty(\gamma_{j\sigma}).$
Then from the definition of the operator ${\cal L}(\theta)^*$, it follows that
$$
  \sum\limits_j\bigl\{({\cal P}_ju_j,\ v_j)_{ K_j}=
     -\sum\limits_{j,\sigma,\mu} ({\cal B}_{j\sigma\mu}u,\ \psi_{j\sigma\mu})_{\gamma_{j\sigma}},\ \mbox{for all }
 u_j\in C_0^\infty(\bar K_j\backslash \{0\}).
$$
The last identity and Green formula~(\ref{eqGrPEta}) imply
\begin{equation}\label{eqKerLThetaAdj2}
\begin{array}{c}
     \sum\limits_j\bigl\{\sum\limits_{\sigma,\mu}
    ({\cal B}_{j\sigma\mu}u,\ F_{j\sigma\mu}v_j|_{\gamma_{j\sigma}}-\psi_{j\sigma\mu})_
    {\gamma_{j\sigma}}
 +\sum\limits_{q,\mu} (B_{jq\mu}u_j|_{\gamma_{jq}},\ {\cal T}_{jq\mu}v)_
      {\gamma_{jq}}\bigr\}=\\
 =\sum\limits_j\bigl\{\sum\limits_t(u_j,\ {\cal Q}_jv_{jt})_{ K_{jt}}+\sum\limits_{\sigma,\mu}
    (B'_{j\sigma\mu}u_j|_{\gamma_{j\sigma}},\ C_{j\sigma\mu}v_j|_{\gamma_{j\sigma}})_{\gamma_{j\sigma}}+\\
 +\sum\limits_{q,\mu} (B'_{jq\mu}u_j|_{\gamma_{jq}},\ {\cal T}_{jq,m+\mu}v)_
        {\gamma_{jq}}\bigr\}.
 \end{array} 
\end{equation}
Putting $\supp u_j\in C_0^\infty(K_{jt})$, from~(\ref{eqKerLThetaAdj2}), we 
obtain ${\cal Q}_jv_{jt}=0,$ $j=1,\ \dots,\ N$; $t=1,\ \dots,\ R_j$.

By Theorem~\ref{thGrP}, the system $\{B_{j\sigma\mu},\ B'_{j\sigma\mu}\}_{\mu=1}^m$ is a Dirichlet system on
$\gamma_{j\sigma}$ ($j=1,\ \dots,\ N;\ \sigma=1,\ R_j+1$) of order $2m$. Therefore, for any system of 
functions $\{\Theta_{j\sigma\nu}\}_{\nu=1}^{2m}\subset C_0^\infty(\gamma_{j\sigma})$ there exist
functions $u_j\in C_0^\infty(\bar K_j\backslash \{0\})$ such that
$$
 \begin{array}{c}
 B_{j\sigma\mu}u_j|_{\gamma_{j\sigma}}=\Theta_{j\sigma\mu},\ 
 B'_{j\sigma\mu}u_j|_{\gamma_{j\sigma}}=\Theta_{j\sigma,\mu+m},\ \mu=1,\ \dots,\ m,\\
 u_j=0\ \mbox{in a neighbourhood of } \gamma_{jq}\ (j=1,\ \dots,\ N;\ q=2,\ \dots,\ R_j)
 \end{array}
$$
(see lemma~2.2 \cite[Chapter 2, \S 2.3]{LM}).
Therefore, taking into account that ${\cal Q}_jv_{jt}=0,$ from~(\ref{eqKerLThetaAdj2}), we obtain
$F_{j\sigma\mu}v_j|_{\gamma_{j\sigma}}-\psi_{j\sigma\mu}=0$ and 
$C_{j\sigma\mu}v_j|_{\gamma_{j\sigma}}=0$.

Similarly, since $\{B_{jq\mu},\ B'_{jq\mu}\}_{\mu=1}^m$ is a Dirichlet system on 
$\gamma_{jq}$ ($j=1,\ \dots,\ N;\ q=2,\ \dots,\ R_j$) of order $2m$, we get
${\cal T}_{jq\nu}v=0$.

Finally, we know that $v_j\in E_{-a}^{0}( K_j)$ by assumption and we showed that 
$v_{jt}\in C^\infty(\bar K_{jt}\backslash \{0\})$; therefore, from Theorem~\ref{thAprH'}, it follows that 
$v_j\in{\cal E}_{-a+2m}^{2m}( K_j).$
\end{proof}

\section{Solvability of nonlocal boundary value problems}\label{sectLSolv}
In this section, we study solvability of nonlocal boundary value problems. In subsection~\ref{subsectLThetaFred},
we establish necessary and sufficient conditions for Fredholm solvability of the nonlocal boundary value problems with
parameter $\theta$ in plane angles. In subsection~\ref{subsectLOneValued}, we study necessary conditions
for Fredholm solvability and sufficient conditions for one--valued solvability of nonlocal boundary value problems in
dihedral angles.

\subsection{Fredholm solvability of nonlocal boundary value problems with parameter $\theta$.}\label{subsectLThetaFred}
\begin{theorem}\label{thLThetaFred}
  Put $a=b+l$. Suppose the line  $\Im\lambda=b+1-2m$ contains no poles of the operator--valued function 
  $\tilde{\cal L}^{-1}(\lambda)$; then the operator
 $$
  {\cal L}(\theta)=\{{\cal P}_j(D_y,\ \theta),\ {\cal B}_{j\sigma\mu}(D_y,\ \theta)\}:
  E_a^{l+2m,\,N}(K)\to  E_a^{l,\,N}(K,\ \gamma)
 $$
 is Fredholm for all  $\theta\in S^{n-3}$.

 If there is a $\theta\in S^{n-3}$ such that the operator ${\cal L}(\theta)$ is Fredholm, then the line
 $\Im\,\lambda=b+1-2m$ contains no poles of the operator--valued function $\tilde{\cal L}^{-1}(\lambda).$
\end{theorem}
\begin{proof}
 Suppose the line  $\Im\lambda=b+1-2m$ contains no poles of $\tilde{\cal L}^{-1}(\lambda)$; then by Theorem~\ref{thAprE2},
 the operator ${\cal L}(\theta)$ has finite dimensional kernel and closed range.

 Let us prove that cokernel of the operator ${\cal L}(\theta)$ is of finite dimension. First, we put $l=0.$
 By Theorems~\ref{lSolvLambda1} and~\ref{lSolvLambda1'}, the operators $\tilde{\cal L}(\lambda)$ and $\tilde{\cal M}(\lambda)$
 are Fredholm and have zero indices. Therefore from Green formula~(\ref{eqGrPLambda}) and Remark~\ref{rGrLambda}, it
 follows that $\lambda_0$ is a pole of $\tilde{\cal L}^{-1}(\lambda)$ iff $\lambda'_0=\bar\lambda_0-2i(m-1)$ is a pole 
 of $\tilde{\cal M}^{-1}(\lambda)$. Hence the line $\Im\,\lambda=(-b+2m)+1-2m$ contains no poles of the operator--valued function
 $\tilde{\cal M}^{-1}(\lambda).$ Now by Theorem~\ref{thAprE2'}, kernel of the operator ${\cal M}(\theta)$ is of finite dimension.
 Finally, Lemma~\ref{lKerLThetaAdj} implies $\dim\ker({\cal L}(\theta)^*)=\dim\ker({\cal M}(\theta))<\infty$.

 Consider the case $l\ge1.$ Suppose $f\in E_{a}^{l,\,N}(K,\ \gamma)$. By the above, there exists a
 $u\in E_{a-l}^{2m,\,N}(K)$ such that ${\cal L}(\theta)u=f$ iff $(f,\ \Psi_i)_{E_{a-l}^{0,\,N}(K,\ \gamma)}=0$ for some linearly
 independent functions 
 $\Psi_i\in E_{a-l}^{0,\,N}(K,\ \gamma)$ ($i=1,\ \dots,\ J$). Here $(\cdot,\ \cdot)_{E_{a-l}^{0,\,N}(K,\ \gamma)}$ is the inner product
 in the Hilbert space $E_{a-l}^{0,\,N}(K,\ \gamma)$. In addition, by Theorem~\ref{thAprE1}, we have $u\in E_{a}^{l+2m,\,N}(K)$.
 
 By virtue of the Schwarz inequality and boundness of the imbeding operator of $E_{a}^{l,\,N}(K,\ \gamma)$ into
 $E_{a-l}^{0,\,N}(K,\ \gamma)$, we have
$$ 
 \begin{array}{c}
 (f,\ \Psi_i)_{E_{a-l}^{0,\,N}(K,\ \gamma)}\le\|f\|_{E_{a-l}^{0,\,N}(K,\ \gamma)}\|\Psi_i\|_{E_{a-l}^{0,\,N}(K,\ \gamma)}\le \\
 \\ k_1 \|f\|_{E_{a}^{l,\,N}(K,\ \gamma)}\|\Psi_i\|_{E_{a-l}^{0,\,N}(K,\ \gamma)}
 \end{array}
$$
for all $f\in E_{a}^{l,\,N}(K,\ \gamma)$. Therefore, by virtue of the Riesz theorem concerning a general form of a linear
functional in a Hilbert space, there exist linearly independent functions 
$\hat\Psi_i\in E_{a}^{l,\,N}(K,\ \gamma)$ ($i=1,\ \dots,\ J$) such that
$$
  (f,\ \Psi_i)_{E_{a-l}^{0,\,N}(K,\ \gamma)}=(f,\ \hat\Psi_i)_{E_{a}^{l,\,N}(K,\ \gamma)}\ \mbox{for all } f\in E_{a}^{l,\,N}(K,\ \gamma).
$$
This means that cokernel of the operator ${\cal L}(\theta)$ is of the same finite dimension $J$ for all $l\ge 0$.

The second part of the Theorem follows from Theorem~\ref{thAprE2}.
\end{proof}

\subsection{Solvability of nonlocal boundary value problems in dihedral angles.}\label{subsectLOneValued}
\begin{theorem}\label{thLSolv}
 Put $a=b+l$. Suppose the line  $\Im\,\lambda=b+1-2m$ contains no poles of the operator--valued function
 $\tilde{\cal L}^{-1}(\lambda)$. Suppose also that for $l=0$, we have $\dim\ker({\cal L}(\theta))=0$ for all
 $\theta\in S^{n-3},$ $\codim{\cal R}({\cal L}(\theta_0))=0$ for some $\theta_0\in S^{n-3}$; then the operator
 $$
  {\cal L}=\{{\cal P}_j(D_y,\ D_z),\ {\cal B}_{j\sigma\mu}(D_y,\ D_z)\}:
  H_a^{l+2m,\,N}(\Omega)\to H_a^{l,\,N}(\Omega,\ \Gamma)
 $$
 is an isomorphism.
\end{theorem}
\begin{proof}
By Theorem~\ref{thAprE2}, we have $\dim\ker({\cal L}(\theta))<\infty$ and range ${\cal R}({\cal L}(\theta))$ is closed in
$E_a^{l,\,N}(K,\ \gamma)$ for all $\theta\in S^{n-3}.$ 

Since the operator ${\cal L}(\theta)$ is bounded and $\dim\ker({\cal L}(\theta))=0$ for $l=0$, we have
\begin{equation}\label{eqLSolv1}
 k_1 \|{\cal L}(\theta)u\|_{E_a^{0,\,N}(K,\ \gamma)}\le
   \|u\|_{E_a^{2m,\,N}(K)}\le k_2\|{\cal L}(\theta)u\|_{E_a^{0,\,N}(K,\ \gamma)},
\end{equation}
where $k_1,\ k_2>0$ are independent of $\theta\in S^{n-3}$ and $u$ ($k_2$ does not depend on $\theta\in S^{n-3}$, since
the sphere $S^{n-3}$ is compact).

By assumption, there exists a $\theta_0\in S^{n-3}$ such that the operator ${\cal L}(\theta_0)$ has a bounded inverse. Therefore,
using estimates~(\ref{eqLSolv1}) and the method of continuation with respect to the parameter 
$\theta\in S^{n-3}$ (see the proof of theorem~7.1 \cite[Chapter 2, \S7]{Lad}), we prove that the operator
${\cal L}(\theta)$ has a bounded inverse for all $\theta\in S^{n-3}.$

Reduce problem~(\ref{eqP}), (\ref{eqB}) to problem~(\ref{eqPTheta}), (\ref{eqBTheta}) doing the Fourier transform
with respect to $z: U(y,\ z)\to \hat U(y,\ \eta)$ and changing variables: $y'=|\eta|\cdot y$. Now repeating the proof of 
lemma~7.3 \cite[\S7]{MP} and applying Theorem~\ref{thAprH} of this work, we complete the proof.
\end{proof}

\begin{theorem}\label{thLNecessCond} 
 Suppose for some $b\in{\mathbb R},$ $l_1\ge0,$ the operator
  $$
  {\cal L}=\{{\cal P}_j(D_y,\ D_z),\ {\cal B}_{j\sigma\mu}(D_y,\ D_z)\}:
  H_{a_1}^{l_1+2m,\,N}(\Omega)\to H_{a_1}^{l_1,\,N}(\Omega,\ \Gamma),\ a_1=b+l_1,
 $$
 is Fredholm; then the operator
 $$
  {\cal L}(\theta)=\{{\cal P}_j(D_y,\ \theta),\ {\cal B}_{j\sigma\mu}(D_y,\ \theta)\}:
  E_a^{l+2m,\,N}(K)\to  E_a^{l,\,N}(K,\ \gamma),\ a=b+l,
 $$
 is an isomorphism for all $\theta\in S^{n-3},$ $l=0,\ 1,\ \dots$
\end{theorem}
\begin{proof} 
1) While proving the Theorem, we shall follow the scheme of the paper \cite[\S8]{MP}. 

Similarly to the proof of lemma~8.1 \cite[\S8]{MP}, one can prove that the operator ${\cal L}$ is an isomorphism for $l=l_1$,
$a=a_1$.
 Therefore we have
$$
 \|U\|_{H_{a_1}^{l_1+2m,\,N}(\Omega)}\le k_1\|{\cal L}U\|_{H_{a_1}^{l_1,\,N}(\Omega,\ \Gamma)}.
$$
Substituting $U^p(y,\ z)=p^{1-n/2}e^{i(\theta,\ z)}\varphi(z/p)u(y)$ ($\varphi\in C_0^\infty({\mathbb R}^{n-2})$,
$u\in E_{a_1}^{l_1+2m,\,N}(K)$, $\theta\in S^{n-3}$)
into the last inequality and passing to the limit as $p\to\infty$,
we get
\begin{equation}\label{eqLNecessCond1}
  \|u\|_{E_{a}^{l+2m,\,N}(K)}\le k_2\|{\cal L}(\theta)u\|_{E_{a}^{l,\,N}(K,\ \gamma)}
\end{equation}
for $l=l_1$, $a=a_1$.
This implies that ${\cal L}(\theta)$ has trivial kernel for $l=l_1$, $a=a_1$. 
But by Theorem~\ref{thAprE1}, kernel of ${\cal L}(\theta)$ does not depend on $l$ and $a=b+l$; therefore the operator
${\cal L}(\theta)$ has trivial kernel for all $l$ and $a=b+l$. 

By Theorem~\ref{thAprE2}, estimate~(\ref{eqLNecessCond1}) implies that the line
$\Im\,\lambda=b+1-2m$ contains no poles of the operator--valued function $\tilde{\cal L}^{-1}(\lambda).$ Hence, by
Theorem~\ref{thLThetaFred}, the operator ${\cal L}(\theta)$ is Fredholm for all $l$ and $a=b+l$. From this and from
triviality of $\ker{\cal L}(\theta)$, it follows that estimate~(\ref{eqLNecessCond1}) is valid for all $l$ and $a=b+l$.

2) Repeating the proof of lemma~7.3 \cite[\S7]{MP}, from estimate~(\ref{eqLNecessCond1}), we get
$$
 \|U\|_{H_{a}^{2m,\,N}(\Omega)}\le k_3\|{\cal L}U\|_{H_{a}^{0,\,N}(\Omega,\ \Gamma)},
$$
where $l=0$, $a=b$. Therefore, the operator ${\cal L}:H_{b}^{2m,\,N}(\Omega)\to H_{b}^{0,\,N}(\Omega,\ \Gamma)$
has trivial kernel and closed range. Let us show that its range coincides with $H_{b}^{0,\,N}(\Omega,\ \Gamma)$.
Indeed, since $H_{b+l_1}^{l_1+2m,\,N}(\Omega)\subset H_{b}^{2m,\,N}(\Omega)$, range ${\cal R}({\cal L})_{b+l_1}$ of the operator
${\cal L}:H_{b+l_1}^{l_1+2m,\,N}(\Omega)\to H_{b+l_1}^{l_1,\,N}(\Omega,\ \Gamma)$ is contained in range ${\cal R}({\cal L})_b$
of the operator ${\cal L}:H_{b}^{2m,\,N}(\Omega)\to H_{b}^{0,\,N}(\Omega,\ \Gamma)$:
$$
 {\cal R}({\cal L})_{b+l_1}\subset{\cal R}({\cal L})_b.
$$
By proved in 1),  
${\cal R}({\cal L})_{b+l_1}=H_{b+l_1}^{l_1,\,N}(\Omega,\ \Gamma)$ which is dense in $H_{b}^{0,\,N}(\Omega,\ \Gamma)$; hence,
${\cal R}({\cal L})_{b}$ is also dense in $H_{b}^{0,\,N}(\Omega,\ \Gamma)$. But ${\cal R}({\cal L})_{b}$ is closed; therefore,
${\cal R}({\cal L})_{b}=H_{b}^{0,\,N}(\Omega,\ \Gamma)$. 

So, we have proved that the operator ${\cal L}:H_{b}^{2m,\,N}(\Omega)\to H_{b}^{0,\,N}(\Omega,\ \Gamma)$ is an isomorphism.

3) Now we shall prove the estimate
\begin{equation}\label{eqLNecessCond2}
 \|V\|_{{\cal H}_{-b+2m}^{2m,\,N}(\Omega)}\le k_4\|{\cal M}V\|_{{\cal H}_{-b+2m}^{0,\,N}(\Omega,\ \Gamma)}.
\end{equation}

Denote by ${\rm P}: H_{b-2m}^{0,\,N}(\Omega)\to  H_b^{0,\,N}(\Omega)$ the unbounded operator corresponding to 
problem~(\ref{eqP}), (\ref{eqB}) with homogeneous nonlocal conditions. The operator ${\rm P}$ is given by
$$
 \begin{array}{c}
 \Dom({\rm P})=\{U\in H_b^{2m,\,N}(\Omega):\ {\cal B}_{j\sigma\mu}(D_y,\ D_z)U=0,\\
 j=1,\ \dots,\ N;\ \sigma=1,\ R_j+1;\ \mu=1,\ \dots,\ m\},
 \end{array}
$$
$$
 {\rm P}U=({\cal P}_1(D_y,\ D_z)U_1,\ \dots,\ {\cal P}_N(D_y,\ D_z)U_N),\quad U\in \Dom({\rm P}).
$$

Denote by ${\rm Q}: H_{-b}^{0,\,N}(\Omega)\to H_{-b+2m}^{0,\,N}(\Omega)$ the unbounded operator
corresponding to problem~(\ref{eqQ})--(\ref{eqT}) with homogeneous boundary conditions and homogeneous
nonlocal transmission conditions. The operator ${\rm Q}$ is given by
$$
 \begin{array}{c}
 \Dom({\rm Q})=\{V\in {\cal H}_{-b+2m}^{2m,\,N}(\Omega):\ {\cal C}_{j\sigma\mu}(D_y,\ D_z)V=0,\ 
 {\cal T}_{jq\nu}(D_y,\ D_z)V=0,\\
  j=1,\ \dots,\ N;\ \sigma=1,\ R_j+1;\ \mu=1,\ \dots,\ m;\\
  q=2,\ \dots,\ R_j;\ \nu=1,\ \dots,\ 2m\}
 \end{array}
$$
$$
 {\rm Q}V=(W_1,\ \dots,\ W_N),\ W_j={\cal Q}_j(D_y,\ D_z)V_{jt}\ \mbox{for } x\in\Omega_{jt},\ V\in \Dom({\rm Q}).
$$
It is clear that $\Dom({\rm P})$ is dense in $H_{b-2m}^{0,\,N}(\Omega)$ and $\Dom({\rm Q})$ is dense in
${\cal H}_{-b}^{0,\,N}(\Omega).$ From Theorems~\ref{thAprH} and~\ref{thAprH'}, it follows that the operators 
${\rm P}$ and ${\rm Q}$ are closed.
Since the operator ${\cal L}:H_{b}^{2m,\,N}(\Omega)\to H_{b}^{0,\,N}(\Omega,\ \Gamma)$ is an isomorphism, the
operator ${\rm P}$ is also an isomorphism from $\Dom({\rm P})$ onto $H_b^{0,\,N}(\Omega)$.

Denote by ${\rm P}^*: H_{-b}^{0,\,N}(\Omega)\to H_{-b+2m}^{0,\,N}(\Omega)$ the operator that is adjoint 
to ${\rm P}$ with respect to the inner product
$\sum\limits_j(U_j,\ V_j)_{\Omega_j}$ in $\prod\limits_j L_2(\Omega_j).$ Since the operator ${\rm P}$ is an isomorphism
from $\Dom({\rm P})$ onto $H_b^{0,\,N}(\Omega)$, the 
operator ${\rm P}^*$ is also an isomorphism from $\Dom({\rm P}^*)$ onto 
$H_{-b+2m}^{0,\,N}(\Omega)$ and its domain $\Dom({\rm P}^*)$ is dense in
$H_{-b}^{0,\,N}(\Omega)$. The operator ${\rm P}^*$ is given by
$$
 \sum\limits_j\bigl({\rm P}_jU_j,\ V_j\bigr)_{\Omega_j}=
 \sum\limits_j\bigr(U_j,\ ({\rm P}^*V)_j\bigr)_{\Omega_j}\ \mbox{for any } U\in\Dom({\rm P}), V\in\Dom({\rm P}^*).
$$

Since the closed operator ${\rm P}^*$ is an isomorphism from $\Dom({\rm P}^*)$ onto 
$H_{-b+2m}^{0,\,N}(\Omega)$, we have
\begin{equation}\label{eqAdjP1}
 \|V\|_{H_{-b}^{0,\,N}(\Omega)}\le k_5\|{\rm P}^*V\|_{H_{-b+2m}^{0,\,N}(\Omega)}
\end{equation}
for all $V\in \Dom({\rm P}^*)$, where $k_5>0$ is independent of $V.$

From Theorem~\ref{thGrP} and Remark~\ref{rGr}, it follows that ${\rm Q}\subset{\rm P}^*$.\footnote{One 
  can prove that ${\rm Q}={\rm P}^*$, but for our
  purposes, it is sufficient to prove the weaker result.} Therefore using~(\ref{eqAdjP1}),
we get 
$$
 \|V\|_{{\cal H}_{-b}^{0,\,N}(\Omega)}\le k_5\|{\rm Q}V\|_{H_{-b+2m}^{0,\,N}(\Omega)}
$$
for all $V\in \Dom({\rm Q})$.
From the last inequality, Lemma~\ref{lHomog'}, and Theorem~\ref{thAprH'}, we obtain estimate~(\ref{eqLNecessCond2}).

4) Substituting $V^p(y,\ z)=p^{1-n/2}e^{i(\theta,\ z)}\varphi(z/p)v(y)$ ($\varphi\in C_0^\infty({\mathbb R}^{n-2})$,
$v\in {\cal E}_{-b+2m}^{2m,\,N}(K)$, $\theta\in S^{n-3}$)
into inequality~(\ref{eqLNecessCond2}) and passing to the limit as $p\to\infty$,
we get
$$
 \|v\|_{{\cal E}_{-b+2m}^{2m,\,N}(K)}\le k_6\|{\cal M}(\theta)v\|_{{\cal E}_{-b+2m}^{0,\,N}(K,\ \gamma)}.
$$
Therefore kernel of the operator ${\cal M}(\theta):{\cal E}_{-b+2m}^{2m,\,N}(K)\to {\cal E}_{-b+2m}^{0,\,N}(K,\ \gamma)$ is trivial.
By virtue of Lemma~\ref{lKerLThetaAdj}, $\dim\ker({\cal L}(\theta)^*)=\dim\ker({\cal M}(\theta))=0$. Combining this with 1),
we see that the operator ${\cal L}(\theta):E_b^{2m,\,N}(K)\to E_b^{0,\,N}(K,\ \gamma)$ is an isomorphism. 
Using Theorem~\ref{thAprE1'}, we prove the Theorem for arbitrary $l$ and $a=b+l$.
\end{proof}

\begin{remark}
  From Theorems~\ref{thLThetaFred} and~\ref{thLNecessCond}, it follows that the operator 
 $
  {\cal L}:H_a^{l+2m,\,N}(\Omega)\to H_a^{l,\,N}(\Omega,\ \Gamma)
 $
is an isomorphism for all $l$ and $a=b+l$ whenever 
 $
  {\cal L}:H_{a_1}^{l_1+2m,\,N}(\Omega)\to H_{a_1}^{l_1,\,N}(\Omega,\ \Gamma)
 $
is Fredholm for some $l_1$ and $a_1=b+l_1$.
\end{remark}
\section{One--valued solvability of nonlocal problems for the Poisson equation in dihedral angles}\label{appendPoisson}
As an application of the results obtained in this work we shall prove the one--valued solvability of
nonlocal problems for the Poisson equation in dihedral angles. For this purpose we need to study corresponding
auxiliary nonlocal problems in plane angles which is done by reducing them to boundary value problems for
differential--difference equations (see \cite{SkMs86, SkDu89, SkBook}).
\subsection{Difference operators in plane angles.}
Put
$$
   K=\{y\in{\mathbb R}^2:\ r>0,\ b_1<\varphi<b_{R+1}\},
$$
$$
  K_{t}=\{y\in{\mathbb R}^2:\ r>0,\ b_t<\varphi<b_{t+1}\}\ (t=1,\ \dots,\ R),
$$
$$
  \gamma_{q}=\{y\in{\mathbb R}^2:\ r>0,\ \varphi=b_{q}\}\ (q=1,\ \dots,\ R+1),
$$
where $R\ge1$ is an integer;
$0<b_1<b_2<\dots<b_R<b_{R+1}<2\pi;$ $b_2-b_1=\dots=b_{R+1}-b_R=d>0.$

Consider the difference operator ${\cal R}: L_2({\mathbb R}^2)\to L_2({\mathbb R}^2)$ given by
$$
  ({\cal R}w)(y)=\sum\limits_{p=-R+1}^{R-1}e_p\cdot w(r,\ \varphi+pd),
$$
where $w(r,\ \varphi)$ is the function $w(y)$ written in the polar coordinates; $e_p\in{\mathbb R}$.

Let $I_K: L_2(K)\to L_2({\mathbb R}^2)$ be the operator of extension by zero outside $K;$
$P_K: L_2({\mathbb R}^2)\to L_2( K)$ be the operator of restriction to $K$. Introduce the operator
${\cal R}_K: L_2(K)\to L_2( K)$ given by
$$
 {\cal R}_K=P_K{\cal R}I_K.
$$

The following statement is obvious.
\begin{lemma}\label{lRBoundL2}
 The operators ${\cal R}: L_2({\mathbb R}^2)\to L_2({\mathbb R}^2),$ ${\cal R}_K: L_2( K)\to L_2( K)$ are bounded.

$({\cal R}^*w)(x)=\sum\limits_{p=-R+1}^{R-1} e_p\cdot w(r,\ \varphi-pd);$
 ${\cal R}^*_K=P_K{\cal R}^*I_K.$
\end{lemma}

Introduce an isomorphism of the Hilbert spaces ${\cal U}: L_2(K)\to L_2^R(K_1)$ by the formula
$$
 ({\cal U}w)_t(y)=w(r,\ \varphi+b_t-b_1)\quad (y\in K_1;\ t=1,\ \dots,\ R),
$$
where $L_2^R( K_1)=\prod\limits_{t=1}^R L_2(K_1).$

Denote by ${\cal R}_1$ the matrix of order $R\times R$ with the elements
$$
 r_{p_1p_2}=e_{p_2-p_1}\ (p_1,\ p_2=1,\ \dots,\ R).
$$

\begin{lemma}\label{lRMatrix}
 The operator ${\cal U}{\cal R}_K{\cal U}^{-1}: L_2^R( K_1)\to  L_2^R( K_1)$ is the operator of multiplication by the matrix
  ${\cal R}_1.$
\end{lemma}

\begin{lemma}\label{lRSpectrum}
 Spectrum of the operator ${\cal R}_K: L_2( K)\to L_2( K)$ coincides with spectrum of the matrix
 ${\cal R}_1.$
\end{lemma}

\begin{lemma}\label{lRPosDef}
 The operator ${\cal R}_K+{\cal R}_K^*: L_2( K)\to L_2( K)$ is positive definite if and only if
 the matrix ${\cal R}_1+{\cal R}_1^*$ is positive definite.
\end{lemma}
Lemmas~\ref{lRMatrix}--\ref{lRPosDef} are analogous to lemmas~8.6--8.8 \cite[Chapter 2, \S8]{SkBook}.

Introduce the spaces $W^l(K)$ and $\mathaccent23 W^l( K)$ as a completion of the sets $C_0^\infty(\bar K\backslash\{0\})$ and $C_0^\infty(K)$ correspondingly in the norm
$
 \left(
    \sum\limits_{|\alpha|\le l}\int\limits_K |D_y^\alpha w(y)|^2 dy
                                       \right)^{1/2}.
$
Similarly, we introduce the space $W^l(K_t).$ 

Denote by $w_t$ the restriction of a function $w$ to $K_t$.
Consider the spaces ${\cal W}^l( K)=\bigoplus\limits_{t=1}^R W^l(K_t)$ and ${\cal E}_a^l( K)=\bigoplus\limits_{t=1}^R E_a^l(K_t)$
with the norms
$
 \|w\|_{{\cal W}^l( K)}=\left(\sum\limits_{t=1}^{R}\|w_{t}\|_{W^l( K_{t})}^2\right)^{1/2}
$
and
$
 \|w\|_{{\cal E}_a^l( K)}=\left(\sum\limits_{t=1}^{R}\|w_{t}\|_{E_a^l( K_{t})}^2\right)^{1/2}
$
correspondingly.

\begin{lemma}\label{lRW0toW}
 The operator ${\cal R}_K$ maps continuously $\mathaccent23 W^l( K)$ into $W^l( K)$ and for all
 $w\in\mathaccent23 W^l( K),$
 $$
  D^\alpha{\cal R}_K w={\cal R}_K D^\alpha w\quad (|\alpha|\le l).
 $$
\end{lemma}
Lemma~\ref{lRW0toW} is analogous to lemma~8.13 \cite[Chapter 1, \S8]{SkBook}.

\begin{lemma}\label{lRWtoW}
 The operator ${\cal R}_K$ maps continuously ${\cal W}^l(K)$ into ${\cal W}^l( K)$ and ${\cal E}_a^l(K)$ into ${\cal E}_a^l( K)$.

 If $\det {\cal R}_1\ne 0,$ then the operator ${\cal R}_K^{-1}$ also maps continuously
 ${\cal W}^l( K)$ into ${\cal W}^l( K)$ and ${\cal E}_a^l(K)$ into ${\cal E}_a^l( K)$.
\end{lemma}
The proof follows from Lemmas~\ref{lRMatrix}, \ref{lRSpectrum}.

\subsection{Differential--difference operators in plane angles.}
Consider the differential--difference equation
\begin{equation}\label{eqPR}
 {\cal P}_{\cal R} w=-\sum\limits_{i,\,j=1}^2 ({\cal R}_{ij K}w_{y_j})_{y_i}+
 \sum\limits_{i=1}^2 {\cal R}_{i K}w_{y_i}+{\cal R}_{0 K}w=f(y) \quad (y\in  K)
\end{equation}
with the boundary conditions
\begin{equation}\label{eqBoundHomog}
 w|_{\gamma_1}=w|_{\gamma_{R+1}}=0,
\end{equation}
where ${\cal R}_{ij K}=P_K{\cal R}_{ij}I_K,$
${\cal R}_{i K}=P_K{\cal R}_{i}I_K,$
 ${\cal R}_{0 K}=P_K{\cal R}_{0}I_K;$
$$
  {\cal R}_{ij}w(y)=\sum\limits_{p=-R+1}^{R-1}e_{ijp}\cdot w(r,\ \varphi+pd) \quad (i,\ j=1,\ 2);
$$
$$
  {\cal R}_{i}w(y)=\sum\limits_{p=-R+1}^{R-1}e_{ip}\cdot w(r,\ \varphi+pd) \quad (i=0,\ 1,\ 2);
$$
$e_{ijp},\ e_{ip}\in{\mathbb R};$ $f\in L_2( K)$.

Denote by $(\cdot,\ \cdot)_K$ the inner product in $L_2( K).$

\begin{definition}\label{defSE}
 We shall say that differential--difference equation~(\ref{eqPR}) is strongly elliptic in $\bar K$ if
 for all $w\in C_0^\infty(\bar K\backslash\{0\})$,
 \begin{equation}\label{eqSE}
  \Re({\cal P}_{\cal R}w,\ w)_K\ge c_1\|w\|^2_{W^1( K)}-c_2\|w\|^2_{L_2( K)},
 \end{equation}
 where $c_1>0,\ c_2\ge0$ do not depend on $w.$
\end{definition}

\begin{definition}\label{defSol1}
 A function $w\in\mathaccent23 W^1( K)$ is called a generalized solution for problem~(\ref{eqPR}), (\ref{eqBoundHomog}) if for all
 $u\in\mathaccent23 W^1( K)$,
 $$
  \sum\limits_{i,\,j=1}^2 ({\cal R}_{ij K}w_{y_j},\ u_{y_i})_K+
 \sum\limits_{i=1}^2 ({\cal R}_{i K}w_{y_i},\ u)_K+({\cal R}_{0 K}w,\ u)_K=
 (f,\ u)_K.
 $$
\end{definition}

We define the unbounded operator ${\mathbb P}_{\cal R}: L_2( K)\to L_2( K)$ with domain
$\Dom({\mathbb P}_{\cal R})=\{w\in\mathaccent23 W^1( K):\ {\cal P}_{\cal R}w\in L_2( K)\}$ acting in the space of distributions
$D'( K)$ by the formula
$$
  {\mathbb P}_{\cal R} w=-\sum\limits_{i,\,j=1}^2 ({\cal R}_{ij K}w_{y_j})_{y_i}+
 \sum\limits_{i=1}^2 {\cal R}_{i K}w_{y_i}+{\cal R}_{0 K}w
$$
The operator ${\mathbb P}_{\cal R}$ is called a {\it differential--difference operator.}

It is easy to show that Definition~\ref{defSol1} is equivalent to the following one.
\begin{definition}\label{defSol2}
 A function $w\in D({\mathbb P}_{\cal R})$ is called a generalized solution for problem~(\ref{eqPR}), (\ref{eqBoundHomog})
 if
 $$
  {\mathbb P}_{\cal R}w=f.
 $$
\end{definition}

Denote by $\sigma({\mathbb P}_{\cal R})$ spectrum of the operator
${\mathbb P}_{\cal R}: L_2( K)\to L_2( K).$

Using the strong ellipticity of the operator ${\mathbb P}_{\cal R}$ and Lemmas~\ref{lRBoundL2},
\ref{lRMatrix}, \ref{lRW0toW}, one can prove the following result
(cf. theorem~10.1 \cite[Chapter 2, \S10]{SkBook}).
\begin{theorem}\label{thPRSpectrum}
Suppose differential--difference equation~(\ref{eqPR}) is strongly elliptic; then
 $$
  \sigma({\mathbb P}_{\cal R})\subset\{\lambda\in{\mathbb C}: \Re\lambda>-c_2\},
 $$
 where $c_2\ge0$ is a constant in~(\ref{eqSE}).
\end{theorem}

\begin{example}\label{exDD}
{\rm
Consider the equation
\begin{equation}\label{eqPRExDD}
 -\triangle {\cal R}_K w(y)+{\cal R}_K w(y)=f(y) \quad
  (y\in K=\{y\in{\mathbb R}^2:\ r>0,\ b_1<\varphi<b_3\})
\end{equation}
with the boundary conditions
\begin{equation}\label{eqBoundHomogExDD}
 w|_{\gamma_1}=w|_{\gamma_3}=0,
\end{equation}
where ${\cal R}w(y)=w(r,\ \varphi)-\alpha w(r,\ \varphi+d)-\beta w(r,\ \varphi-d),$ $d=b_3-b_2=b_2-b_1;$
$\alpha,\ \beta\in {\mathbb R}$; $|\alpha+\beta|<2.$

Clearly, the matrix ${\cal R}_1$ has the form
$$
 {\cal R}_1=\left(
 \begin{array}{cc}
  1 & -\alpha\\
  -\beta & 1
 \end{array}
 \right).
$$
Using Lemma~\ref{lRW0toW}, for all $w\in C_0^\infty( K\backslash\{0\})$, we get
$$
 \begin{array}{c}
  \Re(-\triangle {\cal R}_K w+{\cal R}_K w,\ w)_K=\\
  =\frac{\displaystyle 1}{\displaystyle 2}\sum\limits_{i=1}^2 (({\cal R}_K+{\cal R}_K^*) w_{y_i},\
    w_{y_i})_K+\frac{\displaystyle 1}{\displaystyle 2}(({\cal R}_K+{\cal R}_K^*) w,\ w)_K.
 \end{array}
$$
Since $|\alpha+\beta|<2$, the matrix ${\cal R}_1+{\cal R}_1^*$ is positive definite; therefore, by Lemma~\ref{lRPosDef},
the operator ${\cal R}_K+{\cal R}_K^*$ is also positive definite. From this and from the last equality, we obtain
$$
  \Re(-\triangle {\cal R}_K w+{\cal R}_K w,\ w)_K\ge c_1\|w\|^2_{W^1( K)}.
$$
Hence by Theorem~\ref{thPRSpectrum}, boundary value problem~(\ref{eqPRExDD}),
(\ref{eqBoundHomogExDD}) has a unique generalized solution
$w\in\mathaccent23 W^1( K)$ for every $f\in L_2( K).$
}
\end{example}

\subsection{Nonlocal problems for the Poisson equation in dihedral angles.}
Put
$$
 \Omega=\{x=(y,\ z):\ r>0,\ b_1<\varphi<b_3,\ z\in {\mathbb R}^{n-2}\},
$$
$$
  \Omega_{t}=\{x=(y,\ z):\ r>0,\ b_t<\varphi<b_{t+1},\ z\in {\mathbb R}^{n-2}\}\ (t=1,\ 2),
$$
$$
  \Gamma_{q}=\{x=(y,\ z):\ r>0,\ \varphi=b_{q},\ z\in {\mathbb R}^{n-2}\}\ (q=1,\ \dots,\ 3),
$$
where $b_2-b_1=b_{3}-b_2=d>0.$

Consider the nonlocal boundary value problem
\begin{equation}\label{eqPPuasson}
 -\triangle U\equiv -\sum\limits_{i=1}^n U_{x_ix_i}(x)=f(x) \quad
      (x\in\Omega),
\end{equation}
\begin{equation}\label{eqBPuasson}
 \begin{array}{c}
 U|_{\Gamma_{1}}+\alpha U(r,\ \varphi+d,\ z)|_{\Gamma_{1}}=
  g_{1}(x)\ \quad (x\in\Gamma_{1}),\\
 U|_{\Gamma_{3}}+\beta U(r,\ \varphi-d,\ z)|_{\Gamma_{3}}=g_{3}(x) \quad (x\in\Gamma_{3}).
 \end{array}
\end{equation}
Here $U(r,\ \varphi,\ z)$ is the function $U(x)$ written in the cylindrical coordinates; $\alpha,\ \beta\in{\mathbb R};$
$|\alpha+\beta|<2$.

For $n=2$, we put $K=\{y:\ r>0,\ b_1<\varphi<b_3\},$ $K_t=\{y:\ r>0,\ b_t<\varphi<b_{t+1}\},$
$\gamma_q=\{y:\ r>0,\ \varphi=b_{q}\}$. Write the corresponding nonlocal problem in the plane angle $K:$
\begin{equation}\label{eqPThetaPuasson}
 -\triangle u+u\equiv -\sum\limits_{i=1}^2 u_{y_iy_i}(y)+u(y)=f(y) \quad   (y\in K),
\end{equation}
\begin{equation}\label{eqBThetaPuasson}
 \begin{array}{c}
 u|_{\gamma_{1}}+\alpha u(r,\ \varphi+d)|_{\gamma_{1}}=
  g_{1}(y)\ \quad (y\in\gamma_{1}),\\
 u|_{\gamma_{3}}+\beta u(r,\ \varphi-d)|_{\gamma_{3}}=g_{3}(y) \quad (y\in\gamma_{3}).
 \end{array}
\end{equation}

Clearly, the corresponding homogeneous problem with parameter $\lambda$ has the form
\begin{equation}\label{eqPLambdaPuasson}
 -\tilde U_{\varphi\varphi}+\lambda^2 \tilde U=0 \quad   (\varphi\in(b_1,\ b_3)),
\end{equation}
\begin{equation}\label{eqBLambdaPuasson}
 \begin{array}{c}
 \tilde U(\varphi)|_{\varphi=b_{1}}+\alpha \tilde U(\varphi+d)|_{\varphi=b_{1}}=0,\\
 \tilde U(\varphi)|_{\varphi=b_{3}}+\beta \tilde U(\varphi-d)|_{\varphi=b_{3}}=0.
 \end{array}
\end{equation}

One can easily find the eigenvalues of problem~(\ref{eqPLambdaPuasson}), (\ref{eqBLambdaPuasson}). If $\alpha+\beta=0$,
then we have
$$
   \lambda_k=i \frac{\displaystyle \pi}{\displaystyle b_3-b_1}k\quad (k=\pm1,\ \pm2,\ \dots).
$$
If $0<|\alpha+\beta|<2$, then we have
$$
 \begin{array}{c}
  \lambda_{k}=i \frac{\displaystyle 2\pi}{\displaystyle b_3-b_1}k\quad (k=\pm1,\ \pm2,\ \dots),\\
\\
  \lambda_{p}=\left\{
  \begin{array}{c}
   i \frac{\displaystyle \pm 2\arctan\frac{\sqrt{4-(\alpha+\beta)^2}}{\alpha+\beta}}
     {\displaystyle b_3-b_1}+  i\frac{\displaystyle 4\pi p}{\displaystyle b_3-b_1} \quad \mbox{for } -2<\alpha+\beta<0,\\
  i \frac{\displaystyle 2\pi\pm 2\arctan\frac{\sqrt{4-(\alpha+\beta)^2}}{\alpha+\beta}}
     {\displaystyle b_3-b_1}+  i\frac{\displaystyle 4\pi p}{\displaystyle b_3-b_1} \quad \mbox{for } 0<\alpha+\beta<2
  \end{array}\right.\\
      (p=0,\ \pm1,\ \pm2,\ \dots).
 \end{array}
$$

Obviously, the line $\Im\lambda=0$ contains no eigenvalues of problem~(\ref{eqPLambdaPuasson}),
(\ref{eqBLambdaPuasson}). Therefore by Theorem~\ref{thLThetaFred}, the operator
\begin{equation}\label{eqLThetaPuasson}
\begin{array}{c}
 \bigl(-\triangle u+u,\ u|_{\gamma_1}+\alpha u(r,\ \varphi+d)|_{\gamma_{1}},\
 u|_{\gamma_3}+\beta u(r,\ \varphi-d)|_{\gamma_3}\bigr):E_1^2(K)\to\\
 \to E_1^0(K)\times\prod\limits_{\sigma=1,\,3}E_1^{3/2}(\gamma_\sigma)
 \end{array}
\end{equation}
is Fredholm. Let us show that operator~(\ref{eqLThetaPuasson}) has trivial kernel.

Suppose $u\in E_1^2(K)$ is a solution for homogeneous problem~(\ref{eqPThetaPuasson}), (\ref{eqBThetaPuasson}).
Introduce the difference operator ${\cal R}_K=P_K{\cal R}I_K,$ where
$$
 {\cal R}w(y)=w(r,\ \varphi)-\alpha w(r,\ \varphi+d)-\beta w(r,\ \varphi-d).
$$
Put $u={\cal R}_K w.$ Since $|\alpha+\beta|<2$, the matrix
$
  {\cal R}_1=\left(
 \begin{array}{cc}
  1 & -\alpha\\
  -\beta & 1
 \end{array}
 \right)
$
corresponding to the difference operator ${\cal R}_K$ is non-singular and
$$
   {\cal R}_1^{-1}=\frac{\displaystyle 1}{\displaystyle 1-\alpha\beta}\left(
 \begin{array}{cc}
  1 & \alpha\\
  \beta & 1
 \end{array}
 \right).
$$
Therefore, by Lemma~\ref{lRSpectrum}, the operator ${\cal R}_K$ has the bounded inverse ${\cal R}_K^{-1}$ and
$w={\cal R}_K^{-1}u$.

Now we shall show that $w\in {\cal E}_1^2(K)\cap E_1^1(K)$ and
$w|_{\gamma_1}=w|_{\gamma_3}=0.$ Indeed, by Lemma~\ref{lRWtoW}, $w\in {\cal E}_1^2(K)$. Further, using the isomorphism
${\cal U}$, the matrix ${\cal R}_1^{-1}$, and Lemma~\ref{lRMatrix}, we get
\begin{equation}\label{eqWgamma2}
 \begin{array}{c}
 w_1|_{\gamma_2}=[{\cal U}w]_1(r,\ b_2)=
\frac{\displaystyle 1}{\displaystyle 1-\alpha\beta}([{\cal U}u]_1(r,\ b_2)+\alpha[{\cal U}u]_2(r,\ b_2))=\\
=\frac{\displaystyle 1}{\displaystyle 1-\alpha\beta}(u(r,\ b_2)+\alpha u(r,\ b_3)),\\
\\
 w_2|_{\gamma_2}=[{\cal U}w]_2(r,\ b_1)=
\frac{\displaystyle 1}{\displaystyle 1-\alpha\beta}(\beta[{\cal U}u]_1(r,\ b_1)+[{\cal U}u]_2(r,\ b_1))=\\
=\frac{\displaystyle 1}{\displaystyle 1-\alpha\beta}(\beta u(r,\ b_1)+u(r,\ b_2)).
 \end{array}
\end{equation}
But the function $u$ satisfies homogeneous conditions~(\ref{eqBThetaPuasson}) and therefore
$\alpha u(r,\ b_3)=\beta u(r,\ b_1)$.
Combining this with~(\ref{eqWgamma2}), we see that $w_1|_{\gamma_2}=w_2|_{\gamma_2}$, i.e., $w\in E_1^1(K)$.

Similarly,
$$
 \begin{array}{c}
 w_1|_{\gamma_1}=[{\cal U}w]_1(r,\ b_1)=
 \frac{\displaystyle 1}{\displaystyle 1-\alpha\beta}([{\cal U}u]_1(r,\ b_1)+\alpha[{\cal U}u]_2(r,\ b_1))=\\
 =\frac{\displaystyle 1}{\displaystyle 1-\alpha\beta}(u(r,\ b_1)+\alpha u(r,\ b_2))=0,\\
\\
w_2|_{\gamma_3}=[{\cal U}w]_2(r,\ b_2)=
\frac{\displaystyle 1}{\displaystyle 1-\alpha\beta}(\beta[{\cal U}u]_1(r,\ b_2)+[{\cal U}u]_2(r,\ b_2))=\\
=\frac{\displaystyle 1}{\displaystyle 1-\alpha\beta}(\beta u(r,\ b_2)+u(r,\ b_3))=0,
 \end{array}
$$
since the function $u$ satisfies homogeneous conditions~(\ref{eqBThetaPuasson}).

Therefore from the imbedding ${\cal E}_1^2(K)\cap E_1^1(K)\subset W^1(K)$, it follows that
$w\in\mathaccent23 W^1(K)$ and $w$ is a generalized solution to boundary value
problem~(\ref{eqPRExDD}), (\ref{eqBoundHomogExDD}) for $f=0$.
In Example~\ref{exDD}, it is shown that $w=0$ which implies $u={\cal R}_K w=0.$

In order to prove that range of operator~(\ref{eqLThetaPuasson}) coincides with
$E_1^0(K)\times\prod\limits_{\sigma=1,\,3}E_1^{3/2}(\gamma_\sigma),$ we study the problems that are formally adjoint
to problems~(\ref{eqPPuasson}), (\ref{eqBPuasson}) and~(\ref{eqPThetaPuasson}), (\ref{eqBThetaPuasson}) with respect
to the Green formulas. Similarly to Example~\ref{exGr1}, we obtain the following nonlocal transmission problems:
\begin{equation}\label{eqQPuasson}
 -\triangle V_t+V_t=f(x) \quad (x\in \Omega_t;\ t=1,\ 2)
\end{equation}
\begin{equation}\label{eqCPuasson}
 \begin{array}{c}
 V_1|_{\Gamma_1}=g_1(x) \quad (x\in\Gamma_1),\\
 V_2|_{\Gamma_3}=g_3(x) \quad (x\in\Gamma_3),
 \end{array}
\end{equation}
\begin{equation}\label{eqTPuasson}
 \begin{array}{c}
 V_1|_{\Gamma_2}-V_2|_{\Gamma_2}=h_{21}(x) \quad (x\in\Gamma_2),\\
 \frac{\displaystyle\partial V_1}{\displaystyle\partial n_2}\Bigr|_{\Gamma_2}-
      \frac{\displaystyle\partial V_2}{\displaystyle\partial n_2}\Bigr|_{\Gamma_2}+
 \alpha\frac{\displaystyle\partial V_1}{\displaystyle\partial n_1}(r,\ \varphi-d,\ z)\Bigr|_{\Gamma_2}
 +\beta\frac{\displaystyle\partial  V_2}{\displaystyle\partial n_3}(r,\ \varphi+d,\ z)\Bigr|_{\Gamma_2}=\\
 =h_{22}(x) \quad (x\in\Gamma_2)
 \end{array}
\end{equation}
and
\begin{equation}\label{eqQThetaPuasson}
 -\triangle v_t+v_t=f(y) \quad (y\in K_t;\ t=1,\ 2)
\end{equation}
\begin{equation}\label{eqCThetaPuasson}
 \begin{array}{c}
 v_1|_{\gamma_1}=g_1(y) \quad (y\in\gamma_1),\\
 v_2|_{\gamma_3}=g_3(y) \quad (y\in\gamma_3),
 \end{array}
\end{equation}
\begin{equation}\label{eqTThetaPuasson}
 \begin{array}{c}
 v_1|_{\gamma_2}-v_2|_{\gamma_2}=h_{21}(y) \quad (y\in\gamma_2),\\
 \frac{\displaystyle\partial v_1}{\displaystyle\partial n_2}\Bigr|_{\gamma_2}-
      \frac{\displaystyle\partial v_2}{\displaystyle\partial n_2}\Bigr|_{\gamma_2}+
 \alpha\frac{\displaystyle\partial  v_1}{\displaystyle\partial n_1}(r,\ \varphi-d)\Bigr|_{\gamma_2}
 +\beta\frac{\displaystyle\partial v_2}{\displaystyle\partial n_3}(r,\ \varphi+d)\Bigr|_{\gamma_2}=\\
 =h_{22}(y) \quad (y\in\gamma_2).
 \end{array}
\end{equation}
Here $n_1$ is the unit normal vector to $\Gamma_1$ ($\gamma_1$) direct inside $\Omega_1$ ($K_1$);
$n_2$ and $n_3$ are the unit normal vectors to $\Gamma_2$ ($\gamma_2$)
$\Gamma_3$ ($\gamma_3$) correspondingly directed inside $\Omega_2$ ($K_2$).
As we have notices in the proof of Theorem~\ref{thLThetaFred}, $\lambda_0$ is an eigenvalue of
problem~(\ref{eqPLambdaPuasson}), (\ref{eqBLambdaPuasson}) iff  $\lambda'_0=\bar\lambda_0$ is an eigenvalue
of nonlocal transmission problem with parameter $\lambda$ corresponding to
problem~(\ref{eqQPuasson})--(\ref{eqTPuasson}) (which can be written in the obvious way). Hence this problem also has no
eigenvalues on the line $\Im\lambda=0$. Then by Theorem~\ref{thAprE2'}, the operator
\begin{equation}\label{eqMThetaPuasson}
 \begin{array}{c}
 \bigl(-v_\triangle+v,\ v_1|_{\gamma_1},\ v_2|_{\gamma_3},\  v_1|_{\gamma_2}-v_2|_{\gamma_2},\
  \frac{\displaystyle\partial v_1}{\displaystyle\partial n_2}\Bigr|_{\gamma_2}-
      \frac{\displaystyle\partial v_2}{\displaystyle\partial n_2}\Bigr|_{\gamma_2}+\\
 +\alpha\frac{\displaystyle\partial v_1}{\displaystyle\partial n_1}(r,\ \varphi-d)\Bigr|_{\gamma_2}
 +\beta\frac{\displaystyle\partial v_2}{\displaystyle\partial n_3}(r,\ \varphi+d)\Bigr|_
  {\gamma_2}\bigr):\\
 {\cal E}_1^2(K)\to{\cal E}_1^0(K)\times\prod\limits_{\sigma=1,\,3}E_1^{3/2}(\gamma_\sigma)
  \times\prod\limits_{\nu=1}^2 E_1^{2-\nu+1/2}(\gamma_2)
 \end{array}
\end{equation}
has finite dimensional kernel. Here $v_\triangle(y)=\triangle v_t(y)$ for $y\in K_t$, $t=1,\ 2$.
Let us show that kernel of operator~(\ref{eqMThetaPuasson}) is trivial.

Suppose $v\in {\cal E}_1^2(K)$ is a solution for homogeneous problem~(\ref{eqQThetaPuasson})--(\ref{eqTThetaPuasson}).
Consider the adjoint difference operator ${\cal R}_K^*$. The matrix
$
  {\cal R}_1^*=\left(
 \begin{array}{cc}
  1 & -\beta\\
  -\alpha & 1
 \end{array}
 \right)
$
corresponds to the difference operator ${\cal R}_K^*$. Since $|\alpha+\beta|<2$, the matrix ${\cal R}_1^*$ is non-singular and
by Lemma~\ref{lRSpectrum}, there exists the inverse operator $({\cal R}_K^*)^{-1}.$ Put $v=({\cal R}_K^*)^{-1}w$;
hence $w={\cal R}_K^*v$.

Let us show that $w\in E_1^2(K)$ and
$w|_{\gamma_{1}}+\beta w(r,\ \varphi+d)|_{\gamma_{1}}=0$,
$w|_{\gamma_{3}}+\alpha w(r,\ \varphi-d)|_{\gamma_{3}}=0.$
Indeed, by Lemma~\ref{lRWtoW}, $w\in {\cal E}_1^2(K)$. Further using the isomorphism
${\cal U}$, the matrix ${\cal R}_1^*$, and Lemma~\ref{lRMatrix}, we get
\begin{equation}\label{eqWgamma2'}
 \begin{array}{c}
 w_1|_{\gamma_2}=[{\cal U}w]_1(r,\ b_2)=
[{\cal U}v]_1(r,\ b_2)-\beta[{\cal U}v]_2(r,\ b_2)=\\
=v_1(r,\ b_2)-\beta v_2(r,\ b_3)=
v_1(r,\ b_2),\\
\\
 w_2|_{\gamma_2}=[{\cal U}w]_2(r,\ b_1)=
-\alpha[{\cal U}v]_1(r,\ b_1)+[{\cal U}v]_2(r,\ b_1)=\\
=-\alpha v_1(r,\ b_1)+ v_2(r,\ b_2)=v_2(r,\ b_2),
 \end{array}
\end{equation}
since $v$ satisfies homogeneous conditions~(\ref{eqCThetaPuasson}). From~(\ref{eqWgamma2'}) and
homogeneous conditions~(\ref{eqTThetaPuasson}), we get $w_1|_{\gamma_2}=w_2|_{\gamma_2}$.

Similarly,
\begin{equation}\label{eqPartialWgamma2'}
 \begin{array}{c}
\frac{\displaystyle \partial w_1}{\displaystyle \partial \varphi}\Big|_{\gamma_2}
=\frac{\displaystyle \partial v_1}{\displaystyle \partial\varphi} (r,\ b_2)-
\beta \frac{\displaystyle \partial v_2}{\displaystyle \partial\varphi}(r,\ b_3),
\\
 \frac{\displaystyle \partial w_2}{\displaystyle \partial \varphi}\Big|_{\gamma_2}
=-\alpha \frac{\displaystyle \partial v_1}{\displaystyle \partial \varphi}(r,\ b_1)+
\frac{\displaystyle \partial v_2}{\displaystyle \partial \varphi} (r,\ b_2).
 \end{array}
\end{equation}
Taking into account that $\frac{\displaystyle \partial}{\displaystyle\partial n_i}=
\frac{\displaystyle 1}{\displaystyle r}\frac{\displaystyle \partial}{\displaystyle \partial \varphi}$ ($i=1,\ 2$) and
$\frac{\displaystyle \partial}{\displaystyle\partial n_3}=
-\frac{\displaystyle 1}{\displaystyle r}\frac{\displaystyle \partial}{\displaystyle \partial \varphi}$,
from~(\ref{eqPartialWgamma2'}) and homogeneous conditions~(\ref{eqTThetaPuasson}), we obtain
$\frac{\displaystyle \partial w_1}{\displaystyle \partial n_2}\Big|_{\gamma_2}=
\frac{\displaystyle \partial w_2}{\displaystyle \partial n_2}\Big|_{\gamma_2}$. Therefore, $w\in E_1^2(K)$.
Analogously one can show that $w|_{\gamma_{1}}+\beta w(r,\ \varphi+d)|_{\gamma_{1}}=0$,
$w|_{\gamma_{3}}+\alpha w(r,\ \varphi-d)|_{\gamma_{3}}=0.$

This means that $w\in E_1^2(K)$ is a solution for the problem
\begin{equation}\label{eqPThetaPuasson'}
 -\triangle w+w=0 \quad   (y\in K),
\end{equation}
\begin{equation}\label{eqBThetaPuasson'}
 \begin{array}{c}
 w|_{\gamma_{1}}+\beta w(r,\ \varphi+d)|_{\gamma_{1}}=
  0\ \quad (y\in\gamma_{1}),\\
 w|_{\gamma_{3}}+\alpha w(r,\ \varphi-d)|_{\gamma_{3}}=0 \quad (y\in\gamma_{3}).
 \end{array}
\end{equation}
But problem~(\ref{eqPThetaPuasson'}), (\ref{eqBThetaPuasson'}) is a nonlocal boundary value problem of
type~(\ref{eqPThetaPuasson}), (\ref{eqBThetaPuasson}) (one must replace $\alpha$ by $\beta$ and $\beta$ by $\alpha$).
Hence, by the above, $w=0$ if $|\alpha+\beta|<2$. This implies $v=\bar{\cal R}_K w=0.$

From Lemma~\ref{lKerLThetaAdj}, it follows that dimension of cokernal of operator~(\ref{eqLThetaPuasson}) is equal to dimension
of kernel of operator~(\ref{eqMThetaPuasson}). Therefore cokernal of operator~(\ref{eqLThetaPuasson}) is trivial. Finally applying Theorem~\ref{thLSolv}, we obtain that

{\it nonlocal boundary value problem~(\ref{eqPPuasson}), (\ref{eqBPuasson}) has a unique solution
$U\in H_{1+l}^{l+2}(\Omega)$ for every right-hand side $\{f,\ g_1,\ g_3\}\in
H_{1+l}^l(K)\times\prod\limits_{\sigma=1,\,3}H_{1+l}^{l+3/2}(\Gamma_\sigma)$.}

\appendix
\section{A priori estimates for the operator $L^*$ in ${\mathbb R}^n$}\label{appendL*Rn}

\subsection{Some approaches for ordinary differential equations.}
Let ${\cal P}(\xi',\ -i\frac{\displaystyle d}{\displaystyle dx_n})$ and
$B_\nu(\xi',\ -i\frac{\displaystyle d}{\displaystyle dx_n})$ $(\nu=1,\ \dots,\ J;\ J\ge 1)$ be differential operators with
constant coefficients and parameter $\xi'=(\xi_1,\ \dots,\ \xi_{n-1})\in{\mathbb R}^{n-1}$ such that after replacing
$-i\frac{\displaystyle d}{\displaystyle dx_n}$ by $\xi_n$, we get polynomials of orders $2m$ and $m_\nu\le 2m-1$
that are homogeneous with respect to $(\xi',\ \xi_n)$ correspondingly.

Let the following condition hold.
\begin{condition}\label{condEllipPXi}
  ${\cal P}(\xi',\ \xi_n)\ne 0$ for all $(\xi',\ \xi_n)\ne 0.$
\end{condition}

Consider the bounded operator $L_{\xi'}: W^{2m}({\mathbb R})\to L_2({\mathbb R})\times{\mathbb C}^J$ given by
$$
 L_{\xi'}u=({\cal P}(\xi',\ -i\frac{\displaystyle d}{\displaystyle dx_n})u,\ 
                                      B_1(\xi',\ -i\frac{\displaystyle d}{\displaystyle dx_n})u|_{x_n=0},\ \dots,\ 
                                      B_J(\xi',\ -i\frac{\displaystyle d}{\displaystyle dx_n})u|_{x_n=0}).
$$
Introduce the adjoint operator $L_{\xi'}^*: L_2({\mathbb R})\times{\mathbb C}^J\to W^{-2m}({\mathbb R})$ that takes 
$\Psi=(\psi,\ d_1,\ \dots,\ d_J)\in  L_2({\mathbb R})\times{\mathbb C}^J$ to $L_{\xi'}^*\Psi$ by the rule
$$ 
 <u,\ L_{\xi'}^*\Psi>=
 <{\cal P}(\xi',\ -i\frac{\displaystyle d}{\displaystyle dx_n})u,\ \psi>+\sum\limits_{\nu=1}^J
 <B_\nu(\xi',\ -i\frac{\displaystyle d}{\displaystyle dx_n})u|_{x_n=0},\ d_\nu>
$$
for all $u\in W^{2m}({\mathbb R})$.

\begin{lemma}\label{lLXiFred}
 Suppose $n\ge 2$; then for all $\xi'\in{\mathbb R}^{n-1},$
 $\xi'\ne 0,$ the operator $L_{\xi'}$ is Fredholm, its kernel is trivial.
\end{lemma}
\begin{proof}
 Since $n\ge 2,$ condition~\ref{condEllipPXi} implies that
\begin{equation}\label{eqLXiFred1}
 k_1(1+|\xi_n|^2)^{2m}\le |{\cal P}(\xi',\ \xi_n)|^2\le k_2(1+|\xi_n|^2)^{2m}\ \mbox{for } \xi'\ne 0.
\end{equation} 
Here $k_1,\ k_2$ depend on $\xi'$ and do not depend on $\xi_n.$
Multiplying the first inequality in~(\ref{eqLXiFred1}) by $|\tilde u(\xi_n)|^2$ ($\tilde u$ is the Fourier transform of the
function $u$ with respect to $x_n$) and integrating over ${\mathbb R},$ we obtain
$$
 \|u\|_{W^{2m}({\mathbb R})}\le k_3\|{\cal P}(\xi',\ \xi_n)u\|_{L_2({\mathbb R})},
$$
where $k_3>0$ depend only on $\xi'$ and do not depend on $u.$ The last inequality implies that
the operator $L_{\xi'}$ has trivial kernel and closed range.

Let us show that cokernel of the operator $L_{\xi'}$ is of finite dimension. 
Using the Fourier transform and inequality~(\ref{eqLXiFred1}), one can easily check that for
$n\ge 2,\ \xi'\ne 0$, the operator ${\cal P}(\xi',\ -i\frac{\displaystyle d}{\displaystyle dx_n})$ maps
$W^{2m}({\mathbb R})$ onto $L_2({\mathbb R}).$ Therefore the operator $L_{\xi'}$ maps
$W^{2m}({\mathbb R})$ onto $L_2({\mathbb R})\times {\mathbb M}^J,$ where ${\mathbb M}^J$ is a 
closed (since range of $L_{\xi'}$ is closed) subspace of ${\mathbb C}^J.$ But ${\mathbb C}^J$ is a finite dimensional space; 
hence cokernel of the operator $L_{\xi'}$ is also finite dimensional.
\end{proof}

\begin{lemma}\label{lLXiAdjApr}
 Suppose $n\ge 2$; then for all $\xi'\in{\mathbb R}^{n-1},$ $\xi'\ne 0,$ we have

 I) the operator $L_{\xi'}^*$ is Fredholm, its range coincides with $W^{-2m}({\mathbb R});$

 II) for all $\Psi=(\psi,\ d_1,\ \dots,\ d_J)\in  L_2({\mathbb R})\times{\mathbb C}^J$, the following estimate holds:
 \begin{equation}\label{eqLXiAdjApr}
   \|\psi\|_{L_2({\mathbb R})}\le c_{\xi'}\bigl(\|L_{\xi'}^*\Psi\|_{W^{-2m}({\mathbb R})}+\sum\limits_{\nu=1}^J |d_\nu|\bigr),
 \end{equation}
 where $c_{\xi'}>0$ depends on $\xi'$ and does not depend on $\Psi$;

 III) if $\xi'\in {\mathbb K}\subset{\mathbb R}^{n-1},$ where ${\mathbb K}$ is a compactum 
 such that ${\mathbb K}\cap\{0\}=\varnothing$, then
 inequality~(\ref{eqLXiAdjApr}) holds with a constant $c$ that does not depend on $\xi'.$
\end{lemma}
\begin{proof}
 I) follows from Lemma~\ref{lLXiFred}. Let us prove II). Denote by $\ker(L_{\xi'}^*)$ kernel of the operator
 $L_{\xi'}^*.$ Since $L_{\xi'}^*$ is Fredholm, $\ker(L_{\xi'}^*)$ is of finite dimension.

Let us show that in the space $\ker(L_{\xi'}^*)$, we can introduce the norm
$$
 \|\hat \Psi\|_{\ker(L_{\xi'}^*)}=\left(\sum\limits_{\nu=1}^J|\hat d_\nu|^2\right)^{1/2},\ 
 \hat\Psi=(\hat\psi,\ \hat d_1,\ \dots,\ \hat d_J)\in\ker(L_{\xi'}^*)\subset L_2({\mathbb R})\times{\mathbb C}^J,
$$
 which is equivalent to the standart norm in $L_2({\mathbb R})\times{\mathbb C}^J$.
 Among all of the properties of a norm, the following one is not obvious: $\hat\Psi=0$ whenever
$\|\hat \Psi\|_{\ker(L_{\xi'}^*)}=0.$ Check it. Suppose $\|\hat \Psi\|_{\ker(L_{\xi'}^*)}=0$; then
$\hat\Psi=(\hat\psi,\ 0,\ \dots,\ 0).$ Since $\hat\Psi\in\ker(L_{\xi'}^*),$
it follows from the definition of the operator $L_{\xi'}^*$ that
\begin{equation}\label{eqLXiAdjFred1}
 <{\cal P}(\xi',\ -i\frac{\displaystyle d}{\displaystyle dx_n})u,\ \hat\psi>=0
\end{equation}
for all $u\in W^{2m}({\mathbb R}).$

As we have already mentioned in the proof of Lemma~\ref{lLXiFred},  
the operator ${\cal P}(\xi',\ -i\frac{\displaystyle d}{\displaystyle dx_n})$ maps
$W^{2m}({\mathbb R})$ onto $L_2({\mathbb R})$ if $n\ge 2,\ \xi'\ne 0$. From this and~(\ref{eqLXiAdjFred1}), it follows that
$\hat\psi=0$; hence, $\hat\Psi=0.$ 

Now we get that the norm $\|\cdot\|_{\ker(L_{\xi'}^*)}$ is equivalent to the norm $\|\cdot\|_{L_2({\mathbb R})\times{\mathbb C}^J}$, 
since the space $\ker(L_{\xi'}^*)$ is of finite dimension.

The operator $L_{\xi'}^*$ is closed and range of $L_{\xi'}^*$ is closed; hence from theorem~2.3 \cite[\S2]{Kr},
it follows that for any
$\Psi=(\psi,\ d_1,\ \dots,\ d_J)\in L_2({\mathbb R})\times{\mathbb C}^J$, there exists a $\Phi\in L_2({\mathbb R})\times{\mathbb C}^J$ such that
$L_{\xi'}^*\Psi=L_{\xi'}^*\Phi$ and
$$
 \|\Phi\|_{L_2({\mathbb R})\times{\mathbb C}^J}\le k_1\|L_{\xi'}^*\Psi\|_{W^{-2m}({\mathbb R})},
$$
where $k_1>0$ depends on $\xi'$ and does not depend on $\Phi$ and $\Psi.$
But $\Psi=\Phi+\hat\Psi,$ where $\hat\Psi=(\hat\psi,\ \hat d_1,\ \dots,\ \hat d_J)\in\ker(L_{\xi'}^*)$; therefore,
$$
 \|\Psi\|_{L_2({\mathbb R})\times{\mathbb C}^J}\le k_1\|L_{\xi'}^*\Psi\|_{W^{-2m}({\mathbb R})}+
 \|\hat\Psi\|_{L_2({\mathbb R})\times{\mathbb C}^J}.
$$
By proved, the norms $\|\cdot\|_{\ker(L_{\xi'}^*)}$ and $\|\cdot\|_{L_2({\mathbb R})\times{\mathbb C}^J}$ are equivalent; this implies
$$
 \begin{array}{c}
 \|\psi\|_{L_2({\mathbb R})}\le  \|\Psi\|_{L_2({\mathbb R})\times{\mathbb C}^J}\le  
 k_1\|L_{\xi'}^*\Psi\|_{W^{-2m}({\mathbb R})}+k_2\sum\limits_{\nu=1}^J|\hat d_\nu|\le\\
 \le k_1\|L_{\xi'}^*\Psi\|_{W^{-2m}({\mathbb R})}+k_2\sum\limits_{\nu=1}^J|d_\nu|+
 k_2 \|\Phi\|_{L_2({\mathbb R})\times{\mathbb C}^J}\le \\
 \le k_1\|L_{\xi'}^*\Psi\|_{W^{-2m}({\mathbb R})}+
 k_2\sum\limits_{\nu=1}^J|d_\nu|+k_1k_2\|L_{\xi'}^*\Psi\|_{W^{-2m}({\mathbb R})}\le\\
 \le c_{\xi'}(\|L_{\xi'}^*\Psi\|_{W^{-2m}({\mathbb R})}+\sum\limits_{\nu=1}^J|d_\nu|),
 \end{array}
$$
where $c_{\xi'}=\max(k_1+k_1k_2,\ k_2)$.

Let us prove III). Suppose III) does not hold; then there exist sequences
$\{(\xi')^k\}\subset {\mathbb K},$ $\{\Psi_k\}=\{(\psi_k,\ d_1^k,\ \dots,\ d_J^k)\},$ $k=1,\ 2,\ \dots,$ such that
$\|\psi_k\|_{L_2({\mathbb R})}=1,$
\begin{equation}\label{eqLXiAdjFred2}
 \|L_{(\xi')^k}^*\Psi_k\|_{W^{-2m}({\mathbb R})}+\sum\limits_{\nu=1}^J|d_\nu^k|\to 0\ \mbox{as } k\to\infty.
\end{equation}
Choose from $\{(\xi')^k\}$ a subsequence (we shall denote it $\{(\xi')^k\}$ too) that converges to a $(\xi')^0\in {\mathbb K}.$ 
By assumption, $(\xi')^0\ne 0$; therefore by proved, estimate~(\ref{eqLXiAdjApr}) holds for $\xi'=(\xi')^0$.

Notice that
$$
 \begin{array}{c}
  \|L_{(\xi')^0}^*\Psi_k\|_{W^{-2m}({\mathbb R})}\le  \|L_{(\xi')^k}^*\Psi_k\|_{W^{-2m}({\mathbb R})}+\\
  +\|L_{(\xi')^k}^*-L_{(\xi')^0}^*\|_{ L_2({\mathbb R})\times{\mathbb C}^J\to W^{-2m}({\mathbb R})}\cdot
   \|\Psi_k\|_{L_2({\mathbb R})\times{\mathbb C}^J}.
 \end{array}
$$
From~(\ref{eqLXiAdjFred2}), it follows that $\|L_{(\xi')^k}^*\Psi_k\|_{W^{-2m}({\mathbb R})}\to 0$. Further, 
$\|L_{(\xi')^k}^*-L_{(\xi')^0}^*\|_{ L_2({\mathbb R})\times{\mathbb C}^J\to W^{-2m}({\mathbb R})}\to 0$, since $L_{\xi'}$ 
depends on $\xi'$ polynomially. Finally, $\|\Psi_k\|_{L_2({\mathbb R})\times{\mathbb C}^J}$ is uniformly bounded by a constant not
depending on $k$ which follows from~(\ref{eqLXiAdjFred2}) and relation $\|\psi_k\|_{L_2({\mathbb R})}=1$. 
Hence, $ \|L_{(\xi')^0}^*\Psi_k\|_{W^{-2m}({\mathbb R})}\to 0$ as $k\to\infty$. Combining this with~(\ref{eqLXiAdjFred2}), we obtain
\begin{equation}\label{eqLXiAdjFred3}
 \|L_{(\xi')^0}^*\Psi_k\|_{W^{-2m}({\mathbb R})}+\sum\limits_{\nu=1}^J|d_\nu^k|\to 0\ \mbox{as } k\to\infty.
\end{equation}
Now applying estimate~(\ref{eqLXiAdjApr}) to the sequence $\{\Psi_k\}$ and $\xi'=(\xi')^0$, from~(\ref{eqLXiAdjFred3}),
we eventually get
$$
 \|\psi_k\|_{L_2({\mathbb R})}\to 0\ \mbox{as } k\to\infty.
$$
This contradicts the assumption $\|\psi_k\|_{L_2({\mathbb R})}=1.$ 
\end{proof}

\subsection{A priori estimates in ${\mathbb R}^n.$} 
Write a point $x\in{\mathbb R}^n$ ($n\ge2$) in the form
$x=(x',\ x_n),$ where $x'=(x_1,\ \dots,\ x_{n-1})\in{\mathbb R}^{n-1},$ $x_n\in{\mathbb R}.$ Similarly, write a point 
$\xi\in{\mathbb R}^n$ ($n\ge2$) in the form
$\xi=(\xi',\ \xi_n),$ where $\xi'=(\xi_1,\ \dots,\ \xi_{n-1})\in{\mathbb R}^{n-1},$ $\xi_n\in{\mathbb R}.$

Let ${\cal P}(D)={\cal P}(D_{x'},\ D_{x_n}),$  $B_\nu(D)=B_\nu(D_{x'},\ D_{x_n})$ 
$(\nu=1,\ \dots,\ J;\ J\ge 1)$ be differential operators with constant coefficients
such that after replacing $D=(D_{x'},\ D_{x_n})$ by $\xi=(\xi',\ \xi_n)$, we get polynomials 
${\cal P}(\xi)={\cal P}(\xi',\ \xi_n),$ $B_\nu(\xi)=B_\nu(\xi',\ \xi_n)$ of orders $2m$ and $m_\nu\le 2m-1$ correspondingly that are
homogeneous with respect to $\xi=(\xi',\ \xi_n)$.
We shall suppose that the polynomial ${\cal P}(\xi)$ satisfies condition~\ref{condEllipPXi}.

Consider the bounded operator
$$
 L: W^{2m}({\mathbb R}^n)\to L_2({\mathbb R}^n)\times\prod\limits_{\nu=1}^J W^{2m-m_\nu-1/2}({\mathbb R}^{n-1})
$$
given by
$$
 LU=({\cal P}(D)U,\ B_1(D)U|_{x_n=0},\ \dots,\ B_J(D)U|_{x_n=0}).
$$
Notice that the problem corresponding to the operator $L$ is quite artificial. This is not a boundary value problem,
since a solution $U$ is considered in ${\mathbb R}^n$. And this is not a transmission problem,
since we impose the trace conditions on the hyperplane $\{x_n=0\}$, but not transmission conditions. Moreover the
operators $B_J(D)$ do not cover the operator ${\cal P}(D)$ on the hyperplane $\{x_n=0\}$.
Nevertheless we need this problem for getting a priori estimates of solutions to adjoint nonlocal 
problems (\S\ref{sectLMAdj}). This is explained by the specific character of our method, which may be called
``separation of nonlocality".

Introduce the adjoint operator
$L^*: L_2({\mathbb R}^n)\times\prod\limits_{\nu=1}^J W^{-2m+m_\nu+1/2}({\mathbb R}^{n-1})\to W^{-2m}({\mathbb R}^n)$.
The operator $L^*$ takes
$F=(f_0,\ g_1,\ \dots,\ g_J)\in  L_2({\mathbb R}^n)\times\prod\limits_{\nu=1}^J W^{-2m+m_\nu+1/2}({\mathbb R}^{n-1})$ to
$L^*F$ by the rule
$$ 
 \begin{array}{c}
 <U,\ L^*F>= <{\cal P}(D)U,\ f_0>+\sum\limits_{\nu=1}^J
 <B_\nu(D)U|_{x_n=0},\ g_\nu>\ \\
  \mbox{for all } U\in W^{2m}({\mathbb R}^n).
 \end{array}
$$

Denote ${\mathbb R}_+^n=\{x\in{\mathbb R}^n:\ x_n>0\},$ ${\mathbb R}_-^n=\{x\in{\mathbb R}^n:\ x_n<0\}.$
Consider the space ${\cal W}^l({\mathbb R}^n)=W^l({\mathbb R}_+^n)\oplus W^l({\mathbb R}_-^n)$ with the norm
$
 \|U\|_{{\cal W}^l({\mathbb R}^n)}=\left(\|U_+\|_{W^l({\mathbb R}_+^n)}^2+\|U_-\|_{W^l({\mathbb R}_-^n)}^2\right)^{1/2}.
$

\begin{theorem}\label{thLAdjAprRn}
 Suppose
 $$
 F=(f_0,\ g_1,\ \dots,\ g_J)\in  L_2({\mathbb R}^n)\times\prod\limits_{\nu=1}^J W^{-2m+l+m_\nu+1/2}({\mathbb R}^{n-1}),
 $$
 $$
  L^*F\in\left\{
  \begin{array}{l}
    W^{-2m+l}({\mathbb R}^n)\ \mbox{for } l<2m,\\
    {\cal W}^{-2m+l}({\mathbb R}^n)\ \mbox{for } l\ge2m;
  \end{array}
  \right.
 $$
 then $f_0\in{\cal W}^l({\mathbb R}^n)$ and
 \begin{equation}\label{eqLAdjAprRn}
  \|f_0\|_{{\cal W}^l({\mathbb R}^n)}\le c_l \bigl(\|L^*F\|_{-2m+l}+\|f_0\|_{W^{-1}({\mathbb R}^n)}+
  \sum\limits_{\nu=1}^J\|g_\nu\|_{W^{-2m+l+m_\nu+1/2}({\mathbb R}^{n-1})}\bigr),
 \end{equation}
 where 
 $
  \|\cdot\|_{-2m+l}=\left\{
  \begin{array}{l}
    \|\cdot\|_{W^{-2m+l}({\mathbb R}^n)}\ \mbox{for } l<2m,\\
    \|\cdot\|_{{\cal W}^{-2m+l}({\mathbb R}^n)}\ \mbox{for } l\ge2m,
  \end{array}
  \right.
 $
 $c_l>0$ depends on $l\ge0$ and does not depend on $F.$
\end{theorem}
\begin{proof}
Suppose $l=0$. Then using Fourier transform of the functions
$f_0,$ $g_\nu$ and $L^*F$ with respect to $x'$ we derive estimate~(\ref{eqLAdjAprRn}) from Lemma~\ref{lLXiAdjApr} 
(in the same way as estimate~(4.27) \cite[Chapter 2, \S4.4]{LM} follows from~(4.18) \cite[Chapter 2, \S 4.2]{LM}, see the proof of
theorem~4.1 \cite[Chapter 2, \S 4.4]{LM}).

If $l\ge 1$, then we prove that $f_0\in{\cal W}^l({\mathbb R}^n)$ and obtain estimate~(\ref{eqLAdjAprRn}) 
using~(\ref{eqLAdjAprRn}) for $l=0$, the finite difference method, and condition~\ref{condEllipPXi} (in the same way as
estimate~(4.40) \cite[Chapter 2, \S 4.5]{LM} is derived from~(4.40$'$) \cite[Chapter 2, \S 4.5]{LM}, see the proof
of theorem~4.3 \cite[Chapter 2, \S 4.5]{LM}).
\end{proof}

\begin{remark}
Unlike model problems in ${\mathbb R}^n$ (see \cite[Chapter 2, \S3]{LM}), our operator
$L^*$ contains distributions with support on the hyperplane $\{x_n=0\}.$ That is why smoothness of the function $f_0$ 
can be violated on the hyperplane $\{x_n=0\}$ even if
$L^*F$ is infinitely smooth in ${\mathbb R}^{n}$. Moreover, Theorem~\ref{thLAdjAprRn} shows that if we want
the function $f_0$ to be more smooth in ${\mathbb R}_+^{n}$ and ${\mathbb R}_-^{n}$, then we must consider
more smooth function $L^*F$ and {\bf more smooth distributions $g_\nu$ as well}.
\end{remark}
\section{Some properties of weighted spaces}\label{appendWeightSpace}
Introduce the space $H_a^l(\Omega)$ as a completion of the set $C_0^\infty(\bar\Omega\backslash M)$ in the norm
$$
 \|U\|_{H_a^l(\Omega)}=\left(
    \sum_{|\alpha|\le l}\int\limits_\Omega r^{2(a+|\alpha|-l)} |D_x^\alpha U(x)|^2 dx
                                       \right)^{1/2},
$$
where $\Omega=\{x=(y,\ z):\ r>0,\ 0<b_1<\varphi<b_2<2\pi,\ z\in{\mathbb R}^{n-2}\},$ 
$M=\{x=(y,\ z):\ y=0,\ z\in{\mathbb R}^{n-2}\}$. Denote by $H_a^{l-1/2}(\Gamma)$ ($l\ge1$) the space of traces
on the $(n-1)$-dimensional half-plane
$\Gamma=\{x=(y,\ z):\ r>0,\ \varphi=b,\ z\in{\mathbb R}^{n-2}\}$ ($b_1\le b\le b_2$) with the norm
$$
 \|\Psi\|_{H_a^{l-1/2}(\Gamma)}=\inf\|U\|_{H_a^l(\Omega)} \quad
  (U\in H_a^l(\Omega):\  U|_{\Gamma} = \Psi).
$$

Introduce the space $E_a^l(K)$ as a completion of the set $C_0^\infty(\bar K\backslash \{0\})$ in the norm
$$
 \|u\|_{E_a^l(K)}=\left(
    \sum_{|\alpha|\le l}\int\limits_K r^{2a}(r^{2(|\alpha|-l)}+1) |D_y^\alpha u(y)|^2 dy \right)^{1/2},
$$
where $K=\{y\in{\mathbb R}^2:\ r>0,\ 0<b_1<\varphi<b_2<2\pi\}.$ 
By $E_a^{l-1/2}(\gamma)$ ($l\ge1$) we denote the space of traces on the ray 
$\gamma=\{y:\ r>0,\ \varphi=b\}$ ($b_1\le b\le b_2$) with the norm
$$
 \|\psi\|_{E_a^{l-1/2}(\gamma)}=\inf\|u\|_{E_a^l(K)} \quad (u\in E_a^l(K):\  u|_{\gamma} = \psi).
$$

Our aim is to prove the following two theorems. 
\begin{theorem}\label{thTraceH}
 For all $\Psi\in H_a^{l-1/2}(\Gamma)$, we have
 $$
  \left(\int\limits_{\Gamma}r^{2(a-(l-1/2))}|\Psi|^2\,d\Gamma\right)^{1/2} 
    \le c\|\Psi\|_{H_a^{l-1/2}(\Gamma)},
 $$
 where $c>0$ is independent of $\Psi.$
\end{theorem}

\begin{theorem}\label{thTraceE}
 For all $\psi\in E_a^{l-1/2}(\gamma)$, we have
 $$
  \left(\int\limits_{\gamma}r^{2(a-(l-1/2))}|\psi|^2\,d\gamma\right)^{1/2} 
     \le c\|\psi\|_{E_a^{l-1/2}(\gamma)},
 $$
 where $c>0$ is independent of $\psi.$
\end{theorem}

At first, let us formulate two lemmas (see~\cite[Chapter 6, \S 1.3]{NP}.
\begin{lemma}\label{lHEquv}
 The norm $\|U\|_{H_a^l(\Omega)}$ is equivalent to the norm
 $$ 
  \left(\int\limits_{{\mathbb R}^{n-2}}|\eta|^{2(l-a)-2}\|W(\cdot,\ \eta)\|_{E_a^l(K)}^2\,d\eta\right)^{1/2},
 $$
 where $W(y,\ \eta)=\hat U(|\eta|^{-1}y,\ \eta),$ $\hat U(y,\ \eta)$ is the Fourier transform of $U(y,\ z)$ 
 with respect to $z.$
\end{lemma}

\begin{lemma}\label{lEEquv}
 The norm $\|u\|_{E_a^l(\Omega)}$ is equivalent to the norm
 $$ 
  \left(\sum\limits_{k=0}^l\int\limits_{0}^\infty r^{2(a-(l-1/2))}
  \sum_{j=0}^{l-k}(1+r)^{2(l-k-j)}\|(rD_r)^k u(r,\ \cdot)\|_{W^j(b_1,\ b_2)}^2 \,dr\right)^{1/2},
 $$
 $u(r,\ \varphi)$ is the function $u(y)$ written in the polar coordinates.
\end{lemma}

Let us prove Theorem~\ref{thTraceE}. Take a function $u\in E_a^l(K)$ such that
$u|_{\gamma}=\psi,$ $\|u\|_{E_a^l(K)}\le 2\|\psi\|_{E_a^{l-1/2}(\gamma)}.$ 
Since $u(r,\ \varphi)|_{\varphi=b}=\psi(r)$ and the trace operator in Sobolev spaces is bounded, we have
$|\psi(r)|^2\le k_1\|u(r,\ \cdot)\|^2_{W^l(b_1,\ b_2)}.$ Therefore by Lemma~\ref{lEEquv}, we get
\begin{equation}\label{thTraceE1}
 \int\limits_{\gamma}r^{2(a-(l-1/2))}|\psi|^2\,d\gamma \le 
 k_1\int\limits_{0}^\infty r^{2(a-(l-1/2))}\|u(r,\ \cdot)\|^2_{W^l(b_1,\ b_2)}\,dr\le
 k_2\|u\|_{E_a^l(K)}^2.
\end{equation}
Now Theorem~\ref{thTraceE} follows from~(\ref{thTraceE1}) and the 
inequality $\|u\|_{E_a^l(K)}\le 2\|\psi\|_{E_a^{l-1/2}(\gamma)}.$  

\bigskip

Let us prove Theorem~\ref{thTraceH}. Take a function $U\in H_a^l(\Omega)$ such that
$U|_{\Gamma}=\Psi,$ $\|U\|_{H_a^l(\Omega)}\le 2\|\Psi\|_{H_a^{l-1/2}(\Gamma)}.$ 
Using the Fourier transform with respect to $z$ and the Parseval equality, we have
$$
   \int\limits_{\Gamma}r^{2(a-(l-1/2))}|\Psi|^2\,d\Gamma=
 \int\limits_{{\mathbb R}^{n-2}}\int\limits_{{\mathbb R}^1} r^{2(a-(l-1/2))}|\hat\Psi(r,\ \eta)|^2\,drd\eta,
$$
where $\hat\Psi(r,\ \eta)$ is the Fourier transformation of the function $\Psi(r,\ z)$ with respect to $z.$ 
Doing change of variables $r=|\eta|^{-1}r'$ in the last integral and using~(\ref{thTraceE1}), we obtain
\begin{equation}\label{thTraceH1}
 \begin{array}{c}
  \int\limits_{\Gamma}r^{2(a-(l-1/2))}|\Psi|^2\,d\Gamma=\\
  \le \int\limits_{{\mathbb R}^{n-2}}\int\limits_{{\mathbb R}^1}|\eta|^{-2(a-(l-1/2))-1} (r')^{2(a-(l-1/2))}
  |\hat\Psi(\eta^{-1}r',\ \eta)|^2\,dr'd\eta\le \\
 \le k_2 \int\limits_{{\mathbb R}^{n-2}}|\eta|^{2(l-a)-2}\|W(\cdot,\ \eta)\|^2_{E_a^l(K)}\,d\eta,
 \end{array}
\end{equation}
where $W(y,\ \eta)=\hat U(|\eta|^{-1}y,\ \eta).$ Now Theorem~\ref{thTraceH} follows from~(\ref{thTraceH1}), 
Lemma~\ref{lHEquv}, and the inequality $\|U\|_{H_a^l(\Omega)}\le 2\|\Psi\|_{H_a^{l-1/2}(\Gamma)}.$

The author is grateful to professor A.L. Skubachevskii for constant attention to this work.

\end{document}